\documentclass[11pt]{article}
\usepackage{amsfonts}
\usepackage{mathrsfs}
\usepackage{amsthm}
\usepackage{amsmath, cases}
\usepackage{amssymb}
\usepackage{latexsym}
\usepackage{graphicx}
\usepackage{color,varioref}
\usepackage{longtable}
\usepackage{rotating,multirow,array}
\usepackage{enumerate}
\usepackage{caption}
\usepackage{float}
\usepackage{subfigure}
\usepackage{epstopdf}
\usepackage{color}
\usepackage[switch]{lineno}
\usepackage{algorithm}
\usepackage{algorithmic}
\usepackage{multirow}
\definecolor{red}{rgb}{0.9,0,0}
\definecolor{green}{rgb}{0,0.9,0}
\definecolor{blue}{rgb}{0,0,0.9}
\definecolor{pink}{rgb}{0.8,0,0.4}

\newtheorem{theorem}{Theorem}[section]

\begin{document}

\title{\bf  A Dual Symmetric Gauss-Seidel  Alternating Direction Method of Multipliers for Hyperspectral Sparse Unmixing }

\author{
Longfei Ren \thanks{School of Information Science and technology, and the Provincial Key Lab of Information Coding and Transmission, Southwest Jiaotong University, No. 999, Xian Road, West Park, High-tech Zone, Chengdu 610031, China (renlf@my.swjtu.edu.cn).},
Chengjing Wang\thanks{{\bf Corresponding author}, School of Mathematics, Southwest Jiaotong University, No. 999, Xian Road, West Park, High-tech Zone, Chengdu 611756, China  ({\tt renascencewang@hotmail.com}).},
Peipei Tang\thanks{School of Computer and Computing Science, Zhejiang University City College, Hangzhou 310015, China (tangpp@zucc.edu.cn). This author's research is supported by the Natural Science Foundation of Zhejiang Province of China under Grant No. LY19A010028 and the Science $\&$ Technology Development Project of Hangzhou, China under Grant No. 20170533B22, 20162013A08.},
and Zheng Ma\thanks{School of Information Science and technology, and the Provincial Key Lab of Information Coding and Transmission, Southwest Jiaotong University, No. 999, Xian Road, West Park, High-tech Zone, Chengdu 610031, China (zma@home.swjtu.edu.cn).
}}

\maketitle

\begin{abstract}
Since sparse unmixing has emerged as a promising approach to hyperspectral unmixing, some spatial-contextual information in the hyperspectral images has been exploited to improve the performance of the unmixing recently.
The total variation (TV) has been widely used to promote the spatial homogeneity as well as the smoothness between adjacent pixels.
However, the computation task for hyperspectral sparse unmixing with a TV regularization term is heavy.
Besides, the convergence of the primal alternating direction method of multipliers (ADMM) for the hyperspectral sparse unmixing with a TV regularization term  has not been explained in details.
In this paper, we design an efficient and convergent dual symmetric Gauss-Seidel ADMM (sGS-ADMM) for hyperspectral sparse unmixing with a TV regularization term.
We also present the global convergence and local linear convergence rate analysis for this algorithm.
As demonstrated in numerical experiments, our algorithm can obviously improve the efficiency of the unmixing compared with the state-of-the-art algorithm. More importantly, we can obtain images with higher quality.
\end{abstract}


\section{Introduction}
\label{intro}
In recent years, the hyperspectral remote sensing technology has been developed significantly.
However, the spatial resolution of hyperspectral images is low and the mixed pixels are widespread in the observed hyperspectral data. The reason lies in the complexity of the ground surface, the limitation of the spectral acquisition approach as well as the restriction of the hyperspectral imaging instruments. How to extract and separate the pure spectral signatures (endmembers) from the mixed pixels and determine the corresponding proportions (abundances) becomes the key issue for the hyperspectral images analysis and its quantification application.

Hyperspectral unmixing, which decomposes mixed pixels into endmembers and corresponding abundances, has obtained much attention in recent decades. It has many practical applications in environmental monitoring, mine detection, agricultural industry, and so on.
The hyperspectral mixture models can be divided into the linear mixing model (LMM) and nonlinear mixing models (NLMMs) \cite{Bioucas2012}.
In the LMM, we assume that the effects of the secondary reflection and multiple scattering have the least influence on the spectral signature.
In the NLMMs, we assume that the mixed spectral signature is synthesized by the endmembers according to some nonlinear relationship.
For the LMM, each pixel in a hyperspectral image can be linearly decomposed into a number of endmembers weighted by their corresponding abundances.
Since the LMM exhibits some practical merits such as its flexibility in different applications and it is also an acceptable approximation of the light scattering mechanisms in many real scenarios \cite{Bioucas2012}, we will focus on the LMM which is also the mainstream of current research on hyperspectral unmixing.

In the literature, the traditional unmixing based on the LMM includes the geometrical based algorithms and the statistical based algorithms \cite{Bioucas2012}.
 The geometrical based algorithms generally require the assumption that all the reflection spectrum curves belong to the same geometrical simplex set. The vertices of the simplex set represent the corresponding endmembers. So identifying the endmembers is equivalent to searching for the vertices.
 The representative algorithms of this class include the vertex component analysis (VCA) algorithm \cite{VCA}, the pixel purity index (PPI) algorithm \cite{PPI}, the simplex growing algorithm (SGA) \cite{SGA}, the minimum volume enclosing simplex (MVES) algorithm \cite{MVES}, the iterative constrained endmembers (ICE) algorithm \cite{ICE} and the minimum volume transform-nonnegative matrix factorization (MVC-NMF) algorithm \cite{MVC-NMF}.
The statistical algorithms, such as the Bayesian techniques, are based on the priori information of the abundances of endmembers for the variability modeling in a natural framework \cite{Dobigeon2009}.

With the explosive development of compressive sensing \cite{Bruckstein2009}, the sparsity based approaches have recently emerged as a promising alternative for hyperspectral unmixing.
The sparsity based approaches aim at finding the optimal subset of a (potentially very large) spectral library in a semisupervised way. The optimal subset is also the best one that can simulate each pixel of a given hyperspectral image.
As shown in \cite{MUSIC-CSR}, the sparsity based approaches have attracted many interests as they do not require the presence of pure pixels in a given scene and do not need to estimate the number of endmembers in the data, which are two obstacles of the traditional unmixing methods.
In practice, the number of endmembers in the real scenarios is far less than that in the spectral library.
This means that the abundances corresponding to the spectral library are sparse.
As a result, the sparsity based approaches are related to the linear sparse regression techniques.

 Iordache  \emph{et al.} \cite{SUnSAL} first added the sparsity constraint to the hyperspectral unmixing model and proposed a sparse unmixing by the variable splitting and augmented Lagrangian (SUnSAL) algorithm. This opens a new gate so that the abundance estimation neither depends on the purity of the spectra nor a good endmember extraction algorithm.
Subsequently,  Iordache  \emph{et al.} \cite{SUnSAL-TV} proposed the sparse unmixing via the variable splitting augmented Lagrangian and total variation (SUnSAL-TV) algorithm to deal with the hyperspectral unmixing model plus a total variation (TV).
The TV, which promotes the spatial homogeneity as well as the smoothness between adjacent pixels, has been widely used in image processing \cite{Rudin1992}\cite{Zhao2013}\cite{Zakharova2017}. So the SUnSAL-TV algorithm which exploits the spatial information in the hyperspectral images can significantly improve the performance of unmixing.
Meanwhile, as shown in \cite{Zhang17} and \cite{Zhang16}, the unmixing performance based on nonisotropic TV is better than that on isotropic TV.
From the other point of view, the collaborative sparse unmixing by the variable splitting and augmented Lagrangian (CLSUnSAL) algorithm in \cite{CLSUnSAL} takes into account the entire abundances matrix globally. The global row sparsity to all pixels in the hyperspectral images is considered as a constraint.
The collaborative SUnSAL-TV (CLSUnSAL-TV) algorithm \cite{CLSUnSAL-TV} takes the combination of the spatial correlation and the global row sparsity into consideration.
Although the numerical experiments show that the SUnSAL-TV algorithm and the CLSUnSAL-TV algorithm work well, the convergences of the two algorithms have not been guaranteed in theory. Therefore we need to design an efficient and convergent algorithm.

As we all know, all the algorithms we mentioned above are essentially special cases of the alternating direction method of multipliers (ADMM) applied to the primal problem. So we may call it primal ADMM. The classical ADMM was originally proposed by Glowinski and Gabay in 1970s. One may see \cite{GlowinskiM} and \cite{GabayM} for details. If the semi-proximal term of the semi-proximal ADMM (SPADMM) \cite{FazelPST} vanishes, then it is the classical ADMM. We refer the readers to \cite{EcksteinY} and \cite{Glowinski} for a better understanding of the historical development of the classic ADMM. We also refer the readers to \cite{FazelPST} and \cite{linear-rate} for the global convergence and linear convergence rate of the SPADMM for convex problems. The ADMM can solve a great deal of problems successfully. However the convergence of the ADMM is only guaranteed for those problems with 2 blocks. (Note that we regard those separable and independent variables as one block).
As for the 3-block (and beyond) problems, the extended version of the ADMM may not converge. One may see a counter example in \cite{CHYY}.
Recently in \cite{LiST2016} and \cite{LiST2019}, a symmetric Gauss-Seidel (sGS) method was designed for the multi-block convex problems with one nonsmooth block.
This opened a new avenue and brought a new insight for handling the problems with nonsmooth blocks.
Most recently, in \cite{Chen2017} and \cite{Chen2018}, the inexact sGS based ADMM (sGS-ADMM) type methods were proposed for solving a class of convex composite optimization problems with two nonsmooth blocks.

The main contributions of this paper are as follows.
Firstly, we design a dual sGS-ADMM to solve the hyperspectral sparse unmixing with a TV regularization term.
Secondly, we present the global convergence and local linear convergence rate of the primal ADMM and the dual sGS-ADMM. As shown in the numerical experiments, the dual sGS-ADMM can obviously improve the efficiency of the unmixing compared with the primal ADMM.

The remaining parts of this paper are organized as follows.
In the next section, we will introduce some basic notations and definitions.
In Section 3, we will describe the model of hyperspectral unmixing.
In Section 4, we will recall the primal ADMM and present the global convergence and local linear convergence rate of the primal ADMM.
In Section 5, we will propose the dual sGS-ADMM. Then we also present the global convergence and local linear convergence rate of the dual sGS-ADMM.
Numerical experiments will demonstrate the efficiency of the proposed algorithm in Section 6.
The conclusion will be discussed in Section 7.

\section{Preliminaries}
\label{sec:1}
In this section, we introduce some basic notations and definitions in convex analysis. We refer the reader to a bible book of convex analysis \cite{Rockafellar1970} for more in-depth contents.

Let $\mathcal{X}$ be a finite dimensional real Hilbert space.
Let $\mathbf{C}$ be a subset of $\mathcal{X}$. The indicator function of $\mathbf{C}$ is defined by ${\delta _\mathbf{C}}(\mathbf{x})$, i.e., ${\delta _\mathbf{C}}(\mathbf{x}) = 0$ if $\mathbf{x} \in \mathbf{C} $ and ${\delta _\mathbf{C}}(\mathbf{x}) = +\infty $ if $\mathbf{x }\notin \mathbf{C}$.
For $\mathbf{X} \in {R^{m \times n}}$, the Frobenius norm of $\mathbf{X}$ is defined by $\|\mathbf{X}\|{_F} = \sqrt {{\rm{trace}}(\mathbf{X}{\mathbf{X}^T})} $.
The ${l_{1,1}}$ and ${l_{2,1}}$ norms of $\mathbf{X}$ are defined by $\|\mathbf{X}\|{_{1,1}}:=\|\mathbf{X}\|_{1} = \sum _{j = 1}^n\sum _{i = 1}^m|{\mathbf{X}_{ij}}|$
and $\|\mathbf{X}|{|_{2,1}}:= \sum _{k = 1}^m\|{\mathbf{X}^k}|{|_2}$, respectively, where $\mathbf{X}^k$ is the $k$-th row of $\mathbf{X}$.
For any given self-adjoint positive semidefinite linear operator $\mathcal{M}:\mathcal{X}\rightarrow\mathcal{X}$, $\textrm{dist}_{\mathcal{M}}(\mathbf{x},\mathbf{S}):=\textrm{inf}_{\mathbf{x}'\in \mathbf{S}}\|\mathbf{x}-\mathbf{x}'\|_{\mathcal{M}}$ for all $\mathbf{x}\in \mathcal{X}$ and $\mathbf{S}\in\mathcal{X}$, where $\|\mathbf{x}\|_{\mathcal{M}}:=\sqrt{\langle \mathbf{x}, \mathcal{M}\mathbf{x} \rangle}$.
The symbols $\mathcal{I}$ denotes the identity mapping.

\emph{Definition 1} \cite[Section 12]{Rockafellar1970}: For any convex function $p: R^{n}\rightarrow R\cup \{+\infty\}$, the conjugate function of $p$ is defined as
\begin{equation*}
{p^*}(\mathbf{y}):= \sup\limits_{\mathbf{x}} \{{\mathbf{y}^T}\mathbf{x} - p(\mathbf{x})\}.
\end{equation*}

\emph{Definition 2} \cite[Definition 2.2.1]{Sun1986}: A set which can be expressed as the intersection of finitely many closed half spaces of $R^n$ is called a convex polyhedron.
A polyhedral set is the union of finitely many convex polyhedrals.
A function is called piecewise quadratic (linear) if its domain is a polyhedral set and it is quadratic (affine) on each of the convex polyhedral which constitutes its domain.

\emph{Definition 3} \cite[Section 31]{Rockafellar1970}: For a given closed proper convex function ${p}:\chi  \to ( - \infty , + \infty ]$, the proximal mapping ${\rm{Prox}}_{p}( \cdot )$ associated with $p$ is defined by
\begin{equation*}
{\rm{Prox}}{_p}(\mathbf{x}): = \mathop {\rm argmin }\limits_\mathbf{u}\left\{ p(\mathbf{u}) + \frac{1}{2}\| {\mathbf{u} - \mathbf{x}} \|_2^2\right\} ,\forall\, \mathbf{x} \in \mathcal{X}.
\end{equation*}

\emph{Definition 4} \cite[Section 24]{Rockafellar1970}: Let $\mathbf{H}$ be a real Hilbert space with an inner product $\langle\cdot,\cdot\rangle$. A multifunction $F:\mathbf{H}\rightrightarrows \mathbf{H} $ is said to be a monotone operator if
\begin{equation*}
\langle \mathbf{z}-\mathbf{z}',\mathbf{w}-\mathbf{w}'\rangle\geq 0,\quad  \forall\, \mathbf{w}\in F(\mathbf{z}),\mathbf{w}'\in F(\mathbf{z}').
\end{equation*}
It is said to be maximal monotone if, in addition, the graph
\begin{equation*}
\textrm{gph}(F):=\left\{(\mathbf{z},\mathbf{w})\in \mathbf{H}\times \mathbf{H}\,|\,\mathbf{w}\in F(\mathbf{z})\right\}
\end{equation*}
is not properly contained in the graph of any other monotone operator $F':\mathbf{H}\rightrightarrows \mathbf{H}$.

\emph{Definition 5} \cite[Section 3.8]{Rockafellar2009}: Let $(\mathbf{x}^0,\mathbf{y}^0)\in \textrm{gph}(F)$. The multi-valued mapping $F:\mathcal{X}\rightrightarrows\mathcal{Y}$ is said to be calm at $\mathbf{x}^0$ for $\mathbf{y}^0$ with modulus $\kappa_0\geq0$ if there exist a neighborhood $V$ of $\mathbf{x}^0$ and a neighborhood $W$ of $\mathbf{y}^0$ such that
\begin{equation*}
F(\mathbf{x})\cap W\subseteq F(\mathbf{x}^0)+\kappa_0\|\mathbf{x}-\mathbf{x}^0\|B_{\mathbf{y}},\quad \forall\ \mathbf{x}\in V,
\end{equation*}
where $B_{\mathbf{y}}$ is the unit ball in $\mathcal{Y}$.

\section{The system model}
\label{sec:2}
For the LMM, we assume that the spectrum of each mixed pixels can be represented as a linear combination of each endmember spectrum in any given spectral band (see \cite{Bioucas2012}). That is, for each mixed pixels, the linear model can be written as
\begin{equation*}
\mathbf{y} = \mathbf{Ms} + \mathbf{n},
\end{equation*}
where $\mathbf{y} := {[{\mathbf{y}_1},{\mathbf{y}_2},...,{\mathbf{y}_L}]^T}$ denotes the measured spectra of the mixed pixels and $L$ denotes the number of bands, $\mathbf{M} := [{\mathbf{m}_1},{\mathbf{m}_2},...{\rm{,}}{{\mathbf{m}}_q}]$ denotes the endmembers matrix, $q$ denotes the number of endmembers and each ${\mathbf{m}_j} := {[{\mathbf{m}_{1j}},{\mathbf{m}_{2j}},...{\rm{,}}{{\mathbf{m}}_{Lj}}]^T}$ denotes the spectra signature of the $j$-th endmembers, $\mathbf{s} := {[{\mathbf{s}_1},{\mathbf{s}_2},...{\rm{,}}{{\mathbf{s}}_q}]^T}$ denotes the abundances of the endmembers and $\mathbf{n} := {[{\mathbf{n}_1},{\mathbf{n}_2},...{\rm{,}}{{\mathbf{n}}_L}]^T}$ denotes the error vector.
According to the physical meaning of the real scene, the abundances need to satisfy the so called abundance nonnegativity constraint (ANC) and the abundance sum constraint (ASC) \cite{Bioucas2012}. That is,
\begin{eqnarray*}
&&\mathbf{s}_i \ge 0, \ \forall\, i = 1,2,...,q,\\
&&\sum _{i = 1}^q{\mathbf{s}_i} = 1.
\end{eqnarray*}

As mentioned previously, the sparsity based approaches proposed by Iordache \emph{et al.} \cite{Bioucas2012} replace the endmembers matrix $\mathbf{M}$ by a known spectral library $\mathbf{A}\in {R^{L \times m}}$ \cite{SUnSAL}. Unmixing then amounts to finding the optimal subset of signatures in $\mathbf{A}$. Specially, we have
\begin{equation*}
\mathbf{y} = \mathbf{Ax} + \mathbf{n},
\end{equation*}
where $\mathbf{x}\in R^{m\times 1}$ denotes the abundances corresponding to the library $\mathbf{A}$.

As shown in \cite{LiHengChao2017}, the conventional hyperspectral sparse unmixing can be uniformly expressed as the following model:
\begin{equation}
\begin{aligned}
\label{model}
\min_ \mathbf{X} \ & \frac{1}{2}\| {\mathbf{AX} - \mathbf{Y}}\|_F^2 + \lambda {\| \mathbf{X}\|_{\rho ,1}} + {\lambda _{TV}}TV(\mathbf{X})\\
 \mbox{s.t.}\ & \mathbf{X} \ge 0,
\end{aligned}
\end{equation}
where
\begin{equation*}
TV(\mathbf{X}):= \sum \limits_{\{i,j\}\in\varepsilon}\|\mathbf{x}_i-\mathbf{x}_j \|_1
\end{equation*}
is a vector extension of the nonisotropic $TV$, $\varepsilon$ denotes the set of horizontal and vertical neighbors of $\mathbf{X}$, $\lambda,\lambda_{TV}\geq0$ are given parameters, $\rho =1$ or $2$,
$\mathbf{Y} = [{\mathbf{y}_1},{\mathbf{y}_2},...,{\mathbf{y}_n}] \in {R^{L \times n}}$ denotes the observed data, $\mathbf{X} = [{\mathbf{x}_1},{\mathbf{x}_2},...,{\mathbf{x}_n}] \in {R^{m \times n}}$ denotes the abundances matrix and $n$ denotes the number of the pixels.
When $\rho=1,\lambda_{TV}=0$, (\ref{model}) is reduced to the sparse unmixing (SUn) model. When $\rho=2,\lambda_{TV}=0$, (\ref{model}) is actually the collaborative sparse unmixing (CLSUn) model. Similarly, we refer to (\ref{model}) with $\rho=1$ as the sparse unmixing with TV (SUnTV) model and $\rho=2$ as the collaborative sparse unmixing with TV (CLSUnTV) model.

For the design of the algorithm, we need to give a more detailed characterization of the TV norm.

Suppose $n = n_r\times n_c$, where $n_r$ and $n_c$ denote the dimensions of the rows and columns of the pixels, respectively.
	Let
	\begin{equation*}\begin{aligned}
	\mathbf{X}:=\left[\mathbf{x}_1,\mathbf{x}_2,\cdots,\mathbf{x}_{n_r},\mathbf{x}_{n_r+1},\mathbf{x}_{n_r+2},\cdots,\mathbf{x}_{2n_r},\right.
	\cdots,\left.\mathbf{x}_{n-n_r+1},\mathbf{x}_{n-n_r+2},\cdots,\mathbf{x}_n\right],
	\end{aligned}\end{equation*}
	
	\begin{equation*}
	\mathbf{X}': = \left( {\begin{array}{*{20}{c}}
		{{\mathbf{x}_1}}&{{\mathbf{x}_{n_r  + 1}}}& \cdots &{{\mathbf{x}_{n - n_r  + 1}}}\\
		{{\mathbf{x}_2}}&{{\mathbf{x}_{n_r  + 2}}}& \cdots &{{\mathbf{x}_{n - n_r  + 2}}}\\
		\vdots & \vdots & \ddots & \vdots \\
		{{\mathbf{x}_{n_r }}}&{{\mathbf{x}_{2n_r }}}& \cdots &{{\mathbf{x}_n}}
		\end{array}} \right)\in R^{(m\times n_r)\times n_c} ,
	\end{equation*}
	
	\begin{equation*}
	{\mathbf{X}''}: = \left( {\begin{array}{*{20}{c}}
		{{\mathbf{x}_1}}&{{\mathbf{x}_2}}& \cdots &{{\mathbf{x}_{n_r }}}\\
		{{\mathbf{x}_{n_r  + 1}}}&{{\mathbf{x}_{n_r  + 2}}}& \cdots &{{\mathbf{x}_{2n_r }}}\\
		\vdots & \vdots & \ddots & \vdots \\
		{{\mathbf{x}_{n - n_r  + 1}}}&{{\mathbf{x}_{n - n_r  + 2}}}& \cdots &{{\mathbf{x}_n}}
		\end{array}} \right)\in R^{(m\times n_c)\times n_r},
	\end{equation*}
	
	Define two linear operators $\mathcal{B}:R^{(m\times n_r)\times n_c}\rightarrow R^{(m\times n_r)\times (n_c-1)} $ and $\mathcal{C}:R^{(m\times n_c)\times n_r}\rightarrow R^{(m\times n_c)\times (n_r-1)} $ to compute the horizontal differences between the neighboring pixels of $\mathbf{X}'$ and $\mathbf{X}''$.
	\begin{equation*}
	\mathcal{B}{{\mathbf{X}'}}: = \left( {\begin{array}{*{20}{c}}
		\mathbf{c}_1&\mathbf{c}_{n_r  + 1}& \cdots &\mathbf{c}_{n - 2n_r  + 1}\\
		\mathbf{c}_2&\mathbf{c}_{n_r  + 2}& \cdots &\mathbf{c}_{n - 2n_r  + 2}\\
		\vdots & \vdots & \ddots & \vdots \\
		\mathbf{c}_{n_r }&\mathbf{c}_{2n_r }& \cdots &\mathbf{c}_{n-{n_r}}
		\end{array}} \right),
	\end{equation*}
	
	\begin{equation*}
	\mathcal{C}{{\mathbf{X}''}}: = \left( {\begin{array}{*{20}{c}}
		\mathbf{e}_1&\mathbf{e}_2& \cdots &\mathbf{e}_{n_r-1}\\
		\mathbf{e}_{n_r+1}&\mathbf{e}_{n_r  + 2}& \cdots &\mathbf{e}_{2n_r-1}\\
		\vdots & \vdots & \ddots & \vdots \\
		\mathbf{e}_{n-n_r+1 }&\mathbf{e}_{n-n_r+2}& \cdots &\mathbf{e}_{n-1}
		\end{array}} \right),
	\end{equation*}
	where $\mathbf{c}_i=\mathbf{x}_{i+n_r}-\mathbf{x}_{i}$ $(i=1,\cdots,n-n_r)$ and $\mathbf{e}_i=\mathbf{x}_{i+1}-\mathbf{x}_{i}$  $(i=1,\cdots,n_r-1, n_{r}+1,\cdots, 2 n_r-1,\cdots, n-n_r+1 ,\cdots, n-1)$.
	Define two linear operators $\mathcal{\hat{H}}_v:R^{m\times n}\rightarrow R^{m\times (n-n_r)}$ and $\mathcal{\hat{H}}_h:R^{m\times n}\rightarrow R^{m\times (n-n_c)}$:
\begin{eqnarray*}
\mathcal{\hat{H}}_{v}\mathbf{X}&=&[\mathbf{c}_1,\mathbf{c}_2,\ldots,\mathbf{c}_{n-n_r}],\\
\mathcal{\hat{H}}_{h}\mathbf{X}&=&[\mathbf{e}_1,\ldots,\mathbf{e}_{n_{r}-1},\mathbf{e}_{n_{r}+1},\cdots,\mathbf{e}_{2n_{r}-1},
 \cdots,\mathbf{e}_{n-n_{r}+1},\cdots,\mathbf{e}_{n-1}].
\end{eqnarray*}
Then the problem (\ref{model}) can be reformulated to
\begin{equation}\begin{aligned}\label{model1}
\mathop {\min }\limits_\mathbf{X} \ \frac{1}{2}\| {\mathbf{AX} - \mathbf{Y}} \|_F^2 + \lambda {\| \mathbf{X} \|_{\rho ,1}}  + {\lambda _{TV}}{\| {{\mathcal{\hat{H}}}_{v}\mathbf{X}} \|_1}+{\lambda _{TV}}{\| {{\mathcal{\hat{H}}}_{h}\mathbf{X}} \|_1} + {\delta _ {R_{+}^{m\times n}} }\left( \mathbf{X} \right).
\end{aligned}\end{equation}

We should mention that there is a slight difference between the way we handle the TV norm with that in \cite{SUnSAL-TV}.
In our framework, we assume the (standard) reflexive boundary condition holds, which means that the rightmost (lowest) boundaries have no right (lower) neighboring pixels.
Note that the reflexive boundary condition for the hyperspectral images in this paper is an extension of the reflexive boundary condition for the 2- dimensional images. One may see \cite{Beck2009} and \cite{Zuo2011} for the reflexive boundary condition of the 2-dimensional images.
Actually it is more reasonable since it is what happens in the real scenarios.
Instead, in \cite{SUnSAL-TV}, they assume periodic boundaries, i.e. the neighboring pixels of the rightmost (lowest) boundaries are the leftmost (highest) boundaries pixels only for the convenience of adopting the fast Fourier transform. In order to be consistent, we keep the way in \cite{SUnSAL-TV} when using the primal ADMM. In \cite{SUnSAL-TV}, they define two linear operators $\mathcal{H}_h:R^{m\times n}\rightarrow R^{m\times n}$ and $\mathcal{H}_v:R^{m\times n}\rightarrow R^{m\times n}$ to compute the horizontal and vertical differences between the neighboring pixels of $\mathbf{X}$ as follows
\begin{eqnarray*}
\mathcal{H}_{h}\mathbf{X}=\left[\mathbf{d}_1,\mathbf{d}_2,\ldots,\mathbf{d}_n\right],\
\mathcal{H}_{v}\mathbf{X}=\left[\mathbf{b}_1,\mathbf{b}_2,\ldots,\mathbf{b}_n\right],
\end{eqnarray*}
where $\mathbf{d}_i=\mathbf{x}_i-\mathbf{x}_{i_h}$ and $\mathbf{b}_i=\mathbf{x}_i-\mathbf{x}_{i_v}$ $(i=1,\cdots,n)$, with $i$ denoting an index of a pixel, $i_h$ and $i_v$ denoting the indices of the corresponding horizontal and vertical neighbors. Let
\begin{equation*}
\mathcal{H}\mathbf{X} = \left[ {\begin{array}{*{20}{c}}
{{\mathcal{H}_h}\mathbf{X}}\\
{{\mathcal{H}_v}\mathbf{X}}
\end{array}} \right].
\end{equation*}
An equivalent form of the problem (\ref{model}) is
\begin{equation}
\begin{aligned}
\label{P0}
\min_ \mathbf{X} \ & \frac{1}{2}\| {\mathbf{AX} - \mathbf{Y}}\|_F^2 + \lambda {\| \mathbf{X}\|_{\rho ,1}} +{\lambda _{TV}}\|\mathcal{H}\mathbf{X}\|_{1}\\
 \mbox{s.t.}\ & \mathbf{X} \ge 0.
\end{aligned}
\end{equation}

\section{The primal ADMM}
In this section we first recall the SUnSAL-TV algorithm and the CLSUnSAL-TV algorithm, both of which are essentially the primal ADMM.
Then we present the global convergence and local linear convergence rate of the primal ADMM.
\subsection{The primal ADMM}
We can reformulate (\ref{P0}) equivalently by introducing some slack variables
\begin{equation}
\begin{aligned}
\min_ {\mathbf{\widetilde{D}},{\mathbf{D}_1},{\mathbf{D}_2},{\mathbf{D}_3},{\mathbf{D}_4},{\mathbf{D}_5}}\ & \frac{1}{2}\|{\mathbf{D}_1} - \mathbf{Y}\|_F^2 + \lambda \|{\mathbf{D}_2}\|_{\rho,1} + {\lambda _{TV}}\|{\mathbf{D}_4}\|_{1}  + {\delta _ {R_{+}^{m\times n}} }({\mathbf{D}_5})\\
 \mbox{s.t.}\quad\qquad\ & \mathbf{D}_1=\mathbf{A\widetilde{D}},\,\mathbf{D}_2=\mathbf{\widetilde{D}},\,\mathbf{D}_3=\mathbf{\widetilde{D}},\,
 \mathbf{D}_4=\mathcal{H}\mathbf{D}_3,\,\mathbf{D}_5=\mathbf{\widetilde{D}},\\
\end{aligned}
\end{equation}
where $\mathbf{\widetilde{D}}\in R^{m\times n}$, $\mathbf{D}_1\in R^{L\times n}$, $\mathbf{D}_{2}\in R^{m\times n}$, $\mathbf{D}_3\in R^{m\times n}$, $\mathbf{D}_4\in R^{2m\times n}$, $\mathbf{D}_5\in R^{m\times n}$.

The above problem can also be written in the following form
\begin{equation}
\begin{aligned}
\label{compact P1}
\min_ {\mathbf{M},\mathbf{N}}\ & f(\mathbf{M})+g(\mathbf{N})\\
\mbox{s.t.}\ &\mathbf{F}\mathbf{M}+\mathbf{G}\mathbf{N}=0,
\end{aligned}
\end{equation}
where
\begin{eqnarray*}
\mathbf{M}&=&\left(\mathbf{D}_1^T,\mathbf{D}_2^T,\mathbf{D}_3^T,\mathbf{D}_5^T\right)^T\in R^{(3m+L)\times n},\
\mathbf{N}=\left(\mathbf{\widetilde{D}}^T,\mathbf{D}_4^T\right)^T\in R^{3m\times n},\\
f(\mathbf{M})&=&\frac{1}{2}\|{\mathbf{D}_1} - \mathbf{Y}\|_F^2 + \lambda \|{\mathbf{D}_2}\|_{\rho,1}+ {\delta _ {R_{+}^{m\times n}} }({\mathbf{D}_5}),\
g(\mathbf{N})={\lambda _{TV}}\|{\mathbf{D}_4}\|_{1},
\end{eqnarray*}
\begin{eqnarray*}
\mathbf{F}=\left( {\begin{array}{*{20}{c}}
{  \mathbf{- I}_{L\times L}}&0_{L\times m}&0_{L\times m}&0_{L\times m}\\
0_{m\times L}&{  \mathbf{- I}_{m\times m}}&0_{m\times m}&0_{m\times m}\\
0_{m\times L}&0_{m\times m}&{  \mathbf{- I}_{m\times m}}&0_{m\times m}\\
0_{2m\times L}&0_{2m\times m}&-\mathcal{H}&0_{2m\times m}\\
0_{m\times L}&0_{m\times m}&0_{m\times m}&{  \mathbf{- I}_{m\times m}}
\end{array}} \right),\
\mathbf{G}=\left( {\begin{array}{*{20}{c}}
\mathbf{A}_{L\times m}&0_{L\times 2m}\\
\mathbf{I}_{m\times m}&0_{m\times 2m}\\
\mathbf{I}_{m\times m}&0_{m\times 2m}\\
{{0}_{2m\times m}}&\mathbf{I}_{2m\times 2m}\\
\mathbf{I}_{m\times m}&{ 0_{m\times 2m}}
\end{array}} \right).
\end{eqnarray*}
Let $\sigma>0$ be a given positive number. $\mathbf{\Lambda}:=\left(\mathbf{\Lambda}_1^T,\mathbf{\Lambda}_2^T,\mathbf{\Lambda}_3^T,\mathbf{\Lambda}_4^T,\mathbf{\Lambda}_5^T\right)^T\in R^{(5m+L)\times n}$ denotes the Lagrange multipliers of the constraints. The augmented Lagrangian function for the problem (\ref{compact P1}) is
\begin{equation}\begin{aligned}\label{ALP1}
 L_{\sigma}(\mathbf{M},\mathbf{N};\mathbf{\Lambda})=f(\mathbf{M})+g(\mathbf{N})+\frac{\sigma}{2}\|\mathbf{FM}+\mathbf{GN}-\sigma^{-1}\mathbf{\Lambda}\|_F^2-\frac{1}{2\sigma}\|\mathbf{\Lambda}\|_F^2.
\end{aligned}\end{equation}

The primal ADMM can be presented in Algorithm 1.
\begin{algorithm}[!h]
\caption{The primal ADMM}
\begin{algorithmic}[1]

\REQUIRE
Select an initial point $(\mathbf{M}^0,\mathbf{N}^0;\mathbf{\Lambda}^0)$. Set $k=0$, choose $\sigma>0$ and $\tau\in\left(0,\frac{1+\sqrt{5}}{2}\right)$. Iterate the following steps until the stopping criterion is satisfied:

\item [\bf Step 1.]
Compute:
\begin{eqnarray*}
&&\mathbf{M}^{k+1}=\mathop {{\mathop{\rm argmin}\nolimits} }\limits_{\mathbf{M}} \{ L_{\sigma}(\mathbf{M},{\mathbf{N}^k};{\mathbf{\Lambda}^k})\} ,\\
&&\mathbf{N}^{k+1}=\mathop {{\mathop{\rm argmin}\nolimits} }\limits_{\mathbf{N}} \{ L_{\sigma}(\mathbf{M}^{k+1},{\mathbf{N}};{\mathbf{\Lambda}^k})\}.
\end{eqnarray*}
\item [\bf Step 2.]
Update:
\begin{eqnarray*}
&&\mathbf{\Lambda}^{k+1}=\mathbf{\Lambda}^{k}-\tau\sigma(\mathbf{FM}^{k+1}+\mathbf{GN}^{k+1}).
\end{eqnarray*}
\end{algorithmic}
\end{algorithm}

The reader may refer to \cite{SUnSAL-TV} for a better understanding of the details of Algorithm 1. \footnote{We have to emphasize that in \cite{SUnSAL-TV} the authors deal with the problem (1) in three blocks instead of two blocks when they applied the primal ADMM. In fact, the problem can be regarded as two blocks. So Algorithm 1 in this paper is a little different with Algorithm 1 in \cite{SUnSAL-TV}. It not only is faster but also has mathematically guaranteed convergence theory. Actually the code for the algorithm in \cite{SUnSAL-TV} is in accordance with our two-block primal ADMM.}
\subsection{Convergence analysis }
In this subsection, we discuss the global convergence and the local linear convergence rate of the primal ADMM.

Suppose that $(\mathbf{M},\mathbf{N})\in R^{(3m+L)\times n}\times R^{3m\times n}$ is an optimal solution to the problem (\ref{compact P1}).
If there exists $\mathbf{\Lambda}\in R^{(5m+L)\times n}$ such that $(\mathbf{M},\mathbf{N},\mathbf{\Lambda})$ satisfies the following KKT system
\begin{equation}\label{compactKKT0}
\left\{ {\begin{array}{*{20}{l}}
	{0 \in \partial f(\mathbf{M}) - \mathbf{F}^T\mathbf{\Lambda}},\\
	{0 \in \partial g(\mathbf{N}) - \mathbf{G}^T\mathbf{\Lambda}},\\
	{\mathbf{F}\mathbf{M}+\mathbf{G}\mathbf{N}=0},
	\end{array}} \right.
\end{equation}
then $(\mathbf{M},\mathbf{N},\mathbf{\Lambda})$ is a KKT point for the problem (\ref{compact P1}), where $\partial {f}$ and $\partial {g}$ are the subdifferential mappings of $f$ and $g$.
Let $\Omega$ be the solution set of the KKT system (\ref{compactKKT0}).
Let $\mathcal{B}:R^{(5m+L)\times n}\rightarrow R^{(3m+L)\times n}\times R^{3m\times n}$ be a linear operator such that its adjoint $\mathcal{B}^{*}(\mathbf{M},\mathbf{N})=\mathbf{F}\mathbf{M}+\mathbf{G}\mathbf{N}$.
	For any $\mathbf{u}:=(\mathbf{M},\mathbf{N},\mathbf{\Lambda})$, the KKT mapping is defined by
	\begin{equation*}
	\mathbf{Q}(\mathbf{u}): = \left( {\begin{array}{*{20}{c}}
		{{\mathbf{M}} - {{{\mathop{\rm Prox}\nolimits} }_{f}}(\mathbf{M}+\mathbf{F}^T\mathbf{\Lambda})}\\
		{{\mathbf{N}} - {{{\mathop{\rm Prox}\nolimits} }_{g}}(\mathbf{N}+\mathbf{G}^T\mathbf{\Lambda})}\\
		\mathbf{F}\mathbf{M}+\mathbf{G}\mathbf{N}
		\end{array}} \right).
	\end{equation*}

Since the subdifferential mappings of the proper closed convex function $ f$ and $g$ are maximally monotone \cite[Theorem 12.17]{Rockafellar1998}, there exist two self-adjoint and positive semidefinite linear operators $\Sigma_{f}$ and $\Sigma_{g}$ such that for all $\mathbf{s}, \mathbf{s}'\in \textrm{dom}(f)$, $\varsigma\in\partial f(\mathbf{s})$ and $\varsigma'\in\partial f(\mathbf{s}')$
\begin{equation*}
f(\mathbf{s})\geq f(\mathbf{s}')+\langle\varsigma',\mathbf{s}-\mathbf{s}'\rangle+\frac{1}{2}\|\mathbf{s}-\mathbf{s}'\|_{\Sigma_f}^2,\
\langle\varsigma-\varsigma',\mathbf{s}-\mathbf{s}'\rangle\geq\|\mathbf{s}-\mathbf{s}'\|_{\Sigma_f}^2,
\end{equation*}
and for all $\mathbf{t}, \mathbf{t}'\in \textrm{dom}(g)$, $\upsilon\in\partial g(\mathbf{t})$ and $\upsilon'\in\partial g(\mathbf{t}')$
\begin{equation*}
g(\mathbf{t})\geq g(\mathbf{t}')+\langle\upsilon',\mathbf{t}-\mathbf{t}'\rangle+\frac{1}{2}\|\mathbf{t}-\mathbf{t}'\|_{\Sigma_g}^2,\
\langle\upsilon-\upsilon',\mathbf{t}-\mathbf{t}'\rangle\geq\|\mathbf{t}-\mathbf{t}'\|_{\Sigma_g}^2.
\end{equation*}

For the convergence of the primal ADMM, we need the following assumption.

\emph{Assumption 1}: The KKT system (\ref{compactKKT0}) has a non-empty solution set.\label{Assumption: KKT solution0}

Now we present the following global convergence and linear convergence rate of the primal ADMM. It is a special case of the result derived by the combination of \cite[Theorem B.1]{FazelPST} and \cite[Theorem 2]{linear-rate}.
\begin{theorem}
		Suppose that Assumption 1 holds. Let $\tau\in(0,\frac{1+\sqrt{5}}{2})$,
		then there exists a KKT point $\bar{\mathbf{u}}:=(\overline{\mathbf{M}},\overline{\mathbf{N}}, \overline{\mathbf{\Lambda}})\in{\Omega}$ such that the sequence $\{(\mathbf{M}^k,\mathbf{N}^k,\mathbf{\Lambda}^k)\}$ generated by Algorithm 1 converges to $\bar{\mathbf{u}}$.
		Assume that $\mathbf{Q}^{-1}$ is calm at the origin for $\bar{\mathbf{u}}$ with modulus $\eta>0$, i.e., there exists $r>0$ such that
		\begin{equation*}
		\textrm{dist}(\mathbf{u},{\Omega})\leq\eta\|\mathbf{Q}(\mathbf{u})\|,\quad \forall\, \mathbf{u}\in\{\mathbf{u}:\|\mathbf{u}-\bar{\mathbf{u}}\|\leq r\}.
		\end{equation*}
		Then there exists an integer $\bar{k}\geq 1$ such that for all $k\geq \bar{k}$,
		\begin{equation*}\begin{aligned}
		&\textrm{dist}_{\mathcal{M}}^{2}(\mathbf{u}^{k+1},{\Omega}) \leq \mu \textrm{dist}_{\mathcal{M}}^{2}(\mathbf{u}^{k},{\Omega}),
		\end{aligned}\end{equation*}
		where $\mu\in(0,1)$ and
		\begin{equation*}
		\mathcal{M}:=\textrm{Diag}\left(\Sigma_f,\Sigma_g+\sigma\mathcal{I},(\tau\sigma)^{-1}\mathcal{I}\right)+\frac{\sigma}{4}\mathcal{B}\mathcal{B}^{*},
		\end{equation*}
		Moreover, there exists a positive number $\varsigma\in[\mu,1)$ such that for all $k\geq 1$
		\begin{equation*}\begin{aligned}
		&\textrm{dist}_{\mathcal{M}}^{2}(\mathbf{u}^{k+1},{\Omega}) \leq \varsigma  \textrm{dist}_{\mathcal{M}}^{2}(\mathbf{u}^{k},{\Omega}).
		\end{aligned}\end{equation*}
\end{theorem}

\emph{Proof}: We note that Algorithm 1 is actually the SPADMM. According to \cite[Theorem B.1]{FazelPST}, we only need to prove the following two conditions hold.
\begin{equation}
\label{cond:psd1}
\Sigma_{f}+\sigma \mathbf{F}^T\mathbf{F}\succ0,
\end{equation}
\begin{equation}
\label{cond:psd2}
\Sigma_{g}+\sigma \mathbf{G}^T\mathbf{G}\succ0.
\end{equation}
Since $\mathbf{F}$ and $\mathbf{G}$ are full column rank,  we know that $\mathbf{F}^T\mathbf{F}$ and $\mathbf{G}^T\mathbf{G}$ are positive definite.
In addition, $\Sigma_{f}$ and $\Sigma_{g}$ are positive semidefinite, so we can obtain that \eqref{cond:psd1} and  \eqref{cond:psd2} hold directly.

\emph{Remark 1}: Theorem 1 just gives a general result on the linear rate of convergence for Algorithm 1. It is obvious that the core assumption of Theorem 1 is the calmness condition, which is often too strict. For the case of $\rho=1$ in the problem (\ref{model}), $f$ and $g$ are piecewise linear-quadratic functions. From \cite[Proposition 2.24]{Sun1986}, we know that $\mathbf{Q}^{-1}$ is piecewise polyhedral, and furthermore the calmness condition holds automatically by \cite{Robinson1981}.

\section{A dual sGS-ADMM}
In this section, we first propose an sGS-ADMM for the problem \eqref{model1},
then we present the global convergence and local linear convergence rate of the algorithm.

\subsection{A dual sGS-ADMM}
Now we introduce a more efficient dual sGS-ADMM to solve the problem (\ref{model1}).

Let
\begin{eqnarray*}
p(\mathbf{X})=\lambda {\| \mathbf{X} \|_{\rho ,1}} + {\lambda _{TV}}{\| {{\mathcal{\hat{H}}}_{v}\mathbf{X}} \|_1}+ {\delta _ {R_{+}^{m\times n}} }( \mathbf{X} ),\
q(\mathbf{X})= {\lambda _{TV}}{\| {{\mathcal{\hat{H}}}_{h}\mathbf{X}} \|_1}.
\end{eqnarray*}
Then the problem \eqref{model1} can be written in a simple form as below
\begin{equation}\label{P}
\mathop {\min\ }\limits_\mathbf{X} \frac{1}{2}\| {\mathbf{AX} - \mathbf{Y}} \|_F^2+p(\mathbf{X})+q(\mathbf{X}).
\end{equation}

By introducing three slack variables $\mathbf{U}_1$, $\mathbf{U}_2$ and $\mathbf{U}_3$, \eqref{P} can be written as
\begin{equation}
\begin{aligned}
\label{P2}
\min_{\mathbf{X},{\mathbf{U}_1},{\mathbf{U}_2},{\mathbf{U}_3}} \ & \frac{1}{2}\| {{\mathbf{U}_3}} \|_F^2 + p({\mathbf{U}_1}) + q(\mathbf{U}_2) \\
\mbox{s.t.}\qquad & \mathbf{X}-\mathbf{U}_1=0,\,\mathbf{X}-\mathbf{U}_2=0,\,
\mathbf{AX}-\mathbf{Y}-\mathbf{U}_3=0,\\
\end{aligned}
\end{equation}
where $\mathbf{U}_1\in R^{m\times n}$, $\mathbf{U}_2\in R^{m\times n}$, $\mathbf{U}_3\in R^{L\times n}$.

The dual of the problem (\ref{P2}) is
\begin{equation}
\begin{aligned}
\label{D2}
\min_{{\mathbf{V}_1},{\mathbf{V}_2},{\mathbf{V}_3}}\ & {p^*}( - {\mathbf{V}_1}) + {q^*}( - {\mathbf{V}_2}) + \frac{1}{2}\| {{\mathbf{V}_3}} \|_F^2 - \left\langle {{\mathbf{V}_3},\mathbf{Y}}\right\rangle \\
\mbox{s.t.}\quad\ &-{\mathbf{V}_1} - {\mathbf{V}_2} - {\mathbf{A}^T}{\mathbf{V}_3}=0,
\end{aligned}
\end{equation}
where $\mathbf{V}_1\in R^{m \times n}$, $\mathbf{V}_2\in R^{m \times n}$, $\mathbf{V}_3\in R^{L \times n}$.

Let $\sigma>0$, the augmented Lagrangian function for the problem (\ref{D2}) is
\begin{equation*}\begin{aligned}
&{L_\sigma }({\mathbf{V}_1},{\mathbf{V}_2},{\mathbf{V}_3};\mathbf{X}) = {p^*}( - {\mathbf{V}_1}) + {q^*}( - {\mathbf{V}_2}) + \frac{1}{2}\| {{\mathbf{V}_3}} \|_F^2 - \langle {{\mathbf{V}_3},\mathbf{Y}} \rangle +\\&\frac{\sigma }{2}\| {-{\mathbf{V}_1} - {\mathbf{V}_2} - {\mathbf{A}^T}{\mathbf{V}_3} - {\sigma ^{ - 1}}\mathbf{X}} \|_F^2  - \frac{1}{{2\sigma }}\| \mathbf{X} \|_F^2.
\end{aligned}\end{equation*}

Suppose that $\hat{\delta}_{\mathbf{V}_3}$ and ${\delta}_{\mathbf{V}_3}$ are given tolerance vectors.
The dual sGS-ADMM which is a direct application of the sGS-ADMM in \cite{LiST2016} or \cite{LiST2019} can be presented in Algorithm 2.\footnote{In Algorithm 2, $\mathbf{V}_3$ is updated twice because we use the sGS decomposition to solve the subproblem. For more details, one may refer to \cite{LiST2019}.}
\begin{algorithm}[!t]
\caption{A dual sGS-ADMM}
\begin{algorithmic}[2]

\REQUIRE
Select an initial point $ \left(\mathbf{V}_1^0,\mathbf{V}_2^0,\mathbf{V}_3^0;\mathbf{X}^0\right)$. Set $k=0$, choose $\sigma>0$, $\tau\in\left(0,\frac{1+\sqrt{5}}{2}\right)$ and $\{\tilde\epsilon_k\}_{k\geq0}$ be a summable sequence of nonnegative numbers. Iterate the following steps until the stopping criterion is satisfied:

\item [\bf Step 1.]
Compute:
\begin{eqnarray*}
&&\hat{\mathbf{V}}_3^{k}  = \mathop {{\mathop{\rm argmin}\nolimits} }\limits_{\mathbf{V}_3} \{ {L_\sigma }(\mathbf{V}_1^k,\mathbf{V}_2^k,\mathbf{V}_3;{\mathbf{X}^k})-\langle\hat{\delta}^{k}_{\mathbf{V}_3},\mathbf{V}_{3}\rangle\},\\
&&{\mathbf{V}}_1^{k+1}  = \mathop {{\mathop{\rm argmin}\nolimits} }\limits_{\mathbf{V}_1} \{ {L_\sigma }(\mathbf{V}_1,\mathbf{V}_2^k,\hat{\mathbf{V}}_3^{k};{\mathbf{X}^k})\},\\
&&\mathbf{V}_3^{k+1}  = \mathop {{\mathop{\rm argmin}\nolimits} }\limits_{\mathbf{V}_3} \{ {L_\sigma }(\mathbf{V}_1^{k+1},\mathbf{V}_2^k,\mathbf{V}_3;\mathbf{X}^k) -\langle\delta^{k}_{\mathbf{V}_3},\mathbf{V}_{3}\rangle\},\\
&&\mathbf{V}_2^{k+1}  = \mathop {{\mathop{\rm argmin}\nolimits} }\limits_{\mathbf{V}_2} \{ {L_\sigma }(\mathbf{V}_1^{k+1},\mathbf{V}_2,\mathbf{V}_3^{k+1};{\mathbf{X}^k})\},
\end{eqnarray*}
where
\begin{eqnarray*}
&&\hat{\delta}^{k}_{\mathbf{V}_3}=\nabla_{\mathbf{V}_3}L_{\sigma}({\mathbf{V}_1^k},{\mathbf{V}_2^k},{\hat{\mathbf{V}}_3^k};\mathbf{X}^k)\ {\rm{with}} \ \|\hat{\delta}^{k}_{\mathbf{V}_3}\|_{F}\leq \tilde\epsilon_k,\\
&&{\delta}^{k}_{\mathbf{V}_3}=\nabla_{\mathbf{V}_3}L_{\sigma}({\mathbf{V}_1^{k+1}},{\mathbf{V}_2^k},{{\mathbf{V}}_3^{k+1}};\mathbf{X}^k)\ {\rm{with}} \ \|{\delta}^{k}_{\mathbf{V}_3}\|_{F}\leq \tilde\epsilon_k.\\
\end{eqnarray*}
\item [\bf Step 2.]
Update:
\begin{eqnarray*}
&&\mathbf{X}^{k+1}=\mathbf{X}^{k}-\tau\sigma(-{\mathbf{V}}_1^{k+1} - {\mathbf{V}}_2^{k+1} - {\mathbf{A}^T}{\mathbf{V}}_3^{k+1}).
\end{eqnarray*}
\end{algorithmic}
\end{algorithm}

Now we provide the details of Step 1 in Algorithm 2.
Finding the minimum of $ L_{\sigma}({\mathbf{V}_1},{\mathbf{V}_2},{\mathbf{V}_3};\mathbf{X})$ with respect to $ \mathbf{V}_3 $ is equivalent to solving the following problem.
\begin{equation*}\begin{aligned}
\nabla_{\mathbf{V}_3}L_{\sigma}({\mathbf{V}_1},{\mathbf{V}_2},{\mathbf{V}_3};\mathbf{X})={\mathbf{V}_3} - \mathbf{Y} + \sigma \mathbf{A}({\mathbf{V}_1} + {\mathbf{V}_2} + {\mathbf{A}^T}{\mathbf{V}_3} + {\sigma ^{ - 1}}\mathbf{X})=0.
\end{aligned}\end{equation*}
That is, we compute $ \mathbf{V}_3 $ by solving the following linear system of equations
\begin{equation*}\label{V_3}
(\mathbf{I} + \sigma \mathbf{A}{\mathbf{A}^T}){\mathbf{V}_3} =\mathbf{Y} - \sigma \mathbf{A}{\mathbf{V}_1} - \sigma \mathbf{A}{\mathbf{V}_2} - \mathbf{AX}.
\end{equation*}

Making use of the Moreau identity ${{\mathop{\rm Prox}\nolimits} _{\sigma p}}(\mathbf{x}) + \sigma \mathop{\rm Prox}_{{p^*}/\sigma }(\mathbf{x}/\sigma ) =\mathbf{ x}$,
we can get $\mathbf{V}_1,\mathbf{V}_2$ in closed forms as below
\begin{equation*}\label{V_1}
\mathbf{V}_1=\frac{1}{\sigma }{{\mathop{\rm Prox}\nolimits} _{\sigma p}}(\sigma \mathbf{C}_1)-\mathbf{C}_1,\
\mathbf{V}_2=\frac{1}{\sigma }{{\mathop{\rm Prox}\nolimits} _{\sigma q}}(\sigma \mathbf{C}_2)-\mathbf{C}_2,
\end{equation*}
where
\begin{equation*}
\mathbf{C}_1=\mathbf{V}_2+\mathbf{A}^{T}\mathbf{V}_3+\sigma^{-1}\mathbf{X},\
\mathbf{C}_2=\mathbf{V}_1+\mathbf{A}^{T}\mathbf{V}_3+\sigma^{-1}\mathbf{X}.
\end{equation*}

We discuss different cases for different $p$ and $q$.

(i) When $p(\cdot)=\lambda\|\cdot\|_{1}+\delta_{R_{+}^{m\times n}}(\cdot)+{\lambda _{TV}}{\| {{\mathcal{\hat{H}}}_{v}(\cdot)} \|_1}$.
The following proposition states that the proximal mapping of $\sigma p$ can be decomposed into the composition of the three proximal mappings.

\emph{Proposition 1}: For any $\sigma >0$, it holds that
\begin{equation*}
{{\mathop{\rm Prox}\nolimits} _{{\sigma p_{}}}} =  {{\mathop{\rm Prox}\nolimits} _{\sigma\lambda\|\cdot\|_{1}}}\circ\mathop{\rm Prox}\nolimits_{\delta_{R_{+}^{m\times n}}} \circ {{\mathop{\rm Prox}\nolimits} _{\sigma\lambda _{TV}\| {{\mathcal{\hat{H}}}_{v}(\cdot)} \|_1}}.
\end{equation*}
\emph{Proof}: Let $r(\cdot)=\sigma\lambda\|\cdot\|_{1}+\delta_{R_{+}^{m\times n}}(\cdot)$. From the equivalence of (iii) and (iv) in \cite[Theorem 4]{Yu2013}, we can obtain
\begin{equation*}
{{\mathop{\rm Prox}\nolimits} _{{r}}} ={{\mathop{\rm Prox}\nolimits} _{\sigma\lambda\|\cdot\|_{1}}} \circ \mathop{\rm Prox}\nolimits_{\delta_{R_{+}^{m\times n}}}.
\end{equation*}
Since $r$ is permutation invariant, that is, $r_1(Q\cdot)=r_1(\cdot)$ for all permutation $Q$, then by \cite[Corollary 4]{Yu2013}, we can obtain
\begin{equation}
{{\mathop{\rm Prox}\nolimits} _{{\sigma p_{}}}} =\mathop{\rm Prox}\nolimits_{r} \circ {{\mathop{\rm Prox}\nolimits} _{\sigma\lambda _{TV}\| {{\mathcal{\hat{H}}}_{v}(\cdot)} \|_1}}
\end{equation}
holds directly.

We define $p_{1}(\cdot)=\delta_{R_{+}^{m\times n}}(\cdot)+\lambda_{TV}{\| {{\mathcal{\hat{H}}}_{v}(\cdot)} \|_1}$ for convenience.
From Proposition 1, we can obtain
\begin{equation*}\begin{aligned}
{{\mathop{\rm Prox}\nolimits} _{{\sigma p}}}(\sigma \mathbf{C}_1) &= {{\mathop{\rm Prox}\nolimits} _{\sigma \lambda \| {\cdot} \|_{1}}}({{\mathop{\rm Prox}\nolimits} _{{\sigma p_1}}}(\sigma \mathbf{C}_1)) \\&= \textrm{sign}({{\mathop{\rm Prox}\nolimits} _{{\sigma p_1}}}(\sigma \mathbf{C}_1)) \circ \max (| {{\mathop{\rm Prox}\nolimits} _{{\sigma p_1}}}(\sigma \mathbf{C}_1)| - \sigma\lambda,0).
\end{aligned}\end{equation*}
Note that the proximal mapping of $\sigma p_1$ is
\begin{equation*}\begin{aligned}
{{\mathop{\rm Prox}\nolimits} _{\sigma {p_1}}}(\sigma \mathbf{C}_1) = \mathop {\rm argmin }\limits_{\mathbf{Z}}\Big\{{\delta_{R_{+}^{m\times n}}(\mathbf{Z})} + \sigma {\lambda _{TV}}{\| {{\mathcal{\hat{H}}}_{v}\mathbf{Z}} \|_1} + \frac{1}{2}{\| {\mathbf{Z} - \sigma \mathbf{C}_1} \|_F^2}\Big\}.
\end{aligned}\end{equation*}
Let
\begin{equation*}\begin{aligned}\label{U}
\mathbf{Z}=\left[\mathbf{z}_1,\mathbf{z}_2,\cdots,\mathbf{z}_{n_r},\mathbf{z}_{n_r+1},\mathbf{z}_{n_r+2},\cdots,\mathbf{z}_{2n_r},\right.
\cdots,\left.\mathbf{z}_{n-n_r+1},\mathbf{z}_{n-n_r+2},\cdots,\mathbf{z}_n\right],
\end{aligned}\end{equation*}
\begin{equation*}\begin{aligned}\label{U}
\mathbf{K}=\left[\mathbf{k}_1,\mathbf{k}_2,\cdots,\mathbf{k}_{n_r},\mathbf{k}_{n_r+1},\mathbf{k}_{n_r+2},\cdots,\mathbf{k}_{2n_r},\right.
\cdots,\left.\mathbf{k}_{n-n_r+1},\mathbf{k}_{n-n_r+2},\cdots,\mathbf{k}_n\right],
\end{aligned}\end{equation*}
\begin{equation}\begin{aligned}\label{x_TV}
\mathbf{Z}_{\sigma\lambda_{TV}}(\sigma \mathbf{C}_1):=\mathop {\rm argmin }\limits_{\mathbf{Z}}\left\{\sigma\lambda_{TV}{\| {{\mathcal{\hat{H}}}_{v}\mathbf{Z}} \|_1}+\frac{1}{2}\|\mathbf{Z}-\sigma \mathbf{C}_1\|_F^2 \right\},
\end{aligned}\end{equation}
where $\mathbf{K}=\sigma \mathbf{C}_1$ and $\mathbf{z}_i,\mathbf{k}_i$ $(i=1,2,\cdots,n)$ denote the $i$th columns of $\mathbf{Z},\mathbf{K}$.
In order to compute $\mathbf{Z}_{\sigma\lambda_{TV}}(\sigma \mathbf{C}_1)$, we first denote
\begin{equation*}\begin{aligned}
{\mathbf{Z}'}^{*}:=\mathop {\rm argmin }\limits_{\mathbf{Z}'}\Big\{\sigma\lambda_{TV}\sum\limits_{k=1}\limits^{n_c-1}\|\mathbf{Z}'(:,k+1)-\mathbf{Z}'(:,k)\|_{1}+
\frac{1}{2}\|\mathbf{Z}'-\mathbf{K}'\|_F^2 \Big\},
\end{aligned}\end{equation*}
where
\begin{equation*}
\mathbf{Z}': = \left( {\begin{array}{*{20}{c}}
{{\mathbf{z}_1}}&{{\mathbf{z}_{n_r  + 1}}}& \cdots &{{\mathbf{z}_{n - n_r  + 1}}}\\
{{\mathbf{z}_2}}&{{\mathbf{z}_{n_r  + 2}}}& \cdots &{{\mathbf{z}_{n - n_r  + 2}}}\\
 \vdots & \vdots & \ddots & \vdots \\
{{\mathbf{z}_{n_r }}}&{{\mathbf{z}_{2n_r }}}& \cdots &{{\mathbf{z}_n}}
\end{array}} \right),\
\mathbf{K}': = \left( {\begin{array}{*{20}{c}}
{{\mathbf{k}_1}}&{{\mathbf{k}_{n_r  + 1}}}& \cdots &{{\mathbf{k}_{n - n_r  + 1}}}\\
{{\mathbf{k}_2}}&{{\mathbf{k}_{n_r  + 2}}}& \cdots &{{\mathbf{k}_{n - n_r  + 2}}}\\
 \vdots & \vdots & \ddots & \vdots \\
{{\mathbf{k}_{n_r }}}&{{\mathbf{k}_{2n_r }}}& \cdots &{{\mathbf{k}_n}}
\end{array}} \right),
\end{equation*}

\begin{equation*}
{\mathbf{Z}'}^*: = \left( {\begin{array}{*{20}{c}}
{{\mathbf{z}_1^*}}&{{\mathbf{z}_{n_r  + 1}^*}}& \cdots &{{\mathbf{z}_{n - n_r  + 1}^*}}\\
{{\mathbf{z}_2^*}}&{{\mathbf{z}_{n_r  + 2}^*}}& \cdots &{{\mathbf{z}_{n - n_r  + 2}^*}}\\
 \vdots & \vdots & \ddots & \vdots \\
{{\mathbf{z}_{n_r}^*}}&{{\mathbf{z}_{2n_r }^*}}& \cdots &{{\mathbf{z}_n^*}}
\end{array}} \right).
\end{equation*}
 Since ${\mathbf{Z}'}^{*}$ is separable, we can find ${\mathbf{Z}'}^{*}$ row by row. That is,
\begin{equation*}\begin{aligned}
{\mathbf{Z}'}^{*}(i,:)=&\mathop {\rm argmin }\limits_{{\mathbf{Z}'}(i,:)}\Big\{\sigma\lambda_{TV}\sum\limits_{k=1}\limits^{n_c-1}|{\mathbf{Z}'}(i,k+1)-{\mathbf{Z}'}(i,k)|+\frac{1}{2}\|{\mathbf{Z}'}(i,:)-{\mathbf{K}'}(i,:)\|_2^2 \Big\},\\& i=1,2,\cdots,m\times n_r.
\end{aligned}\end{equation*}

In our numerical experiments, we use the 1D TV Denoising Algorithm \cite{Condat2013} to compute ${\mathbf{Z}'}^{*}$.

Then we can obtain $\mathbf{Z}_{\sigma\lambda_{TV}}(\sigma \mathbf{C}_1)$ from ${\mathbf{Z}'}^{*}$ easily as below
\begin{equation*}\begin{aligned}
\mathbf{Z}_{\sigma\lambda_{TV}}(\sigma \mathbf{C}_1)=&[\mathbf{z}_1^*,\mathbf{z}_2^*,\cdots,\mathbf{z}_{n_r}^*,\mathbf{z}_{n_r+1}^*,\mathbf{z}_{n_r+2}^*,\cdots,\mathbf{z}_{2n_r}^*,
\cdots,\mathbf{z}_{n-n_r+1}^*,\mathbf{z}_{n-n_r+2}^*,\\&\cdots,\mathbf{z}_n^*].
\end{aligned}\end{equation*}

By Proposition 1, we can compute the proximal mapping of $\sigma {p_1}$ by composing the proximal mapping of $\delta_{R_{+}^{m\times n}}(\cdot)$ with the proximal mapping of $\sigma\lambda_{TV}{\| {{\mathcal{\hat{H}}}_{v}(\cdot)} \|_1}$ as below
\begin{equation*}
{{\mathop{\rm Prox}\nolimits} _{\sigma {p_{1}}}}(\sigma \mathbf{C}_1) =\Pi_{R_{+}^{m\times n}} (\mathbf{Z}_{\sigma\lambda_{TV}}(\sigma \mathbf{C}_1)).
\end{equation*}

(ii) When $p(\cdot)=\lambda\|\cdot\|_{2,1}+\delta_{R_{+}^{m\times n}}(\cdot)+\lambda_{TV}{\| {{\mathcal{\hat{H}}}_{v}(\cdot)} \|_1}$.
Similar to Proposition 1, for any $\sigma >0$, it also holds that
\begin{equation*}
{{\mathop{\rm Prox}\nolimits} _{{\sigma p_{}}}} =  {{\mathop{\rm Prox}\nolimits} _{\sigma\lambda\|\cdot\|_{2,1}}}\circ\mathop{\rm Prox}\nolimits_{\delta_{R_{+}^{m\times n}}} \circ {{\mathop{\rm Prox}\nolimits} _{\sigma\lambda _{TV}\| {{\mathcal{\hat{H}}}_{v}(\cdot)} \|_1}}.
\end{equation*}
As its proof is nearly the same as that of Proposition 1, we omit it.

Let $p_{2}(\cdot)=\delta_{R_{+}^{m\times n}}(\cdot)+\lambda_{TV}{\| {{\mathcal{\hat{H}}}_{v}(\cdot)} \|_1}$.
The proximal mapping of $\sigma p$ takes the following form
\begin{equation*}
{{\mathop{\rm Prox}\nolimits} _{{\sigma p}}}(\sigma \mathbf{C}_1)= \mbox{diag}(\alpha_1,\alpha_{2},\ldots,\alpha_{m}) {{\mathop{\rm Prox}\nolimits} _{\sigma {p_{2}}}}(\sigma \mathbf{C}_1),
\end{equation*}
where
\begin{equation*}
\alpha_i = \left( {\frac{{\max \left\{ {{\| {{ {{\mathop{\rm Prox}\nolimits} _{\sigma {p_{2}}}}(\sigma \mathbf{C}_1)(i,:)}} \|}_2} - \sigma\lambda ,0\right\} }}{{\max \left\{ {{\| {{ {{\mathop{\rm Prox}\nolimits} _{\sigma {p_{2}}}}(\sigma \mathbf{C}_1)(i,:)}} \|}_2} - \sigma\lambda ,0\right\}  + \sigma\lambda }}} \right)
\end{equation*}
and
\begin{equation*}
{{\mathop{\rm Prox}\nolimits} _{\sigma {p_{2}}}}(\sigma \mathbf{C}_1) =\Pi_{R_{+}^{m\times n}} (\mathbf{Z}_{\sigma\lambda_{TV}}(\sigma \mathbf{C}_1)).
\end{equation*}

(iii) Note $q(\cdot)=\lambda_{TV}{\| {{\mathcal{\hat{H}}}_{h}(\cdot)} \|_1}$,
let
\begin{equation*}\begin{aligned}
{\mathbf{Z}''}^{*}:=\mathop {\rm argmin }\limits_{{\mathbf{Z}''}}\Big\{\sigma\lambda_{TV}\sum\limits_{k=1}\limits^{n_r-1}\|{\mathbf{Z}''}(:,k+1)-{\mathbf{Z}''}(:,k)\|_{1}+
\frac{1}{2}\|{\mathbf{Z}''}-{\mathbf{K}''}\|_F^2 \Big\},
\end{aligned}\end{equation*}
where
\begin{equation*}
{\mathbf{Z}''}: = \left( {\begin{array}{*{20}{c}}
{{\mathbf{z}_1}}&{{\mathbf{z}_2}}& \cdots &{{\mathbf{z}_{n_r }}}\\
{{\mathbf{z}_{n_r  + 1}}}&{{\mathbf{z}_{n_r  + 2}}}& \cdots &{{\mathbf{z}_{2n_r }}}\\
 \vdots & \vdots & \ddots & \vdots \\
{{\mathbf{z}_{n - n_r  + 1}}}&{{\mathbf{z}_{n - n_r  + 2}}}& \cdots &{{\mathbf{z}_n}}
\end{array}} \right),\
{\mathbf{K}''}: = \left( {\begin{array}{*{20}{c}}
{{\mathbf{k}_1}}&{{\mathbf{k}_2}}& \cdots &{{\mathbf{k}_{n_r }}}\\
{{\mathbf{k}_{n_r  + 1}}}&{{\mathbf{k}_{n_r  + 2}}}& \cdots &{{\mathbf{k}_{2n_r }}}\\
 \vdots & \vdots & \ddots & \vdots \\
{{\mathbf{k}_{n - n_r  + 1}}}&{{\mathbf{k}_{n - n_r  + 2}}}& \cdots &{{\mathbf{k}_n}}
\end{array}} \right),
\end{equation*}
\begin{equation*}
{\mathbf{Z}''}^{*}: = \left( {\begin{array}{*{20}{c}}
{{\mathbf{z}_1^*}}&{{\mathbf{z}_2^*}}& \cdots &{{\mathbf{z}_{n_r }^*}}\\
{{\mathbf{z}_{n_r  + 1}^*}}&{{\mathbf{z}_{n_r  + 2}^*}}& \cdots &{{\mathbf{z}_{2n_r }^*}}\\
 \vdots & \vdots & \ddots & \vdots \\
{{\mathbf{z}_{n - n_r  + 1}^*}}&{{\mathbf{z}_{n - n_r  + 2}^*}}& \cdots &{{\mathbf{z}_n^*}}
\end{array}} \right).
\end{equation*}
Then we can compute $\mathbf{Z}''^{*}$ similarly to that of $\mathbf{Z}'^{*}$, that is,
\begin{equation*}\begin{aligned}
{\mathbf{Z}''}^{*}(i,:)=&\mathop {\rm argmin }\limits_{{\mathbf{Z}''}(i,:)}\Big\{\sigma\lambda_{TV}\sum\limits_{k=1}\limits^{n_r-1}|{\mathbf{Z}''}(i,k+1)-{\mathbf{Z}''}(i,k)|+\frac{1}{2}\|{\mathbf{Z}''}(i,:)-{\mathbf{K}''}(i,:)\|_2^2 \Big\},\\& i=1,2,\cdots,m\times n_c.
\end{aligned}\end{equation*}
We can obtain the proximal mapping of $\sigma q$ from ${\mathbf{Z}''}^{*}$ as below
\begin{equation*}\begin{aligned}
{{\mathop{\rm Prox}\nolimits} _{{\sigma q}}}(\sigma \mathbf{C}_2)=&[\mathbf{z}_1^*,\mathbf{z}_2^*,\cdots,\mathbf{z}_{n_r}^*,\mathbf{z}_{n_r+1}^*,\mathbf{z}_{n_r+2}^*,\cdots,\mathbf{z}_{2n_r}^*,
\cdots,\mathbf{z}_{n-n_r+1}^*,\mathbf{z}_{n-n_r+2}^*,\\&\cdots,\mathbf{z}_n^*].
\end{aligned}\end{equation*}

As for the computational complexity of Algorithm 2, the most expensive step of Algorithm 2 is to compute $\mathbf{V}_3$ whose computational complexity is $O(nL\cdot \max\{m,L\})$. Thus, the computational complexity of Algorithm 2 in each iteration is $O(nL\cdot \max\{m,L\})$.

\subsection{Convergence analysis}
In this subsection, we analyze the global convergence and the local linear convergence rate of the dual sGS-ADMM.
Note that the problem (\ref{D2}) can also be written in the following form
\begin{equation}
\begin{aligned}
\label{compact}
\min_{{\mathbf{W}},{\mathbf{V}_2}}\ &\psi(\mathbf{W})+\varphi(\mathbf{V}_2)\\
\mbox{s.t.}\quad &-\mathcal{A}\mathbf{W}-\mathbf{V}_2=0,
\end{aligned}
\end{equation}
where $\mathbf{W}=\left(\mathbf{V}_1^T,\mathbf{V}_3^T\right)^T\in R^{(m+L)\times n}$, $\mathcal{A}=\left[\mathbf{I}_{m\times m},\mathbf{A}^T\right]$ and
\begin{eqnarray*}
	\psi(\mathbf{W})=p^{*}(-\mathbf{V}_1)+ \frac{1}{2}\| {{\mathbf{V}_3}} \|_F^2 - \left\langle {{\mathbf{V}_3},\mathbf{Y}} \right\rangle,\
	\varphi(\mathbf{V}_2)= {q^*}( - {\mathbf{V}_2}).
\end{eqnarray*}
Suppose that $(\mathbf{W},\mathbf{V}_2)\in R^{(m+L)\times n}\times R^{m\times n}$ is an optimal solution to the problem (\ref{compact}). If there exists $\mathbf{X}\in R^{m\times n}$ such that $(\mathbf{W},\mathbf{V}_2,\mathbf{X})$ satisfies the following KKT system
\begin{equation}\label{compactKKT}
\left\{ {\begin{array}{*{20}{l}}
	{0 \in \partial \psi(\mathbf{W}) + \mathcal{A}^*\mathbf{X}},\\
	{0 \in \partial \varphi(\mathbf{V}_2) + \mathbf{X}},\\
	{-\mathcal{A}\mathbf{W}-\mathbf{V}_2=0},
	\end{array}} \right.
\end{equation}
then $(\mathbf{W},\mathbf{V}_2,\mathbf{X})$ is a KKT point for the problem (\ref{compact}).
Let $\overline{\Omega}$ be the solution set of the KKT system (\ref{compactKKT}) for convenience.
Let $\mathcal{E}:R^{m\times n}\rightarrow R^{(m+L)\times n}\times R^{m\times n}$ be a linear operator such that its adjoint $\mathcal{E}^{*}(\mathbf{W},\mathbf{V}_2)=-\mathcal{A}\mathbf{W}-\mathbf{V}_2$.
	For $\mathbf{v}:=(\mathbf{W},\mathbf{V}_2,\mathbf{X})$, the KKT mapping is defined by
	\begin{equation*}
	\mathbf{R}(\mathbf{v}): = \left( {\begin{array}{*{20}{c}}
		{{\mathbf{W}} - {{{\mathop{\rm Prox}\nolimits} }_{\psi}}(\mathbf{W}-\mathcal{A}^*\mathbf{X})}\\
		{{\mathbf{V}_2} - {{{\mathop{\rm Prox}\nolimits} }_\varphi}({\mathbf{X}} - {\mathbf{V}_2})}\\
		-\mathcal{A}\mathbf{W}-\mathbf{V}_2
		\end{array}} \right).
	\end{equation*}

Since the subdifferential mappings of the closed proper convex functions $ p^{*}$ and $ q^{*}$ are maximally monotone \cite[Theorem 12.17]{Rockafellar1998}, there exist two self-adjoint and positive semidefinite linear operators $\Sigma_{\psi}$ and $\Sigma_{\varphi}$ such that for all $\mathbf{y},\mathbf{y}'\in \textrm{dom}(\psi)$, $\xi\in\partial\psi(\mathbf{y})$ and $\xi'\in\partial\psi(\mathbf{y}')$
\begin{equation*}
\psi(\mathbf{y})\geq\psi(\mathbf{y}')+\langle\xi',\mathbf{y}-\mathbf{y}'\rangle+\frac{1}{2}\|\mathbf{y}-\mathbf{y}'\|_{\Sigma_\psi}^2,\
\langle\xi-\xi',\mathbf{y}-\mathbf{y}'\rangle\geq\|\mathbf{y}-\mathbf{y}'\|_{\Sigma_\psi}^2,
\end{equation*}
and for all $\mathbf{z},\mathbf{z}'\in \textrm{dom}(\varphi)$, $\zeta\in\partial\varphi(\mathbf{z})$ and $\zeta'\in\partial\varphi(\mathbf{z}')$
\begin{equation*}
\varphi(\mathbf{z})\geq\varphi(\mathbf{z}')+\langle\zeta',\mathbf{z}-\mathbf{z}'\rangle+\frac{1}{2}\|\mathbf{z}-\mathbf{z}'\|_{\Sigma_\varphi}^2,\
\langle\zeta-\zeta',\mathbf{z}-\mathbf{z}'\rangle\geq\|\mathbf{z}-\mathbf{z}'\|_{\Sigma_\varphi}^2.
\end{equation*}

For the convergence of the dual sGS-ADMM, we need the following assumption.

\emph{Assumption 2}: The KKT system (\ref{compactKKT}) has a non-empty solution set.\label{Assumption: KKT solution}

We can also obtain the global convergence and linear convergence rate of the sGS-ADMM based on \cite[Theorem B.1]{FazelPST} and \cite[Theorem 2]{linear-rate}.
\begin{theorem}
		Suppose that Assumption 2 holds and
		$\mathbf{S}: = (\sigma^{-1}\mathbf{I}+\mathbf{A}\mathbf{A}^T)-\mathbf{A}[\mathbf{I}+\mathbf{A}^T(\sigma^{-1}\mathbf{I}+\mathbf{A}\mathbf{A}^T)^{-1}\mathbf{A}]^{-1}\mathbf{A}^T$ is positive definite.
		Let $\tau\in(0,\frac{1+\sqrt{5}}{2})$,
		then there exists a KKT point $\bar{\mathbf{v}}:=(\overline{\mathbf{W}},\overline{\mathbf{V}}_2, \overline{\mathbf{X}})\in\overline{\Omega}$ such that the sequence $\{(\mathbf{V}_1^k,\mathbf{V}_2^k,\mathbf{V}_3^k,\mathbf{X}^k)\}$ generated by Algorithm 2 converges to $\bar{\mathbf{v}}$. Assume that $\mathbf{R}^{-1}$ is calm at the origin for $\bar{\mathbf{v}}$ with modulus $\eta'>0$, i.e., there exists $r'>0$ such that
		\begin{equation*}
		\textrm{dist}(\mathbf{u},\overline{\Omega})\leq\eta'\|\mathbf{R}(\mathbf{u})\|,\quad \forall\, \mathbf{u}\in\{\mathbf{u}:\|\mathbf{u}-\bar{\mathbf{u}}\|\leq r'\}.
		\end{equation*}
		Then there exists an integer $\bar{k}\geq 1$ such that for all $k\geq \bar{k}$,
		\begin{equation*}\begin{aligned}
		&\textrm{dist}_{\mathcal{P}}^{2}(\mathbf{u}^{k+1},\overline{\Omega}) \leq \mu' \textrm{dist}_{\mathcal{P}}^{2}(\mathbf{u}^{k},\overline{\Omega}),
		\end{aligned}\end{equation*}
		where $\mu'\in(0,1)$,
		\begin{equation*}
\mathcal{P}:=\textrm{Diag}(\mathcal{T}+\Sigma_\psi,\Sigma_\varphi+\sigma\mathcal{I},(\tau\sigma)^{-1}\mathcal{I})+s_\tau\sigma\mathcal{E}\mathcal{E}^{*},
		\end{equation*}
		and
		\begin{equation*}
		s_\tau:=\frac{5-\tau-3 \min\{\tau,\tau^{-1}\}}{4}.
		\end{equation*}
		Moreover, there exists a positive number $\varsigma'\in[\mu',1)$ such that for all $k\geq 1$
		\begin{equation*}\begin{aligned}
		&\textrm{dist}_{\mathcal{P}}^{2}(\mathbf{u}^{k+1},\overline{\Omega}) \leq \varsigma'  \textrm{dist}_{\mathcal{P}}^{2}(\mathbf{u}^{k},\overline{\Omega}).
		\end{aligned}\end{equation*}
\end{theorem}

\emph{Proof}: As pointed out in \cite{LiST2019}, implementing one cycle of the sGS method is equivalent to solve the associated convex quadratic programming problem plus an extra semiproximal term $\frac{1}{2}\|\mathbf{W}-\mathbf{W}^k\|_{\mathcal{T}}^2$,
where $\mathcal{T}$ is a symmetric positive semidefinite operator related to the sGS decomposition and $\mathbf{W}^{k}$ is the previous iterate.
This indicates that the dual sGS-ADMM is essentially a special case of the SPADMM.
According to \cite[Theorem B.1]{FazelPST}, in order to prove the convergence, we only need to prove that the following two conditions hold.
\begin{equation}
\label{cond:psd4}
\Sigma_\varphi+\sigma \mathbf{I}_{m\times m}\succ 0,
\end{equation}
\begin{equation}
\label{cond:psd3}
\Sigma_\psi +\sigma \mathcal{A}^*\mathcal{A}+\mathcal{T}\succ 0,
\end{equation}
where
\begin{equation*}
\Sigma_\psi: =\left(
                \begin{array}{cc}
                  0_{m\times m} & \ 0_{m\times L} \\
                  0_{L\times m} & \ \mathbf{I}_{L\times L} \\
                \end{array}
              \right),
\end{equation*}
\begin{equation*}\begin{aligned}
\mathcal{T}:&= \left(
                                                           \begin{array}{cc}
                                                           0_{m\times m} & \sigma \mathbf{A}^T \\
                                                            0_{L\times m} & 0_{L\times L} \\
                                                           \end{array}
                                                         \right)\left(
                                                                  \begin{array}{cc}
                                                                    \sigma^{-1}\mathbf{I}_{m\times m} & 0_{m\times L} \\
                                                                    0_{L\times m} & (\mathbf{I}_{L\times L}+\sigma \mathbf{A}\mathbf{A}^T)^{-1} \\
                                                                  \end{array}
                                                                \right)\left(
                                                                         \begin{array}{cc}
                                                                           0_{m\times m} &0_{m\times L} \\
                                                                           \sigma \mathbf{A} & 0_{L\times L} \\
                                                                         \end{array}
                                                                       \right)
\\&=\left(
               \begin{array}{cc}
                 \sigma \mathbf{A}^T(\sigma^{-1}\mathbf{I}_{L\times L}+\mathbf{A}\mathbf{A}^T)^{-1}\mathbf{A} & 0_{m\times L} \\
                 0_{L\times m} & 0_{L\times L} \\
               \end{array}
             \right).
             \end{aligned}\end{equation*}
It is  obvious that \eqref{cond:psd4} holds automatically. By simple calculations, we know that \eqref{cond:psd3} is equivalent to that
\begin{equation*}\begin{aligned}
\left(
  \begin{array}{cc}
    \mathbf{I}+\mathbf{A}^T(\sigma^{-1}\mathbf{I}_{L\times L}+\mathbf{A}\mathbf{A}^T)^{-1}\mathbf{A} & \mathbf{A}^T \\
    \mathbf{A} & \sigma^{-1}\mathbf{I}_{L\times L}+\mathbf{A}\mathbf{A}^T \\
  \end{array}
\right)
\end{aligned}\end{equation*}
is positive definite. From the Schur complement condition \cite{Horn1985}, we only need to require that $\mathbf{S} = (\sigma^{-1}\mathbf{I}+\mathbf{A}\mathbf{A}^T)-\mathbf{A}[\mathbf{I}+\mathbf{A}^T(\sigma^{-1}\mathbf{I}+\mathbf{A}\mathbf{A}^T)^{-1}\mathbf{A}]^{-1}\mathbf{A}^T$ is positive definite.

\emph{Remark 2}: For the case of $\rho=1$ in the problem (\ref{model1}), $p$ and $q$ are piecewise linear-quadratic functions.
We also know that the calmness condition holds automatically by \cite[Proposition 2.24]{Sun1986} and \cite[Corollary]{Robinson1981}.

\section{Numerical experiments}
In this section, we implement some numerical experiments to demonstrate the efficiency of our algorithm. All the experiments were conducted on a PC with Inter (R) Core (TM) i5-4210M CPU @2.60GHz 2.59GHz of 8G memory running 64bit Windows operation system.
The Monte Carlo simulations were used ten times in all our experiments.
All the codes were written in {\sc Matlab} 2017b with some subroutines in C.
For convenience, we use SUnTV-sGSADMM and CLSUnTV-sGSADMM to denote the dual sGS-ADMM applied to the SUnTV and CLSUnTV, respectively.
The {\sc Matlab} codes of SUnSAL-TV and CLSUnSAL-TV are based on \cite{SUnSAL-TV}.\footnote{downloaded from http://www.lx.it.pt/$\sim$bioucas/publications.html.}

We use the signal reconstruction error (SRE) measured in decibel (dB) to measure the performances of different algorithms, which is defined as follows \cite{SUnSAL}:
\begin{equation*}
\textrm{SRE}(\textrm{dB})=10{\log _{10}}\frac{{E[\|\mathbf{x}\|_2^2]}}{{E[\|\mathbf{x} - \hat{\mathbf{x}} \|_2^2]}},
\end{equation*}
where $\mathbf{x}$ and $\hat{\mathbf{x}}$ denote the true abundances and the estimated abundances, respectively and $E[ \cdot ]$ represents the statistical expectation.
In general, a larger SRE indicates a better unmixing performance.
In addition, the probability of success $(p_s)$ is also employed to estimate the probability that the relative error power is less than a certain threshold value \cite{SUnSAL}:
\begin{equation*}
{p_s} = P(||\hat {\mathbf{x}} - \mathbf{x}|{|^2}/||\mathbf{x}{\rm{|}}{{\rm{|}}^2} \le {\rm{threshold}}).
\end{equation*}
The unmixing algorithm can be regarded as being successful when $ ||\hat {\mathbf{x}} - \mathbf{x}|{|^2}/||\mathbf{x}{\rm{|}}{{\rm{|}}^2} \le 0.316(5\ \textrm{dB}) $ in all our experiments.

We measure the accuracy of the solutions by the following relative KKT residual and Error:
\begin{eqnarray*}
&&\mathbf{R}_{P_1}=(\|\mathbf{D}_1-\mathbf{A\widetilde{D}}\|_F+\|\mathbf{D}_2-\mathbf{\widetilde{D}}\|_F+\|\mathbf{D}_3-\mathbf{\widetilde{D}}\|_F+\|\mathbf{D}_4-\mathcal{H}\mathbf{D}_3\|_F+\\&&\quad\quad\quad\|\mathbf{D}_5-\mathbf{\widetilde{D}}\|_F)/(1+\|\mathbf{A}\|_F),\\
&&\mathbf{R}_{D_1}=(\|\mathbf{A}^{T}\mathbf{\Lambda}_1+\mathbf{\Lambda}_2+\mathbf{\Lambda}_3+\mathbf{\Lambda}_5\|_F+ \|\mathbf{\Lambda}_3-\mathcal{H}^T\mathbf{\Lambda}_4\|_F)/(1+\|\mathbf{A}\|_F),\\
&&\textrm{Error}_1=\|\mathbf{\widetilde{D}}^{k+1}-\mathbf{\widetilde{D}}^{k}\|_{F}/\|\mathbf{\widetilde{D}}^{k+1}\|_{F},\\
&&\mathbf{R}_{P_2}=\|\mathbf{AX}-\mathbf{Y}-\mathbf{U}_3\|_F/(1+\|\mathbf{Y}\|_F),\\
&&\mathbf{R}_{D_2}=\|\mathbf{V}_1+\mathbf{V}_2+\mathbf{A}^{T}\mathbf{V}_3\|_F/(1+\|\mathbf{A}\|_F),\\
&&\textrm{Error}_2=\|\mathbf{X}^{k+1}-\mathbf{X}^{k}\|_{F}/\|\mathbf{X}^{k+1}\|_{F}.
\end{eqnarray*}
The stopping criterion for Algorithm 1 is $\mathbf{R}_{P_1}<\textrm{tol}_1$ and $\mathbf{R}_{D_1}<\textrm{tol}_1$ or $\textrm{Error}_1<\textrm{tol}_2$, and the stopping criterion for Algorithm 2 is $\mathbf{R}_{P_2}<\textrm{tol}_1$ and $\mathbf{R}_{D_2}<\textrm{tol}_1$ or $\textrm{Error}_2<\textrm{tol}_2$, where $\textrm{tol}_1$ and $\textrm{tol}_2$ are predefined error tolerances. In our implementation, we empirically set $\textrm{tol}_1=10^{-3}$ and $\textrm{tol}_2=10^{-4}$.
The maximum number of iterations of the SUnSAL-TV and CLSUnSAL-TV are capped by $200$, and the maximum number of iterations of the SUnTV-sGSADMM and CLSUnTV-sGSADMM are capped by $50$.\footnote{Since the SRE value will get worse after the KKT residual is smaller than some value, and Algorithm 2 is obviously faster than Algorithm 1 according to the experiments behind, the settings of the maximal iterations are not the same.}
The optimal regularization parameters $\lambda$ and $\lambda_{TV}$ in all compared algorithms were selected from all possible combinations of the following finite set $\{0.5,0.1,0.05,0.01,0.005,$\\$0.001,0.0005,0.0001,0.00005,0.00001\}$ in order to produce the highest SRE.

\subsection{Numerical results for the simulated data}
For the simulated data experiments, the spectral library $\mathbf{A}\in R^{224\times 240}$ is randomly picked out from the United States Geological Survey (USGS) digital spectral library splib06.\footnote{Available online: http://speclab.cr.usgs.gov/spectral.lib06.}
These spectra have 224 bands and are uniformly distributed between 0.4-2.5 um.
The mutual coherence of $\mathbf{A}$ is very close to 1.
By using this spectral library, we generated the following three simulated hyperspectral data cubes.

\begin{table*}[!tbp]
\scriptsize
\caption{The numerical experiments on the simulated data}
\setlength{\tabcolsep}{0mm}{
\centering
\label{tabel1}
\begin{tabular}{|c|c|c|c|c|c|c|}
\hline
\multicolumn{7}{|c|}{White noise}                                                                                                                                                                                                          \\ \hline
\multicolumn{1}{|c|}{Data cube} & \multicolumn{1}{c|}{SNR(dB)} &\multicolumn{1}{c|}{Parameters} & SUnSAL-TV & SUnTV-sGSADMM & CLSUnSAL-TV & CLSUnTV-sGSADMM  \\ \hline
\multirow{15}{*}{DC1}
                                & \multicolumn{1}{c|}{\multirow{5}{*}{20}} & \multicolumn{1}{c|}{SRE(dB)}  & 7.1449 (0.0216)    & 12.2804 (0.2820)   & 6.9773 (0.0384)   & 11.3766 (0.2387)  \\ \cline{3-7}
                                & \multicolumn{1}{c|}{}                    & \multicolumn{1}{c|}{$p_s$}    & 0.9561 (0.0046)    & 0.9789 (0.0037)  & 0.9505 (0.0015) & 0.9796 (0.0027) \\ \cline{3-7}
                                & \multicolumn{1}{c|}{}                    & \multicolumn{1}{c|}{time(s)}  & 134.8394 (2.1428)   & 15.1354 (0.2630)    & 138.4100 (1.7137)  & 15.4247 (0.1992)   \\ \cline{3-7}
                                & \multicolumn{1}{c|}{}                    & \multicolumn{1}{c|}{$\lambda$}& 0.05    & 0.005   & 0.5    & 0.5    \\ \cline{3-7}
                                & \multicolumn{1}{c|}{}               & \multicolumn{1}{c|}{$\lambda_{TV}$}& 0.05    & 0.1     & 0.05   & 0.1   \\ \cline{2-7}

                                & \multicolumn{1}{c|}{\multirow{5}{*}{30}} & \multicolumn{1}{c|}{SRE(dB)}  & 15.0541 (0.1364)   & 17.5032 (0.2956)   & 14.1178 (0.1306)  & 16.4817 (0.1534)  \\ \cline{3-7}
                                & \multicolumn{1}{c|}{}                    & \multicolumn{1}{c|}{$p_s$}    & 0.9956 (0)  & 0.9956 (0)  & 0.9959 (0.0009) & 0.9996 (0.0013)      \\ \cline{3-7}
                                & \multicolumn{1}{c|}{}                    & \multicolumn{1}{c|}{time(s)}  & 132.7691 (1.0134)   & 15.6431 (0.2305)    & 136.9076 (1.1648)  & 15.6762 (0.1578)  \\ \cline{3-7}
                                & \multicolumn{1}{c|}{}                  & \multicolumn{1}{c|}{$\lambda$}  & 0.005   & 0.001   & 0.5    & 0.1    \\ \cline{3-7}
                                & \multicolumn{1}{c|}{}             & \multicolumn{1}{c|}{$\lambda_{TV}$}  & 0.01    & 0.01    & 0.01   & 0.01   \\ \cline{2-7}

                                & \multicolumn{1}{c|}{\multirow{5}{*}{40}} & \multicolumn{1}{c|}{SRE(dB)}  & 22.4525 (0.1267)   & 23.2758 (0.0850)   & 22.7363 (0.0993)  &23.0622 (0.0606)  \\ \cline{3-7}
                                & \multicolumn{1}{c|}{}                    & \multicolumn{1}{c|}{$p_s$}    & 1.0000 (0)      & 1.0000 (0)       & 1.0000 (0)      & 1.0000 (0)      \\ \cline{3-7}
                                & \multicolumn{1}{c|}{}                    & \multicolumn{1}{c|}{time(s)}  & 132.1023 (0.7609)   & 15.3504 (0.1765)    & 136.3907 (0.5466)  & 15.8081 (0.2889)   \\ \cline{3-7}
                                & \multicolumn{1}{c|}{}                    & \multicolumn{1}{c|}{$\lambda$}& 0.001   & 0.001   & 0.1    & 0.1    \\ \cline{3-7}
                                & \multicolumn{1}{c|}{}               & \multicolumn{1}{c|}{$\lambda_{TV}$}& 0.005   & 0.005   & 0.005  & 0.005  \\ \cline{1-7}
\multicolumn{7}{|l|}{}                                                                                                                                                                                                                     \\ \hline
\multirow{15}{*}{DC2}           & \multicolumn{1}{c|}{\multirow{5}{*}{20}} & \multicolumn{1}{c|}{SRE(dB)}  & 8.7308 (0.0700)    & 8.8712 (0.1462)    & 6.9100 (0.0700)   & 7.1687 (0.0666)   \\ \cline{3-7}
                                & \multicolumn{1}{c|}{}                    & \multicolumn{1}{c|}{$p_s$}    & 0.8573 (0.0118)  & 0.8861 (0.0166)  & 0.7661 (0.0155) & 0.7789 (0.0082)  \\ \cline{3-7}
                                & \multicolumn{1}{c|}{}                    & \multicolumn{1}{c|}{time(s)}  & 196.5701 (4.8571)   & 27.6584 (0.2723)    & 197.3205 (1.3822)  & 27.4746 (0.1534)   \\ \cline{3-7}
                                & \multicolumn{1}{c|}{}                    & \multicolumn{1}{c|}{$\lambda$}& 0.01    & 0.01    & 0.005 & 0.001  \\ \cline{3-7}
                                & \multicolumn{1}{c|}{}               & \multicolumn{1}{c|}{$\lambda_{TV}$}& 0.01    & 0.01    & 0.05   & 0.05   \\ \cline{2-7}

                                & \multicolumn{1}{c|}{\multirow{5}{*}{30}} & \multicolumn{1}{c|}{SRE(dB)}  & 15.9441 (0.0544)   & 16.2301 (0.0401)   & 14.3376 (0.0560)  & 14.3426 (0.0861)  \\ \cline{3-7}
                                & \multicolumn{1}{c|}{}                    & \multicolumn{1}{c|}{$p_s$}    & 0.9967 (0.0003)  & 0.9945 (0.0006)  & 0.9953 (0.0007) & 0.9918 (0.0014)  \\ \cline{3-7}
                                & \multicolumn{1}{c|}{}                    & \multicolumn{1}{c|}{time(s)}  & 195.4667 (1.2963)   & 27.2704 (0.1196)    & 197.3390 (0.6061)  & 28.0106 (0.0646)   \\ \cline{3-7}
                                & \multicolumn{1}{c|}{}                    & \multicolumn{1}{c|}{$\lambda$}& 0.005   & 0.005   & 0.005  & 0.0001  \\ \cline{3-7}
                                & \multicolumn{1}{c|}{}               & \multicolumn{1}{c|}{$\lambda_{TV}$}& 0.005   & 0.005   & 0.005  & 0.005  \\ \cline{2-7}

                                & \multicolumn{1}{c|}{\multirow{5}{*}{40}} & \multicolumn{1}{c|}{SRE(dB)}  & 21.0885 (0.0471)   & 21.6292 (0.0722)   & 19.7971 (0.0458)  & 19.7908 (0.0677)   \\ \cline{3-7}
                                & \multicolumn{1}{c|}{}                    & \multicolumn{1}{c|}{$p_s$}    & 1.0000 (0)  & 0.9997 (0.0001)  & 1.0000 (0) & 0.9997 (0.0001)  \\ \cline{3-7}
                                & \multicolumn{1}{c|}{}                    & \multicolumn{1}{c|}{time(s)}  & 195.4878 (0.8895)   & 27.5456 (0.0784)    & 197.4681 (0.2657)  & 28.2979 (0.1381)   \\ \cline{3-7}
                                & \multicolumn{1}{c|}{}                    & \multicolumn{1}{c|}{$\lambda$}& 0.001   & 0.0005  & 0.01   & 0.005  \\ \cline{3-7}
                                & \multicolumn{1}{c|}{}               & \multicolumn{1}{c|}{$\lambda_{TV}$}& 0.001   & 0.001   & 0.001  & 0.001  \\ \cline{1-7}
                                \multicolumn{7}{|l|}{}                                                                                                                                                                                                                     \\ \hline
\multirow{15}{*}{DC3}           & \multicolumn{1}{c|}{\multirow{5}{*}{20}} & \multicolumn{1}{c|}{SRE(dB)}  & 5.7547 (0.0390)    & 6.1169 (0.0357)    &5.1223 (0.0304)   & 5.1282 (0.0455)   \\ \cline{3-7}
                                & \multicolumn{1}{c|}{}                    & \multicolumn{1}{c|}{$p_s$}    & 0.6235 (0.0101)  & 0.6575 (0.0080)  & 0.5947 (0.0023)  & 0.5839 (0.0055)   \\ \cline{3-7}
                                & \multicolumn{1}{c|}{}                    & \multicolumn{1}{c|}{time(s)}  & 194.7100 (3.0791)   & 27.5685 (0.0489)    & 196.6756 (2.1201)  & 28.5397 (1.7651)   \\ \cline{3-7}
                                & \multicolumn{1}{c|}{}                    & \multicolumn{1}{c|}{$\lambda$}& 0.01    & 0.01    & 0.001  & 0.001  \\ \cline{3-7}
                                & \multicolumn{1}{c|}{}               & \multicolumn{1}{c|}{$\lambda_{TV}$}& 0.01    & 0.01    & 0.05   & 0.05   \\ \cline{2-7}

                                & \multicolumn{1}{c|}{\multirow{5}{*}{30}} & \multicolumn{1}{c|}{SRE(dB)}  & 11.5701 (0.0540)   & 11.6032 (0.0338)   & 9.6904 (0.0426)   & 9.6343 (0.0556)  \\ \cline{3-7}
                                & \multicolumn{1}{c|}{}                    & \multicolumn{1}{c|}{$p_s$}    & 0.9517 (0.0013)   & 0.9444 (0.0020)  & 0.9025 (0.0023) & 0.8872 (0.0036)  \\ \cline{3-7}
                                & \multicolumn{1}{c|}{}                    & \multicolumn{1}{c|}{time(s)}  & 193.4672 (0.4695)   & 27.5775 (0.0398)    & 196.5886 (1.5875)    & 28.0825 (0.1090)  \\ \cline{3-7}
                                & \multicolumn{1}{c|}{}                    & \multicolumn{1}{c|}{$\lambda$}& 0.01    & 0.005   & 0.05   & 0.0001  \\ \cline{3-7}
                                & \multicolumn{1}{c|}{}               & \multicolumn{1}{c|}{$\lambda_{TV}$}& 0.005   & 0.005   & 0.005  & 0.005  \\ \cline{2-7}

                                & \multicolumn{1}{c|}{\multirow{5}{*}{40}} & \multicolumn{1}{c|}{SRE(dB)}  & 18.0549 (0.0288)   & 17.9194 (0.0531)    & 15.9123 (0.0381)   & 15.9489 (0.0526)   \\ \cline{3-7}
                                & \multicolumn{1}{c|}{}                    & \multicolumn{1}{c|}{$p_s$}    & 1.0000 (0)  & 0.9997 (0.0001)  & 0.9983 (0.0003) & 0.9957 (0.0006)  \\ \cline{3-7}
                                & \multicolumn{1}{c|}{}                    & \multicolumn{1}{c|}{time(s)}  & 194.3947 (0.7461)   & 27.8315 (0.0969)    & 196.9297 (0.9201)  & 28.3006 (0.1138)   \\ \cline{3-7}
                                & \multicolumn{1}{c|}{}                    & \multicolumn{1}{c|}{$\lambda$}& 0.005   & 0.001   & 0.005  & 0.0001 \\ \cline{3-7}
                                & \multicolumn{1}{c|}{}               & \multicolumn{1}{c|}{$\lambda_{TV}$}& 0.0005  & 0.0005  & 0.001  & 0.001  \\ \cline{1-7}
\end{tabular}}
\end{table*}

\begin{table*}[!htbp]
\scriptsize
\caption{The numerical experiments on the simulated data}
\setlength{\tabcolsep}{0mm}{
\centering
\label{tabel2}
\begin{tabular}{|c|c|c|c|c|c|c|}
\hline
\multicolumn{7}{|c|}{Correlated noise}                                                                                                                                                                                                          \\ \hline
\multicolumn{1}{|c|}{Data cube} & \multicolumn{1}{c|}{SNR(dB)} & \multicolumn{1}{c|}{Parameters} & SUnSAL-TV & SUnTV-sGSADMM & CLSUnSAL-TV & CLSUnTV-sGSADMM  \\ \hline
\multirow{15}{*}{DC1}
                                & \multicolumn{1}{c|}{\multirow{5}{*}{20}} & \multicolumn{1}{c|}{SRE(dB)}  & 12.1682 (0.0783)   & 18.1366 (0.1336)   & 12.1838 (0.0655)  & 18.0484 (0.1029)  \\ \cline{3-7}
                                & \multicolumn{1}{c|}{}                    & \multicolumn{1}{c|}{$p_s$}    & 0.9961 (0.0012)  & 0.9993 (0.0006)  & 0.9990 (0.0004) & 1.0000 (0)      \\ \cline{3-7}
                                & \multicolumn{1}{c|}{}                    & \multicolumn{1}{c|}{time(s)}  & 133.0032 (1.1713)   & 14.4636 (0.2023)    & 136.3015 (0.2182)  & 14.7316 (0.2039)   \\ \cline{3-7}
                                & \multicolumn{1}{c|}{}                    & \multicolumn{1}{c|}{$\lambda$}& 0.005   & 0.001   & 0.1    & 0.1    \\ \cline{3-7}
                                & \multicolumn{1}{c|}{}               & \multicolumn{1}{c|}{$\lambda_{TV}$}& 0.0001  & 0.01    & 0.0001 & 0.01 \\ \cline{2-7}

                                & \multicolumn{1}{c|}{\multirow{5}{*}{30}} & \multicolumn{1}{c|}{SRE(dB)}  & 20.9386 (0.1049)   & 23.6159 (0.1003)   & 20.8904 (0.0936)  & 22.6764 (0.0920) \\ \cline{3-7}
                                & \multicolumn{1}{c|}{}                    & \multicolumn{1}{c|}{$p_s$}    & 1.0000 (0)       & 1.0000 (0)       &  1.0000 (0)      & 1.0000 (0)   \\ \cline{3-7}
                                & \multicolumn{1}{c|}{}                    & \multicolumn{1}{c|}{time(s)}  & 131.8280 (0.4336)     & 15.1003 (0.2726)    & 136.3561 (0.1915)  & 16.0053 (0.1282)   \\ \cline{3-7}
                                & \multicolumn{1}{c|}{}                    & \multicolumn{1}{c|}{$\lambda$}& 0.001   & 0.0005  & 0.1    & 0.05    \\ \cline{3-7}
                                & \multicolumn{1}{c|}{}               & \multicolumn{1}{c|}{$\lambda_{TV}$}& 0.0001  & 0.001   & 0.0001 & 0.0001 \\ \cline{2-7}

                                & \multicolumn{1}{c|}{\multirow{5}{*}{40}} & \multicolumn{1}{c|}{SRE(dB)}  & 30.6734 (0.0652)   & 33.7144 (0.1447)   & 30.5253 (0.1075)  & 31.7891 (0.0572)  \\ \cline{3-7}
                                & \multicolumn{1}{c|}{}                    & \multicolumn{1}{c|}{$p_s$}    & 1.0000 (0)       & 1.0000 (0)       & 1.0000 (0)      & 1.0000 (0)   \\ \cline{3-7}
                                & \multicolumn{1}{c|}{}                    & \multicolumn{1}{c|}{time(s)}  & 131.9362 (0.4130)   & 14.9363 (0.1920)    & 136.6289 (0.8981)  & 15.3723 (0.2121)   \\ \cline{3-7}
                                & \multicolumn{1}{c|}{}                    & \multicolumn{1}{c|}{$\lambda$}& 0.0005  & 0.0001  & 0.01   & 0.01   \\ \cline{3-7}
                                & \multicolumn{1}{c|}{}               & \multicolumn{1}{c|}{$\lambda_{TV}$}& 0.0001  & 0.0005  & 0.0001 & 0.0005 \\ \cline{1-7}
\multicolumn{7}{|l|}{}                                                                                                                                                                                                                     \\ \hline
\multirow{15}{*}{DC2}           & \multicolumn{1}{c|}{\multirow{5}{*}{20}} & \multicolumn{1}{c|}{SRE(dB)}  & 17.7750 (0.0822)   & 18.3910 (0.0443)   & 17.9272 (0.0792)  & 18.3065 (0.0573)  \\ \cline{3-7}
                                & \multicolumn{1}{c|}{}                    & \multicolumn{1}{c|}{$p_s$}    & 0.9995 (0.0002)  & 0.9999 (0.0001)       & 1.0000 (0) & 1.0000 (0)      \\ \cline{3-7}
                                & \multicolumn{1}{c|}{}                    & \multicolumn{1}{c|}{time(s)}  & 195.8421 (1.0192)  & 27.6035 (0.0601)    & 197.4170 (0.6048)  & 28.2088 (0.0974)   \\ \cline{3-7}
                                & \multicolumn{1}{c|}{}                    & \multicolumn{1}{c|}{$\lambda$}& 0.0005  & 0.0005  & 0.01   & 0.0005 \\ \cline{3-7}
                                & \multicolumn{1}{c|}{}               & \multicolumn{1}{c|}{$\lambda_{TV}$}& 0.0001  & 0.0001  & 0.0001 & 0.0001 \\ \cline{2-7}

                                & \multicolumn{1}{c|}{\multirow{5}{*}{30}} & \multicolumn{1}{c|}{SRE(dB)}  & 24.7523 (0.0260)   & 25.2095 (0.0273)   & 25.3758 (0.0250)  & 25.7898 (0.0419)  \\ \cline{3-7}
                                & \multicolumn{1}{c|}{}                    & \multicolumn{1}{c|}{$p_s$}    & 1.0000 (0)       & 1.0000 (0)       & 1.0000 (0)      & 1.0000 (0)  \\ \cline{3-7}
                                & \multicolumn{1}{c|}{}                    & \multicolumn{1}{c|}{time(s)}  & 194.9738 (0.7467)   & 27.6767 (0.1028)    & 197.7280 (0.9522)  & 27.5906 (0.0906)   \\ \cline{3-7}
                                & \multicolumn{1}{c|}{}                    & \multicolumn{1}{c|}{$\lambda$}& 0.0005  & 0.0005  & 0.01   & 0.005  \\ \cline{3-7}
                                & \multicolumn{1}{c|}{}               & \multicolumn{1}{c|}{$\lambda_{TV}$}& 0.0001  & 0.0001  & 0.0001 & 0.0005 \\ \cline{2-7}

                                & \multicolumn{1}{c|}{\multirow{5}{*}{40}} & \multicolumn{1}{c|}{SRE(dB)}  & 27.0058 (0.0081)   & 27.2726 (0.0091)   & 27.9773 (0.0249)  & 28.3003 (0.0189)  \\ \cline{3-7}
                                & \multicolumn{1}{c|}{}                    & \multicolumn{1}{c|}{$p_s$}    & 1.0000 (0)       & 1.0000 (0)       & 1.0000 (0) & 1.0000 (0)  \\ \cline{3-7}
                                & \multicolumn{1}{c|}{}                    & \multicolumn{1}{c|}{time(s)}  & 195.4797 (1.0635)     & 27.4218 (0.0814)    & 197.6918 (0.3900)  & 27.3723 (0.3065)   \\ \cline{3-7}
                                & \multicolumn{1}{c|}{}                    & \multicolumn{1}{c|}{$\lambda$}& 0.0005  & 0.0005  & 0.01   & 0.005  \\ \cline{3-7}
                                & \multicolumn{1}{c|}{}               & \multicolumn{1}{c|}{$\lambda_{TV}$}& 0.0001  & 0.0001  & 0.0005 & 0.0005 \\ \cline{1-7}
                                \multicolumn{7}{|l|}{}                                                                                                                                                                                                                     \\ \hline
\multirow{15}{*}{DC3}           & \multicolumn{1}{c|}{\multirow{5}{*}{20}} & \multicolumn{1}{c|}{SRE(dB)}  & 16.1060 (0.1039)   & 17.3572 (0.0635)   & 16.0405 (0.0862)  & 16.7953 (0.0507)  \\ \cline{3-7}
                                & \multicolumn{1}{c|}{}                    & \multicolumn{1}{c|}{$p_s$}    & 0.9938 (0.0010)  & 0.9982 (0.0004)  & 0.9962 (0.0007) & 0.9983 (0.0004)  \\ \cline{3-7}
                                & \multicolumn{1}{c|}{}                    & \multicolumn{1}{c|}{time(s)}  & 193.5519 (0.7429)   & 27.6214 (0.0663)    & 196.1045 (0.9418)  & 28.2954 (0.0854)   \\ \cline{3-7}
                                & \multicolumn{1}{c|}{}                    & \multicolumn{1}{c|}{$\lambda$}& 0.001   & 0.001   & 0.01   & 0.0005 \\ \cline{3-7}
                                & \multicolumn{1}{c|}{}               & \multicolumn{1}{c|}{$\lambda_{TV}$}& 0.0001  & 0.0001  & 0.0001 & 0.0001 \\ \cline{2-7}

                                & \multicolumn{1}{c|}{\multirow{5}{*}{30}} & \multicolumn{1}{c|}{SRE(dB)}  & 24.9063 (0.0700)    & 26.6341 (0.0331)   & 24.2799 (0.0487)  & 23.6197 (0.0420)  \\ \cline{3-7}
                                & \multicolumn{1}{c|}{}                    & \multicolumn{1}{c|}{$p_s$}    & 1.0000 (0)       & 1.0000 (0)       & 1.0000 (0)      & 1.0000 (0)   \\ \cline{3-7}
                                & \multicolumn{1}{c|}{}                    & \multicolumn{1}{c|}{time(s)}  & 193.9397 (0.4759)   & 27.8804 (0.2705)    & 196.7967 (0.7968)    & 28.3764 (0.0515)   \\ \cline{3-7}
                                & \multicolumn{1}{c|}{}                    & \multicolumn{1}{c|}{$\lambda$}& 0.0005  & 0.0005  & 0.01   & 0.0005  \\ \cline{3-7}
                                & \multicolumn{1}{c|}{}               & \multicolumn{1}{c|}{$\lambda_{TV}$}& 0.0001  & 0.0001  & 0.0001 & 0.0001 \\ \cline{2-7}

                                & \multicolumn{1}{c|}{\multirow{5}{*}{40}} & \multicolumn{1}{c|}{SRE(dB)}  & 32.1884 (0.0375)   & 33.1761 (0.0307)   & 29.9157 (0.0239)  & 29.7300 (0.0202)  \\ \cline{3-7}
                                & \multicolumn{1}{c|}{}                    & \multicolumn{1}{c|}{$p_s$}    & 1.0000 (0)       & 1.0000 (0)       & 1.0000 (0)      & 1.0000 (0)   \\ \cline{3-7}
                                & \multicolumn{1}{c|}{}                    & \multicolumn{1}{c|}{time(s)}  & 194.3170 (0.5253)   & 27.7275 (0.0925)    & 197.2733 (1.1988)  & 28.3184 (0.0790)   \\ \cline{3-7}
                                & \multicolumn{1}{c|}{}                    & \multicolumn{1}{c|}{$\lambda$}& 0.0005   & 0.0005  & 0.01   & 0.0005  \\ \cline{3-7}
                                & \multicolumn{1}{c|}{}               & \multicolumn{1}{c|}{$\lambda_{TV}$}& 0.0001  & 0.0001  & 0.0001 & 0.0001 \\ \cline{1-7}
\end{tabular}}
\end{table*}
\begin{figure*}[!h]
\centering
\mbox{
{\includegraphics[width=1.4in]{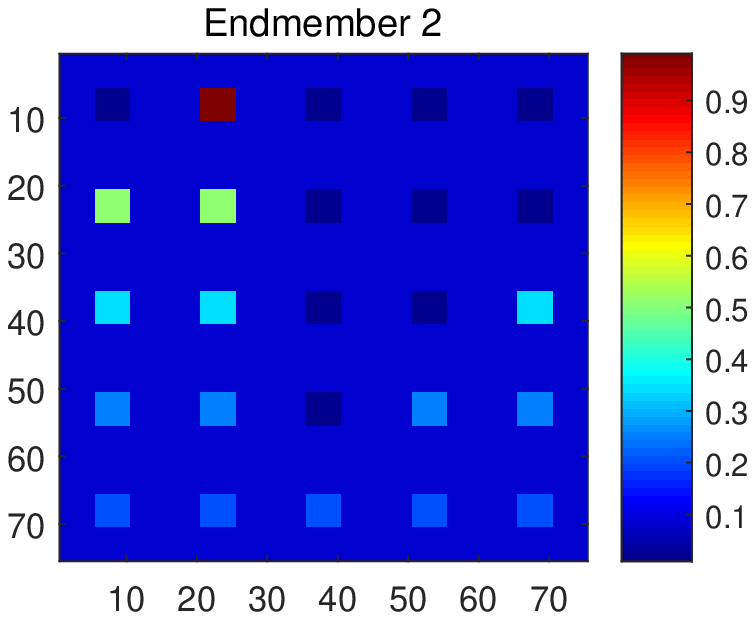}}
{\includegraphics[width=1.4in]{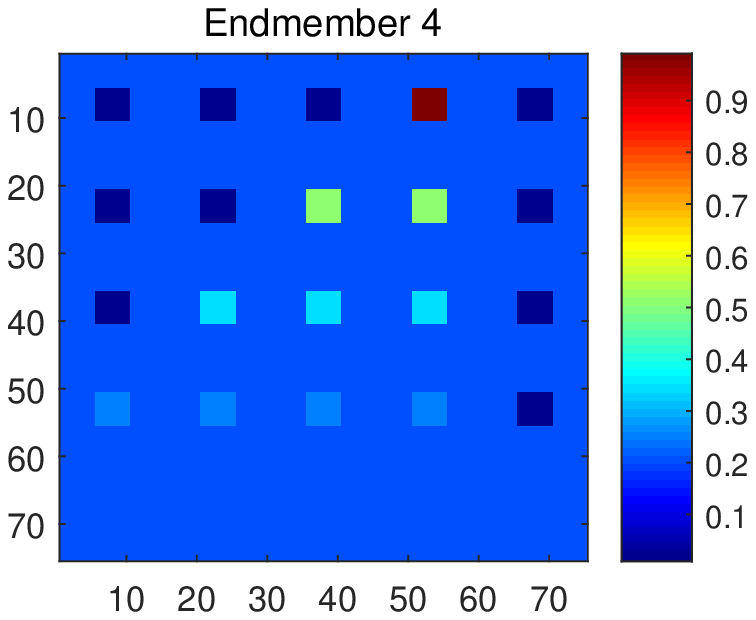}}
{\includegraphics[width=1.4in]{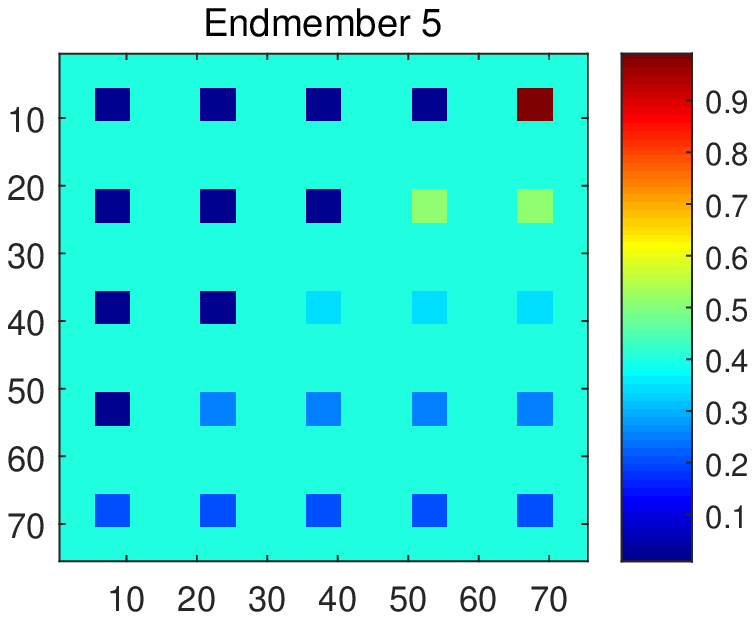}}
}\\(a)

\mbox{
{\includegraphics[width=1.4in]{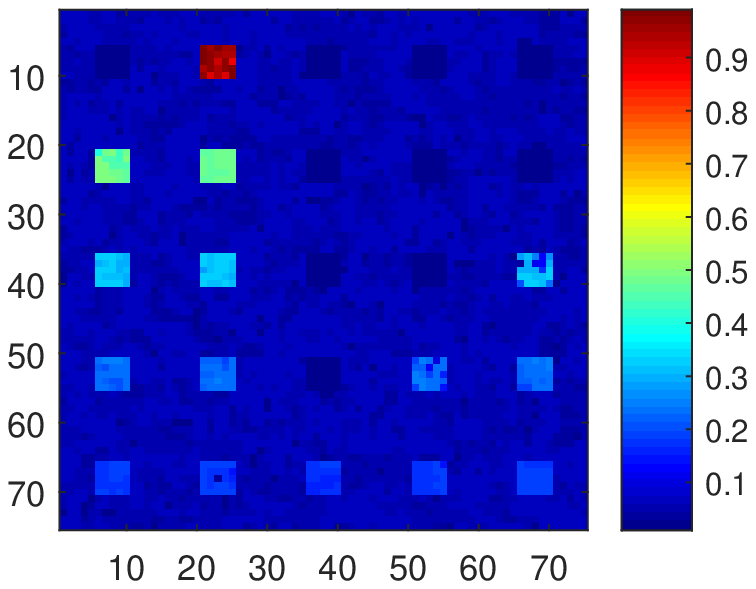}}
{\includegraphics[width=1.4in]{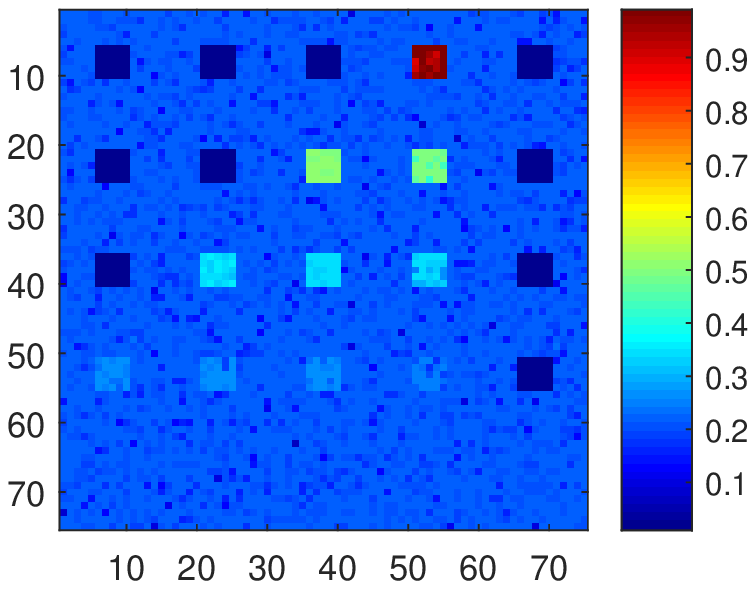}}
{\includegraphics[width=1.4in]{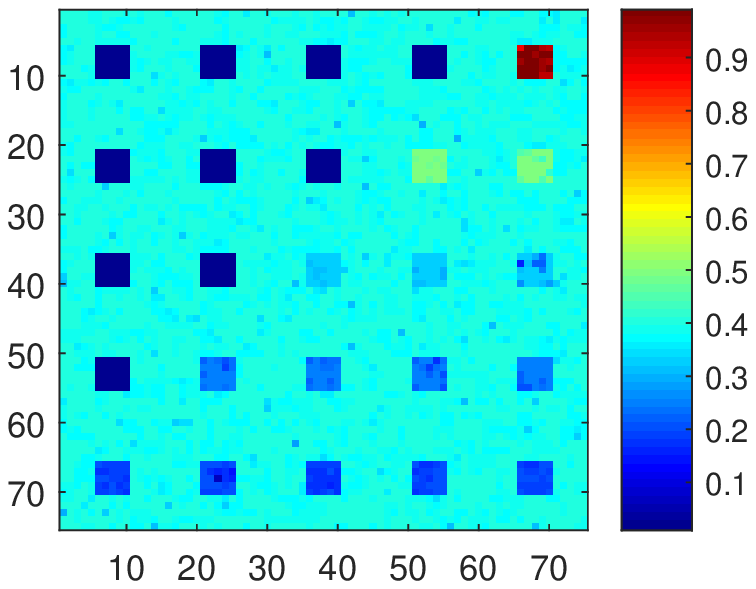}}
}\\(b)

\mbox{
{\includegraphics[width=1.4in]{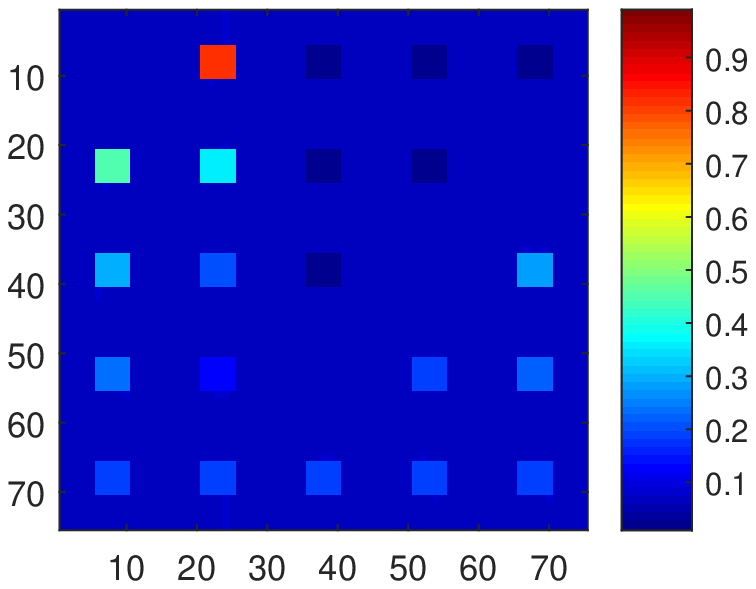}}
{\includegraphics[width=1.4in]{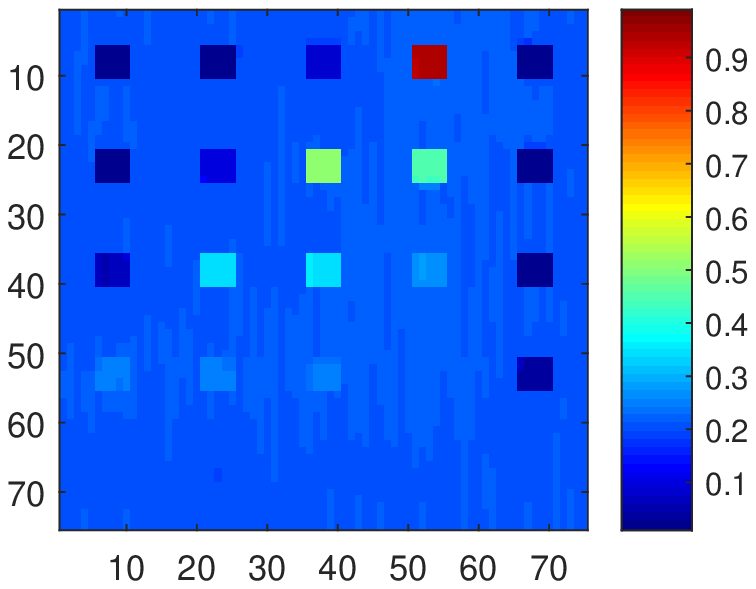}}
{\includegraphics[width=1.4in]{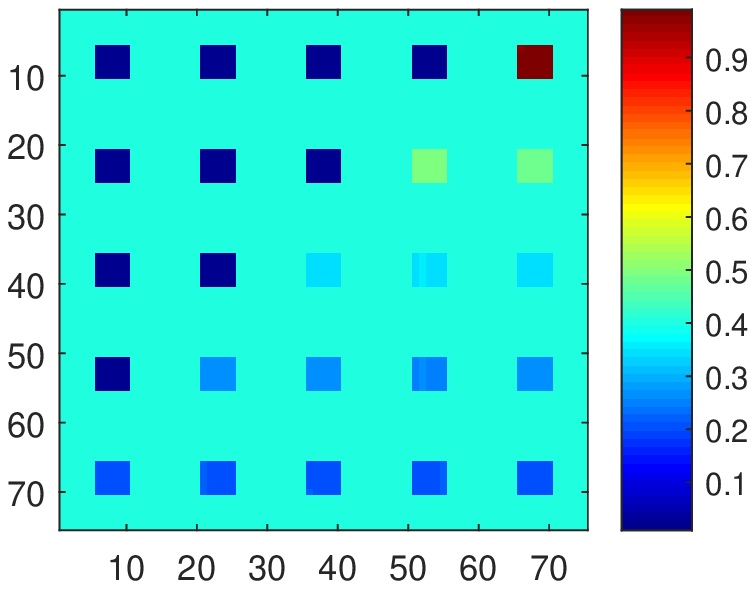}}
}\\(c)

\mbox{
{\includegraphics[width=1.4in]{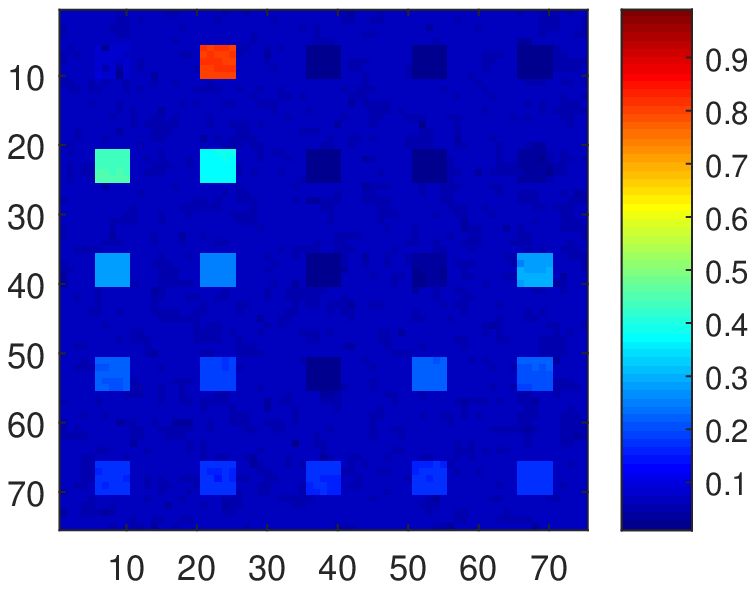}}
{\includegraphics[width=1.4in]{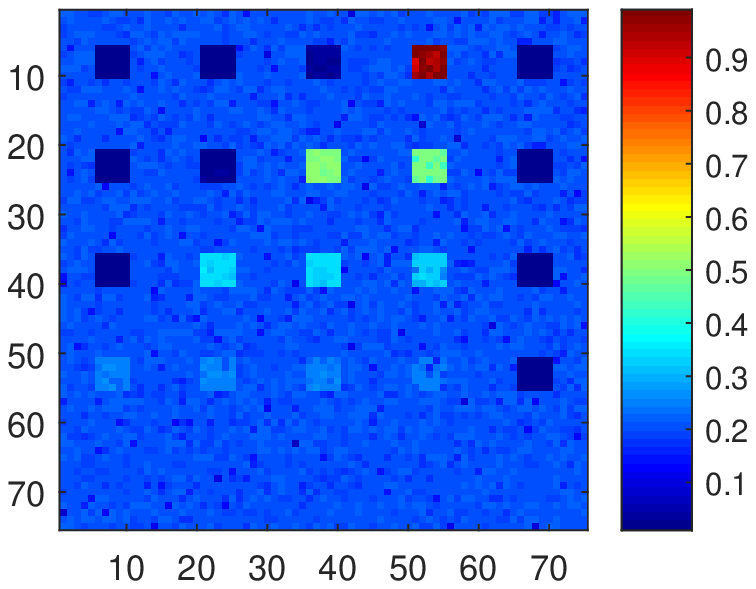}}
{\includegraphics[width=1.41in]{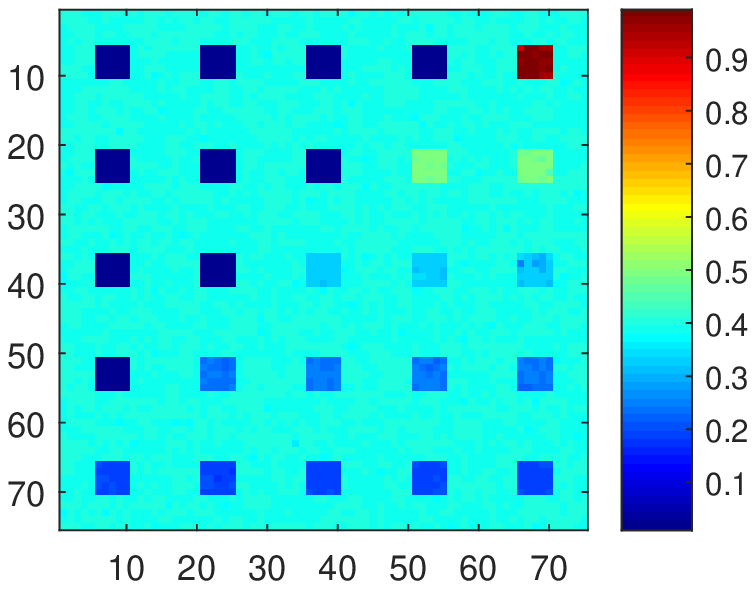}}
}\\(d)

\mbox{
{\includegraphics[width=1.4in]{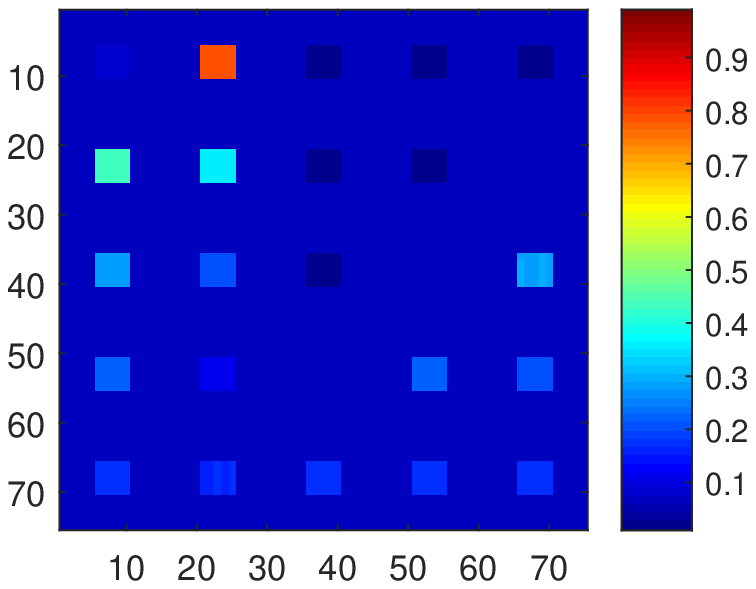}}
{\includegraphics[width=1.4in]{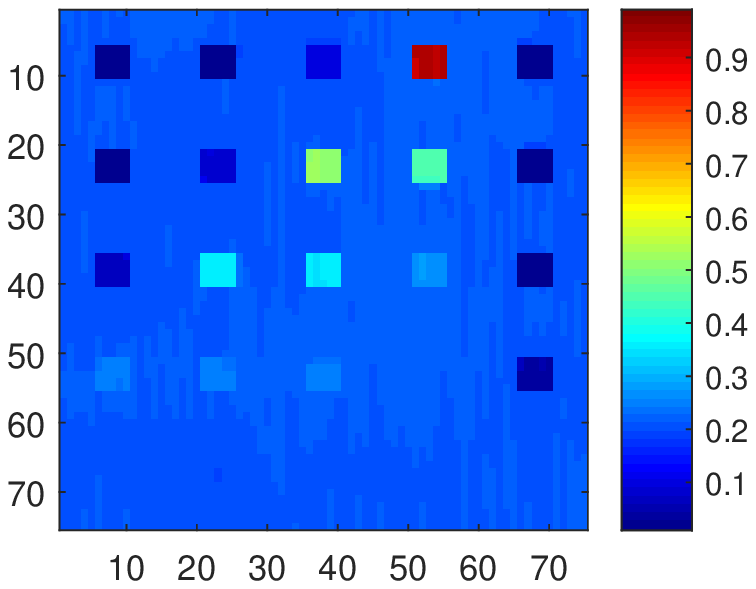}}
{\includegraphics[width=1.4in]{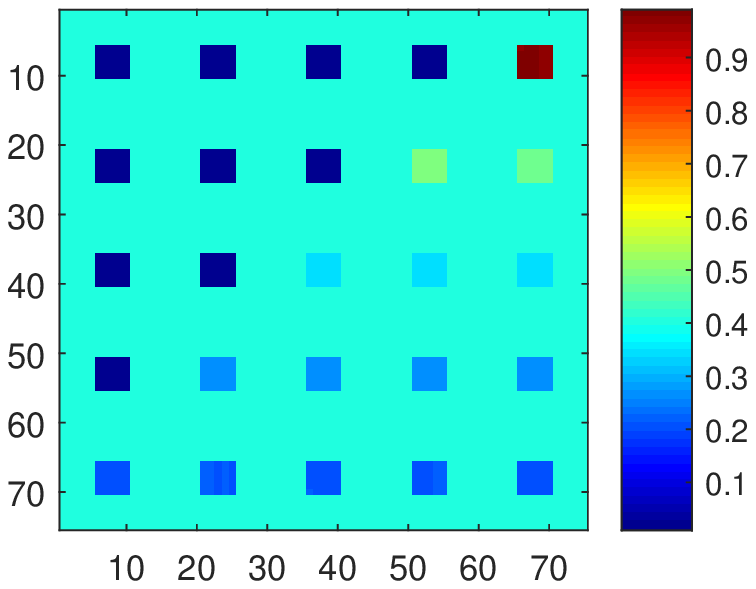}}
}\\(e)

\caption{The reference abundances and the estimated abundances obtained by the different unmixing algorithms for endmember 2, endmember 4, endmember 5 in the DC1 with SNR=20 dB (correlated noise). From top to bottom: (a) Reference abundances; (b) SUnSAL-TV ($\lambda=0.005, \lambda_{TV}=0.0001$); (c) SUnTV-sGSADMM ($\lambda=0.001, \lambda_{TV}=0.01$); (d) CLSUnSAL-TV ($\lambda=0.1, \lambda_{TV}=0.0001$); (e) CLSUnTV-sGSADMM ($\lambda=0.1, \lambda_{TV}=0.01$). }
\label{DC1R}
\end{figure*}

(1) Simulated Data Cube 1 (DC1): the size of the DC1 is $75\times 75$, and each pixel contains 224 bands.
We randomly chose five endmembers from $\mathbf{A}$ and generated the abundances of endmembers following the methodology of \cite{SUnSAL}.
Then  the white noise and the correlated noise (resulting from low-pass filtering i.i.d. Gaussian noise, using a normalized cutoff frequency of $5\pi/L$) with signal-to-noise ratio ($\textrm{SNR} = E[\|\mathbf{Ax}\|_2^2]/E[\|\mathbf{n}\|_2^2]$) of 20 dB, 30 dB and 40 dB were added to the DC1, respectively.

(2) Simulated Data Cube 2 (DC2): the size of the DC2 is $100\times100$, and each pixel contains 224 bands.
9 endmembers were randomly chosen from $\mathbf{A}$. The abundances of the DC2 satisfy the ANC and the ASC, and it was generated based on the Gaussian fields method whose type is Mattern \cite{Mattern1999}. Similar to the DC1, the white noise and the correlated noise were added to the DC2 with different SNR.

(3) Simulated Data Cube 3 (DC3): the way we generated the DC3 is like the way of generating the DC2 in \cite{SUnSAL-TV}.
Similarly, 9 endmembers were randomly chosen from $\mathbf{A}$, and the abundances of each endmembers satisfy the ANC and the ASC.
Similar to the DC1 and DC2, the white noise and the correlated noise were added to the DC3.

In Table \ref{tabel1} and Table \ref{tabel2}, we report the mean and standard deviations (in brackets) of SRE values, $p_s$ values, computing times of the different unmixing algorithms on the simulated data and the optimal regularization parameters.
As can be seen in Table \ref{tabel1} and Table \ref{tabel2}, with regard to the DC1, the unmixing based on the dual sGS-ADMM is 8 to 9 times faster than that based on the primal ADMM in the white noise case; while the unmixing based on the dual sGS-ADMM is 9 to 10 times faster than that based on the primal ADMM in the correlated noise case.
For the DC2 and the DC3, both the SUnTV-sGSADMM and the CLSUnTV-sGSADMM are about 7 to 8 times faster than the SUnSAL-TV and the CLSUnSAL-TV for both the white noise and the correlated noise under different SNR levels.

Furthermore, we can see that almost all the SRE values and the $p_s$ values based on the dual sGS-ADMM are relatively higher than those based on the primal ADMM.
Especially for the DC1, the SRE value by the SUnTV-sGSADMM increases by more than 5 dB compared with that by the SUnSAL-TV in the Gaussian noise case with SNR=20 dB. And in the correlated noise case with SNR=20 dB, the computing time of the CLSUnTV-sGSADMM is only one ninth of the CLSUnSAL-TV, meanwhile its SRE value has increased by nearly 6 dB.

\begin{figure*}[!t]
	\centering
	\mbox{
		\subfigure[]{\includegraphics[width=1.5in]{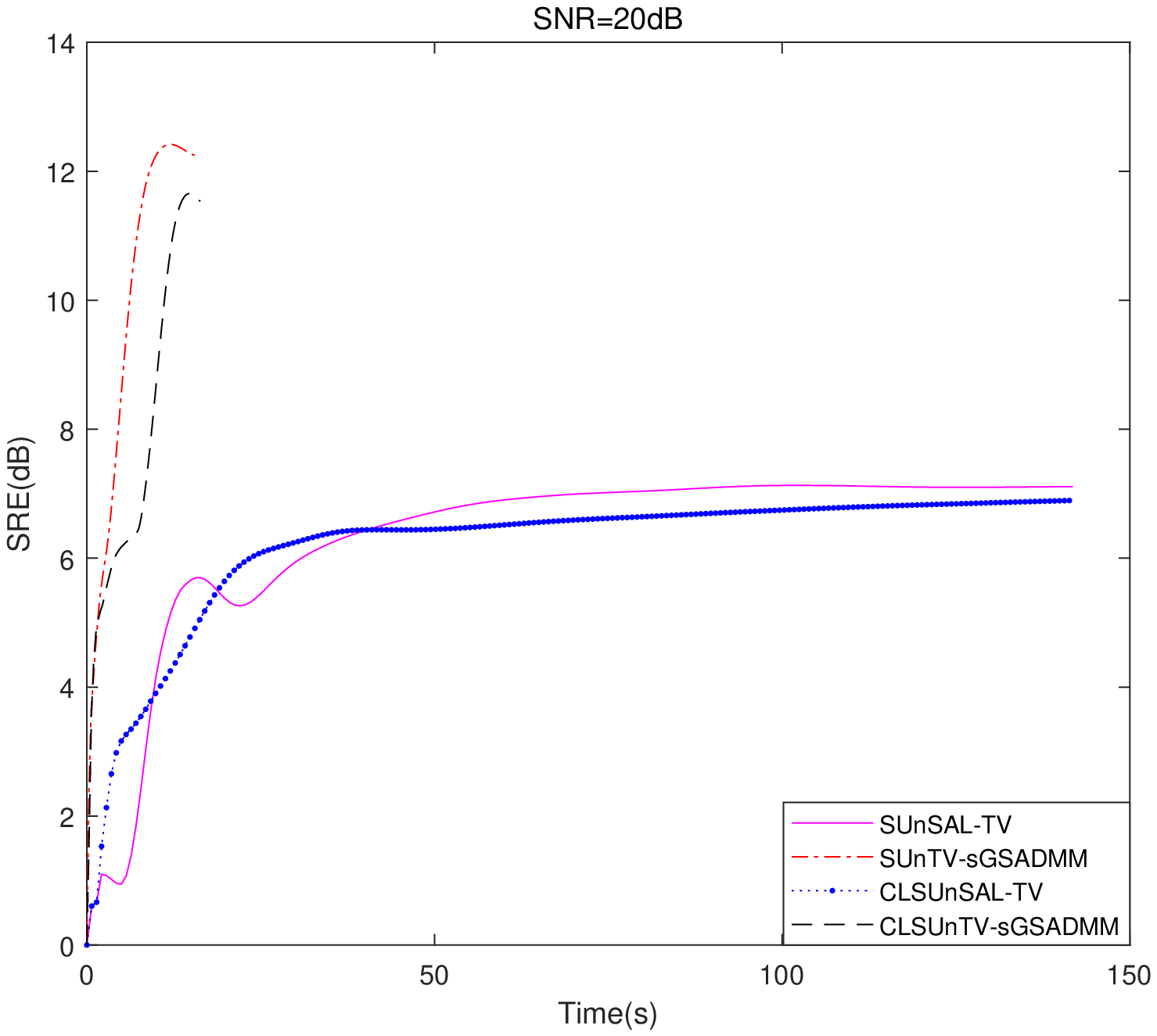}}
        \subfigure[]{\includegraphics[width=1.5in]{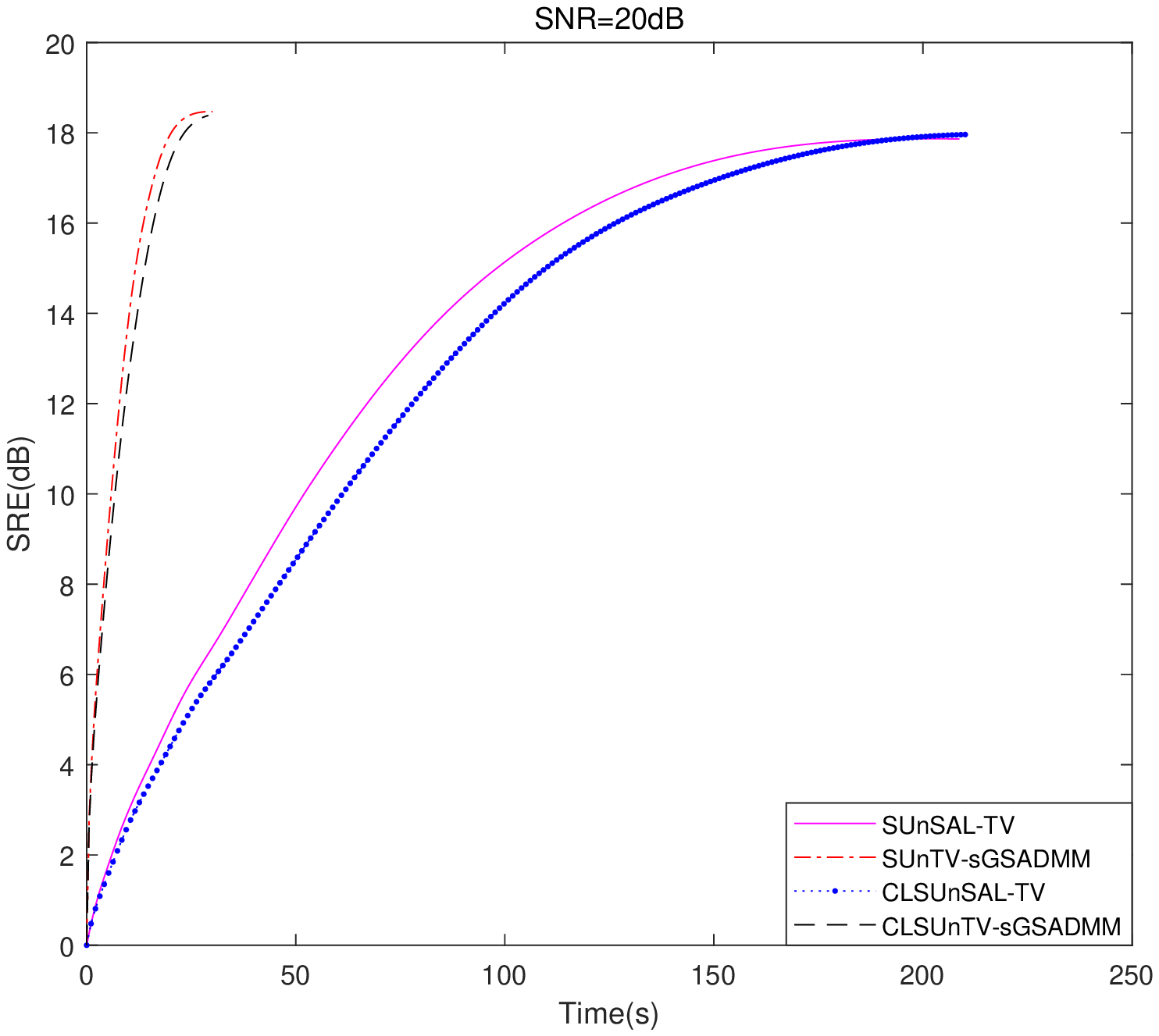}}
		\subfigure[]{\includegraphics[width=1.5in]{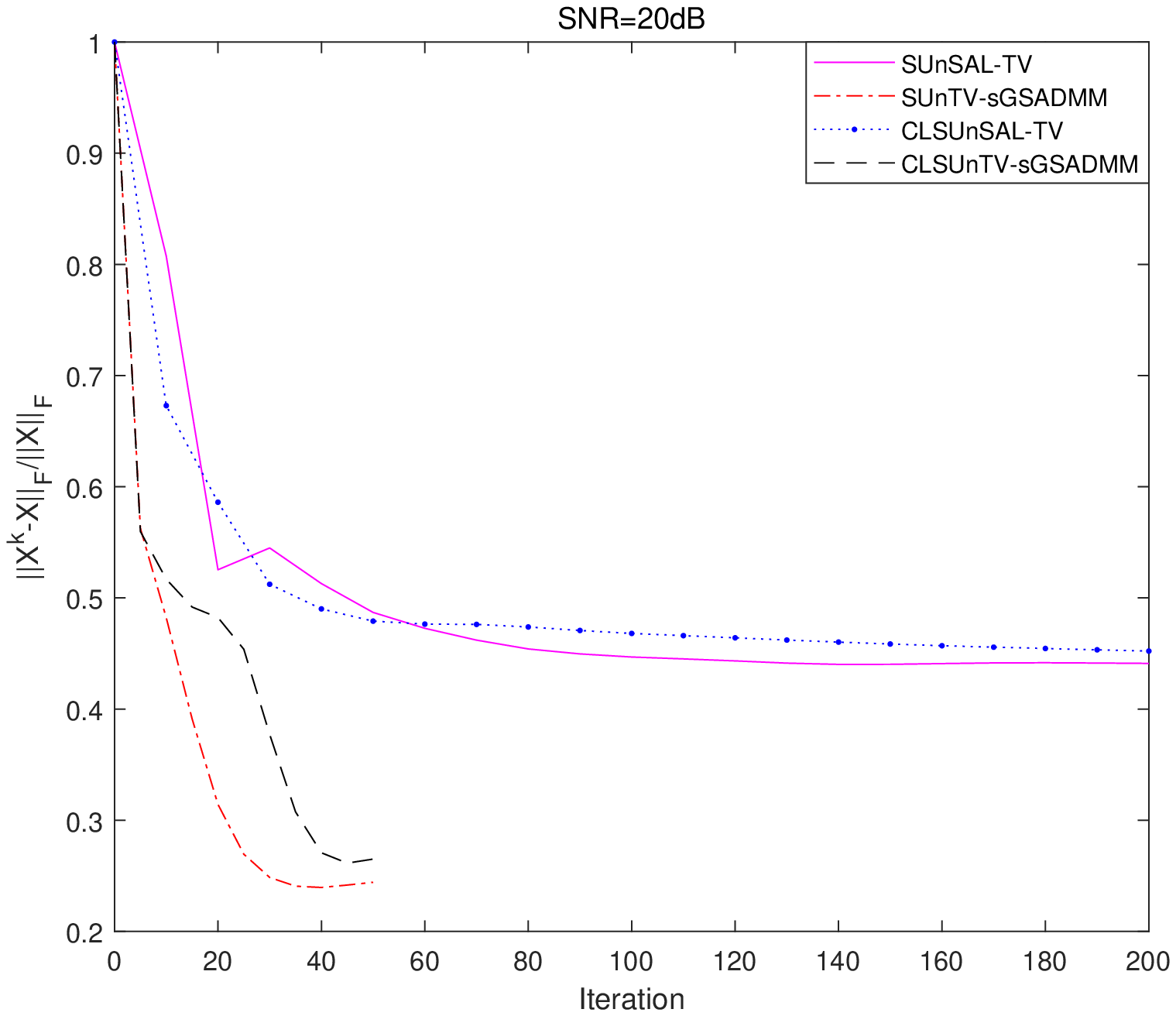}}
		
	}\\
\mbox{
		\subfigure[]{\includegraphics[width=1.5in]{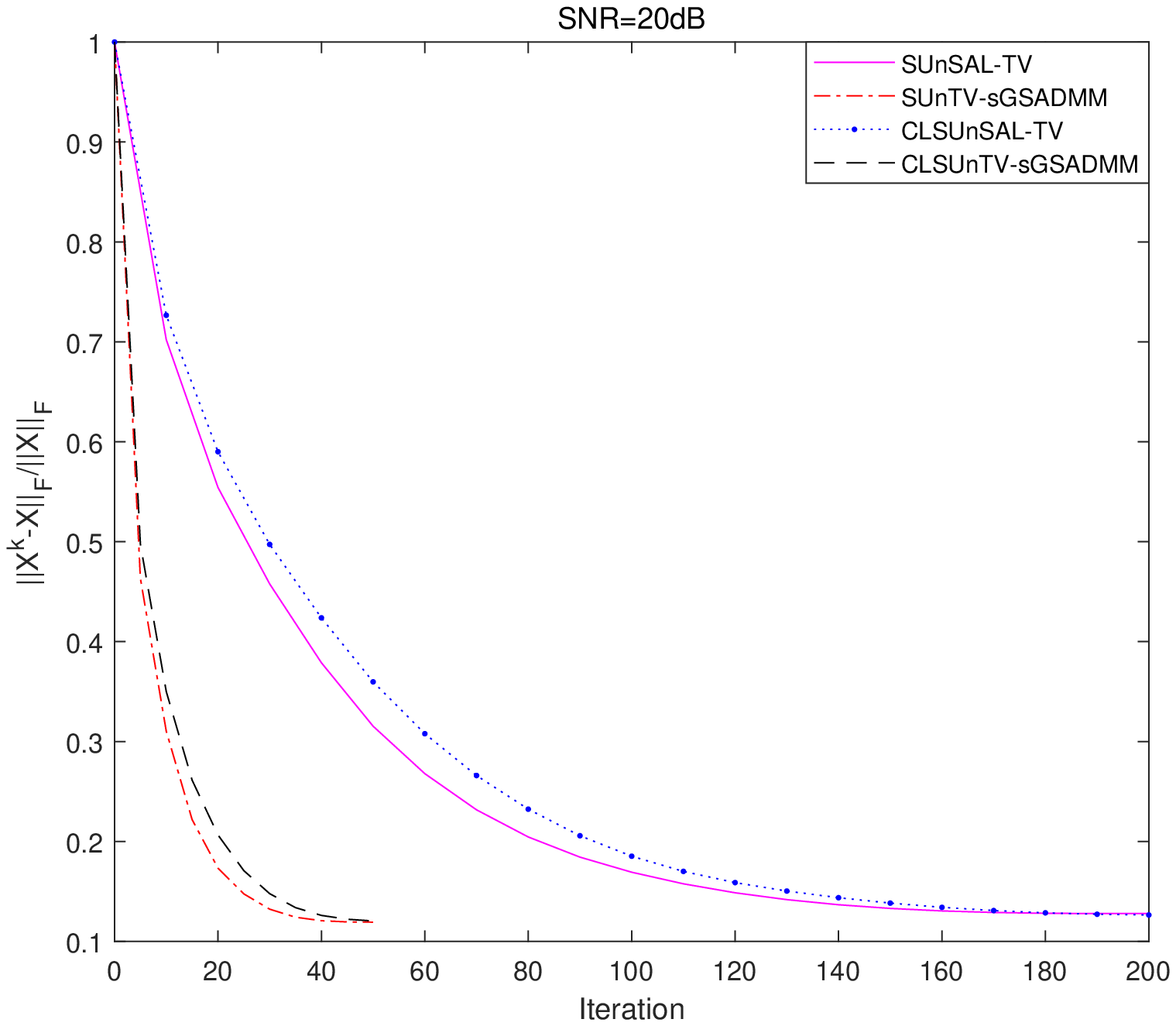}}
        \subfigure[]{\includegraphics[width=1.5in]{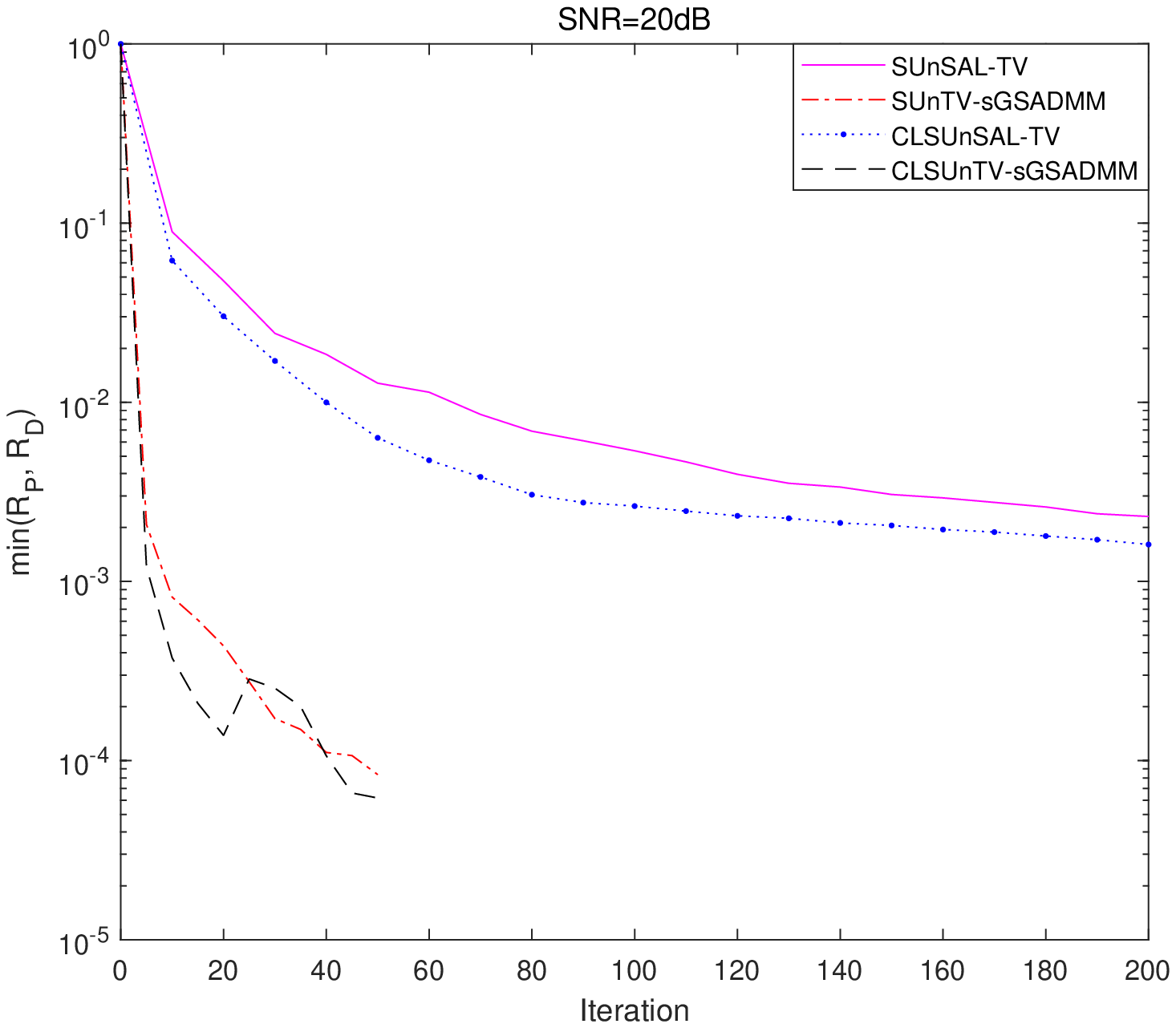}}
		\subfigure[]{\includegraphics[width=1.5in]{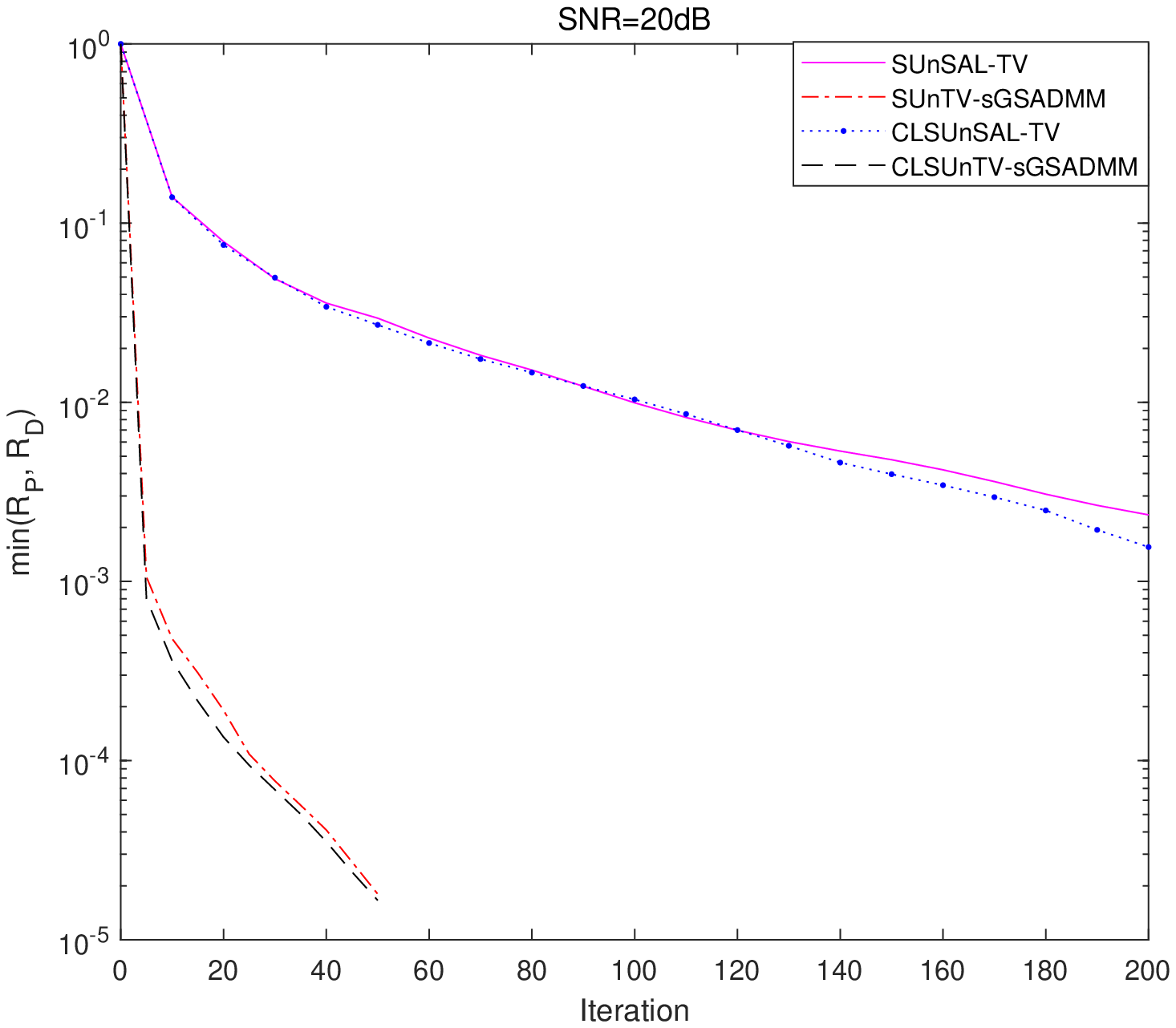}}
	}\\
	\caption{(a) (b): Variations of the SRE (dB) values with respect to the time for different algorithms when dealing with the DC1 with white noise (SNR=20 dB) and the DC2 with correlated noise (SNR=20 dB);
(c) (d): Variations of the $\|\mathbf{X}^\text{k}-\mathbf{X}\|_F/\|\mathbf{X}\|_F$ with respect to the iterations for different algorithms when dealing with the DC1 with white noise (SNR=20 dB) and the DC2 with correlated noise (SNR=20 dB);
(e) (f): Variations of $\min(\mathbf{R}_P,\mathbf{R}_D)$ with respect to the iterations  for different algorithms when dealing with the DC1 with white noise (SNR=20 dB) and the DC2 with correlated noise (SNR=20 dB).}
	\label{ConvergenceAnalysis}
\end{figure*}

For visual comparison, Fig. \ref{DC1R} shows the reference abundances and the estimated abundances obtained by the different unmixing algorithms for endmember 2, endmember 4, and endmember 5 in the DC1 with SNR = 20 dB (correlated noise).
We can see from Fig. \ref{DC1R} that the estimated abundances obtained by the primal ADMM seem to be noisy. In contrast, the estimated abundances obtained by the sGS-ADMM are more accurate and have a better visual effect.
\begin{figure*}[!t]
	\centering
	\mbox{
		\subfigure[]{\includegraphics[width=1.15in]{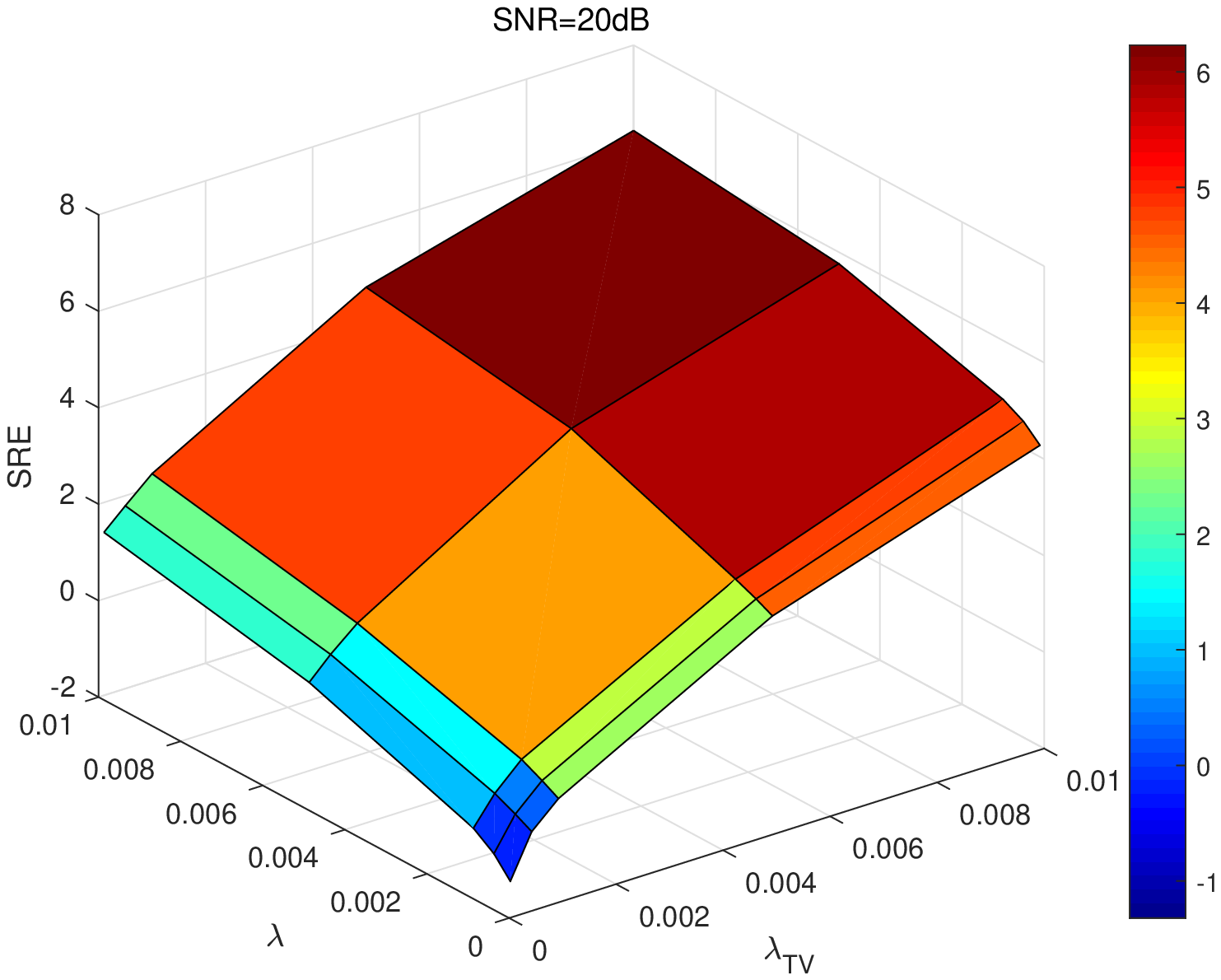}}
		\subfigure[]{\includegraphics[width=1.15in]{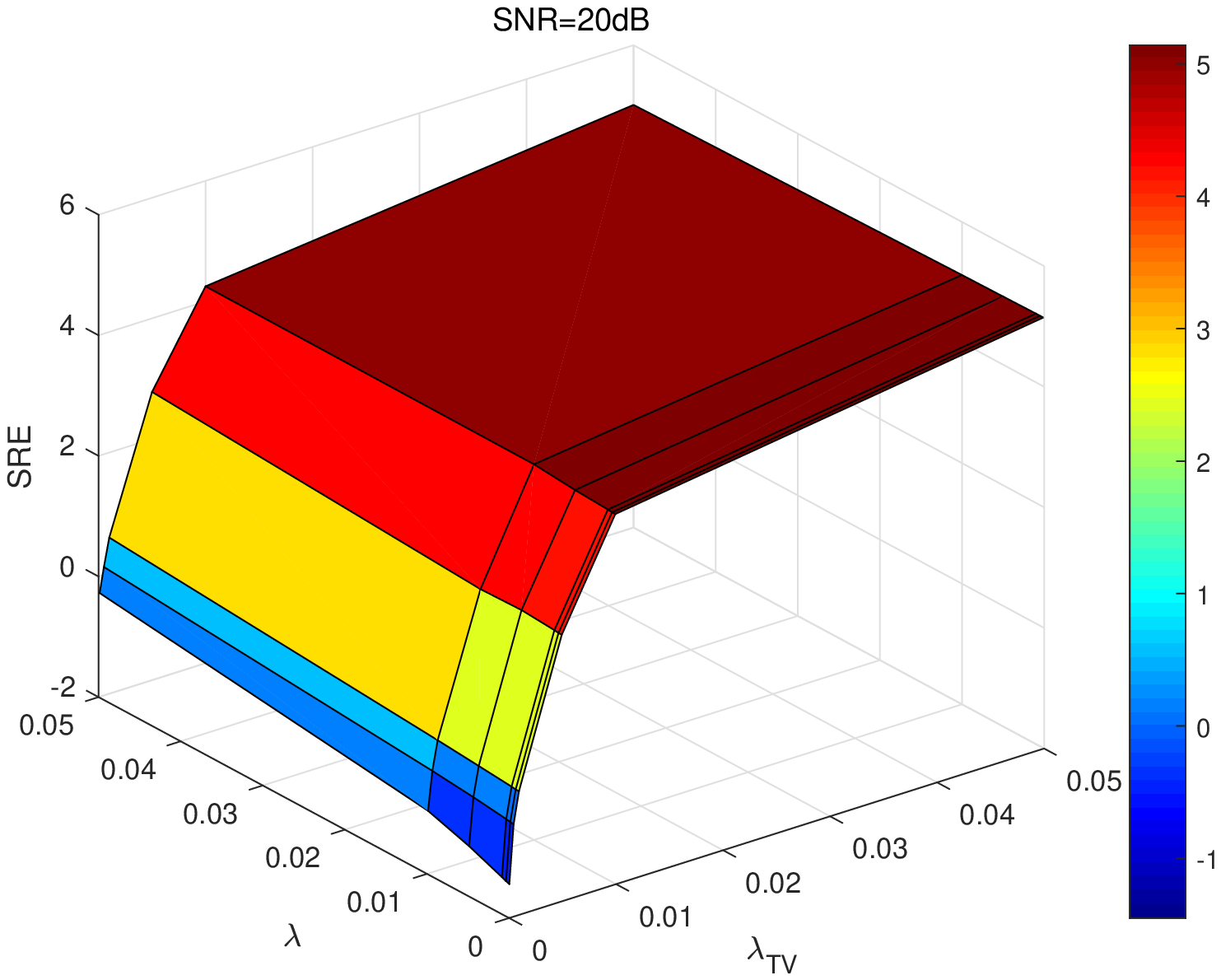}}
		\subfigure[]{\includegraphics[width=1.15in]{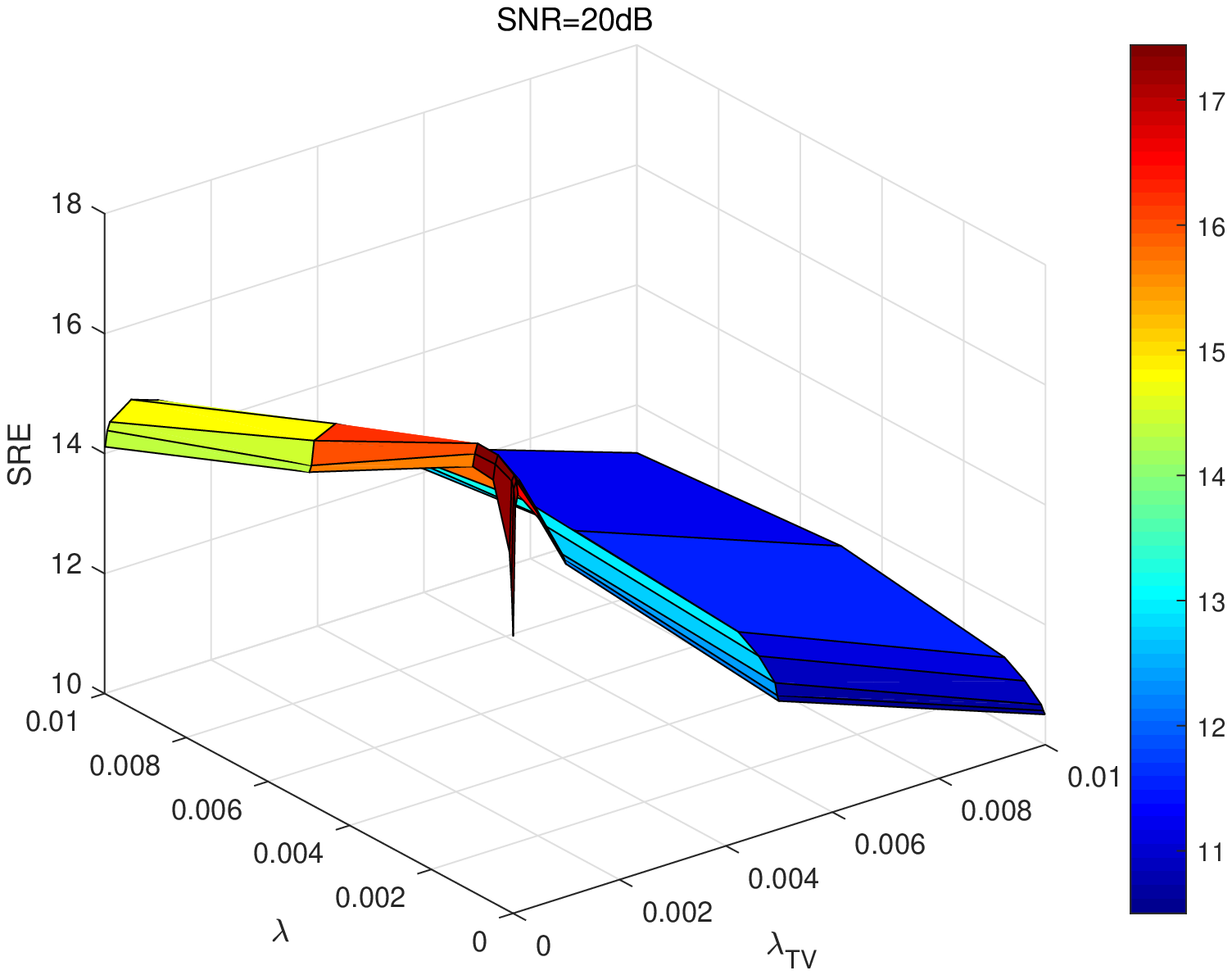}}
        \subfigure[]{\includegraphics[width=1.15in]{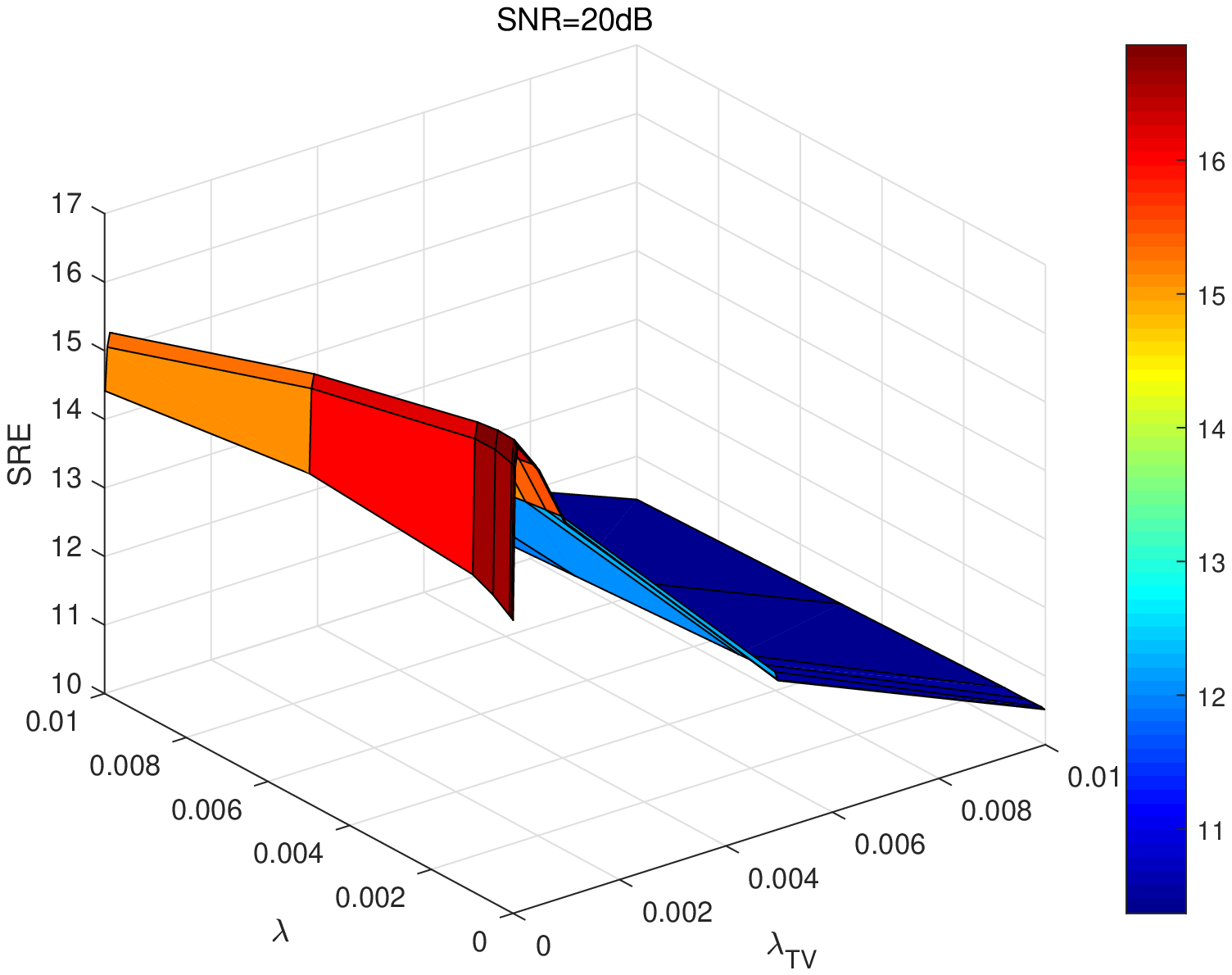}}
	}\\

	\mbox{
		\subfigure[]{\includegraphics[width=1.15in]{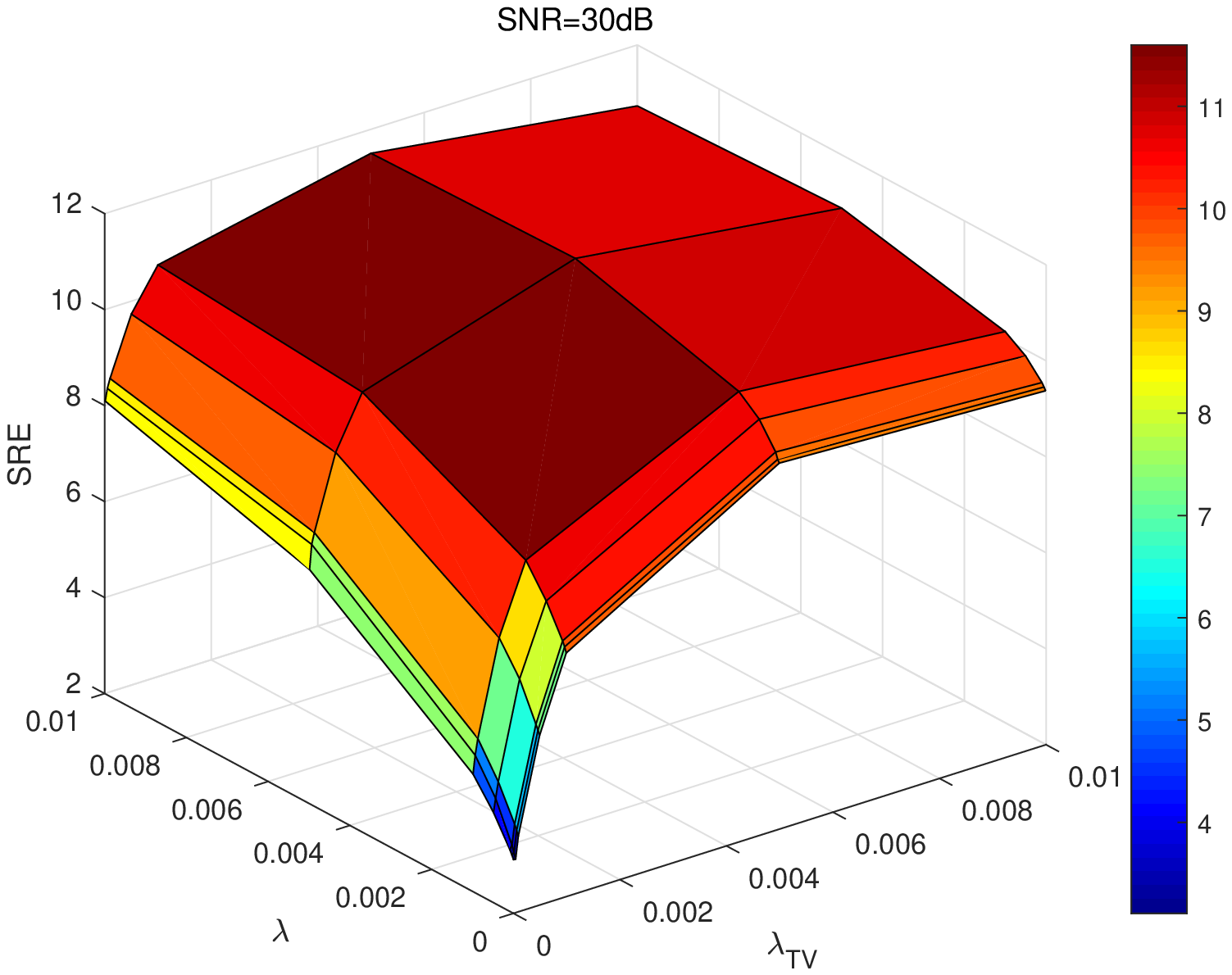}}
		\subfigure[]{\includegraphics[width=1.15in]{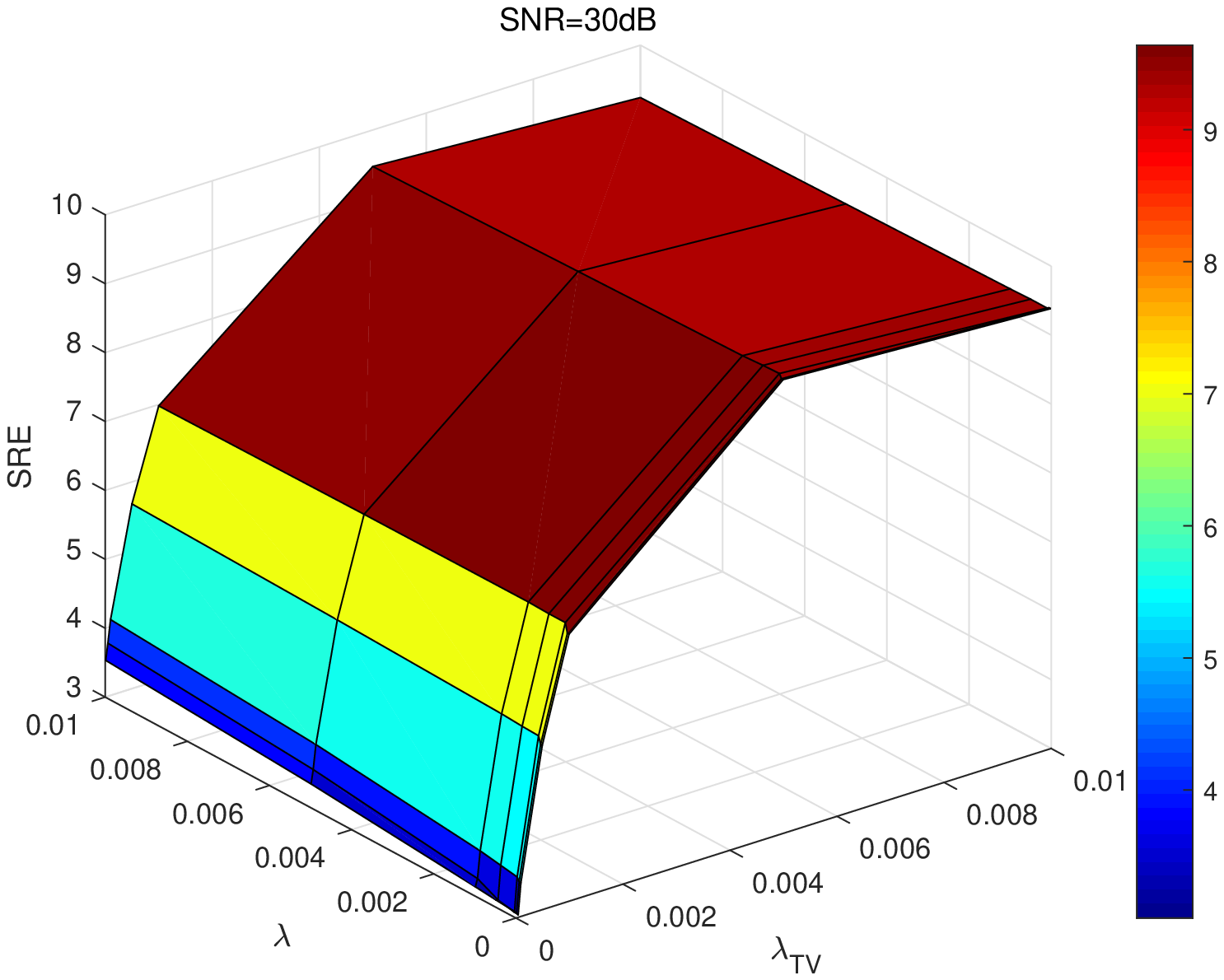}}
		\subfigure[]{\includegraphics[width=1.15in]{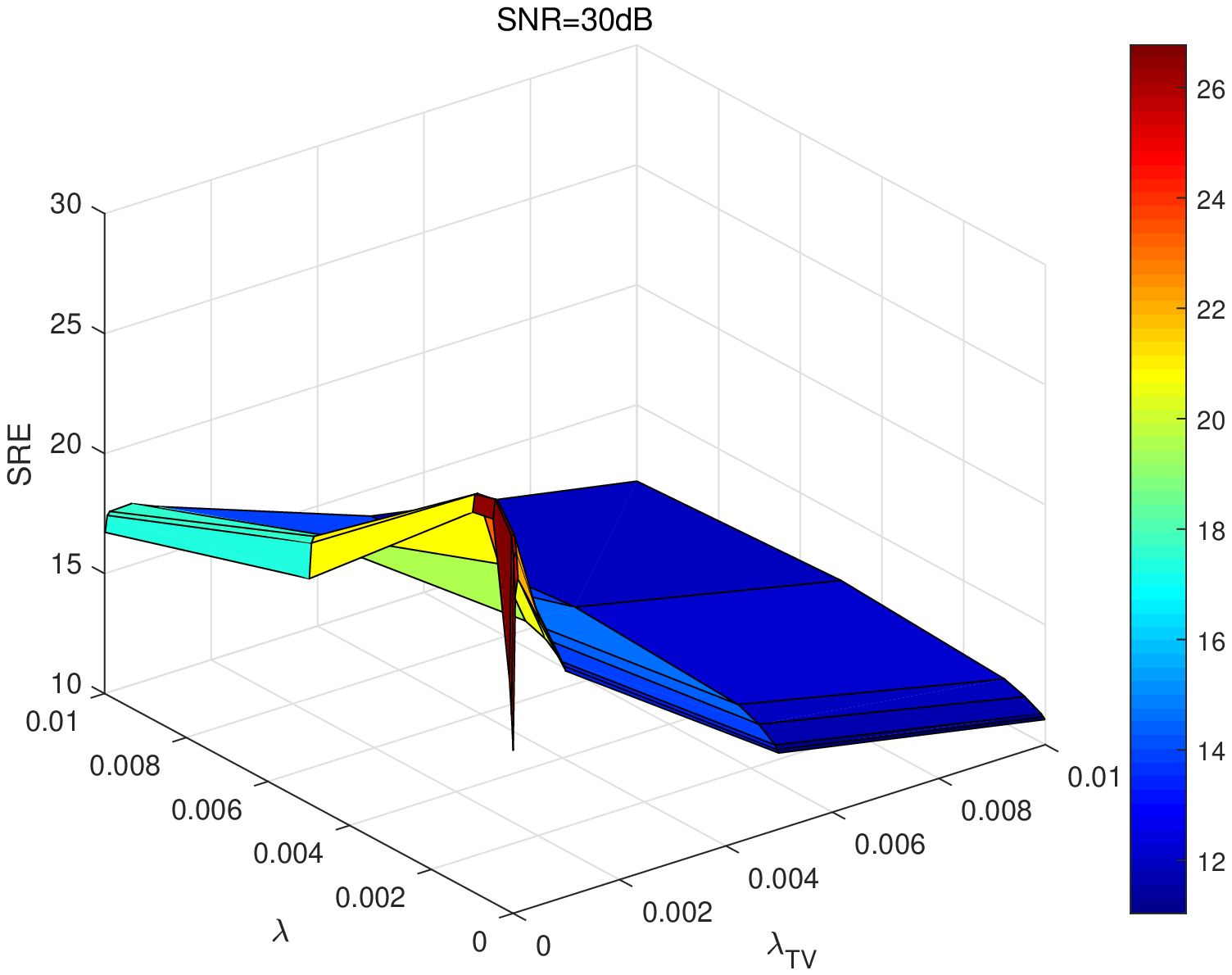}}
        \subfigure[]{\includegraphics[width=1.15in]{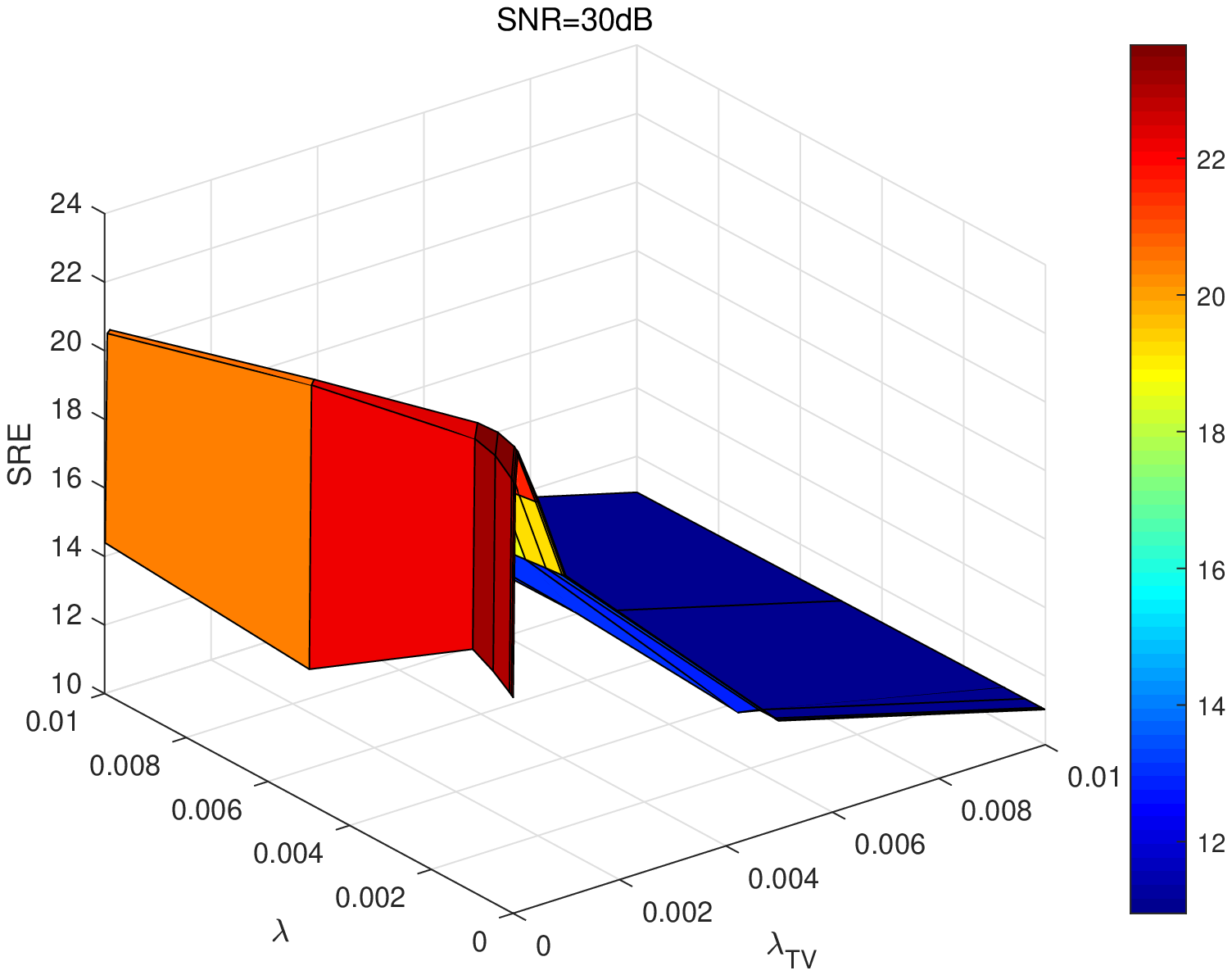}}
	}\\

\mbox{
		\subfigure[]{\includegraphics[width=1.15in]{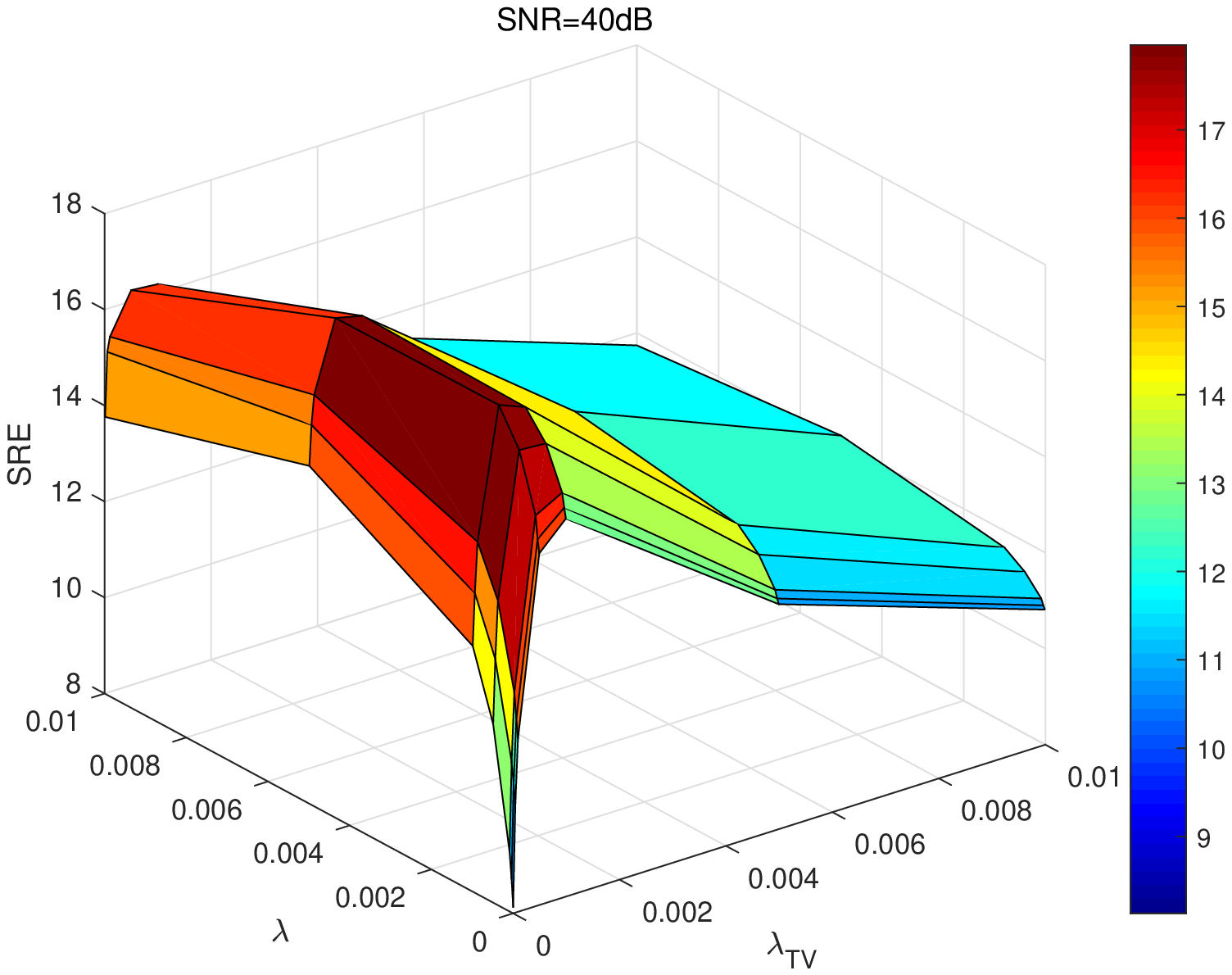}}
		\subfigure[]{\includegraphics[width=1.15in]{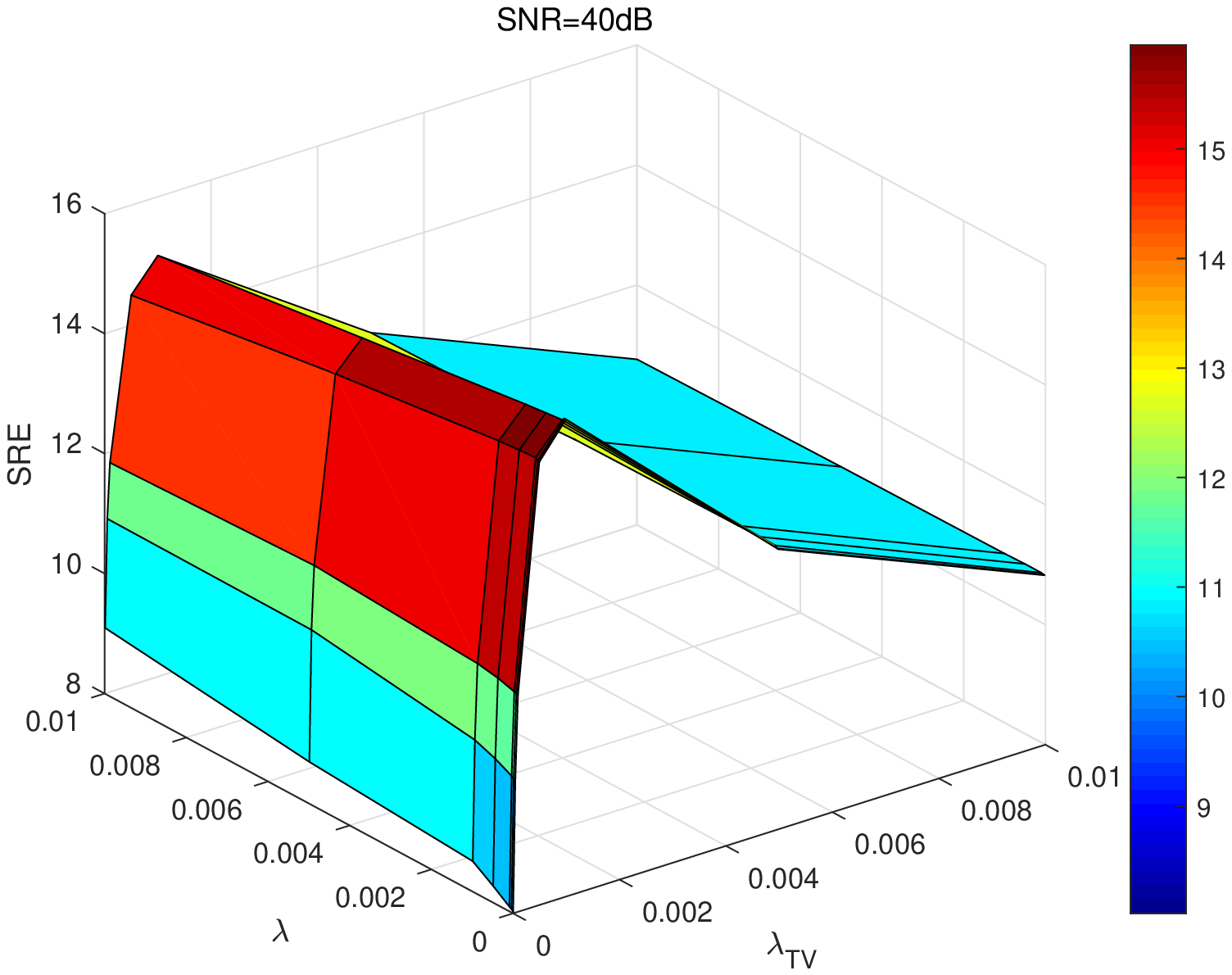}}
		\subfigure[]{\includegraphics[width=1.15in]{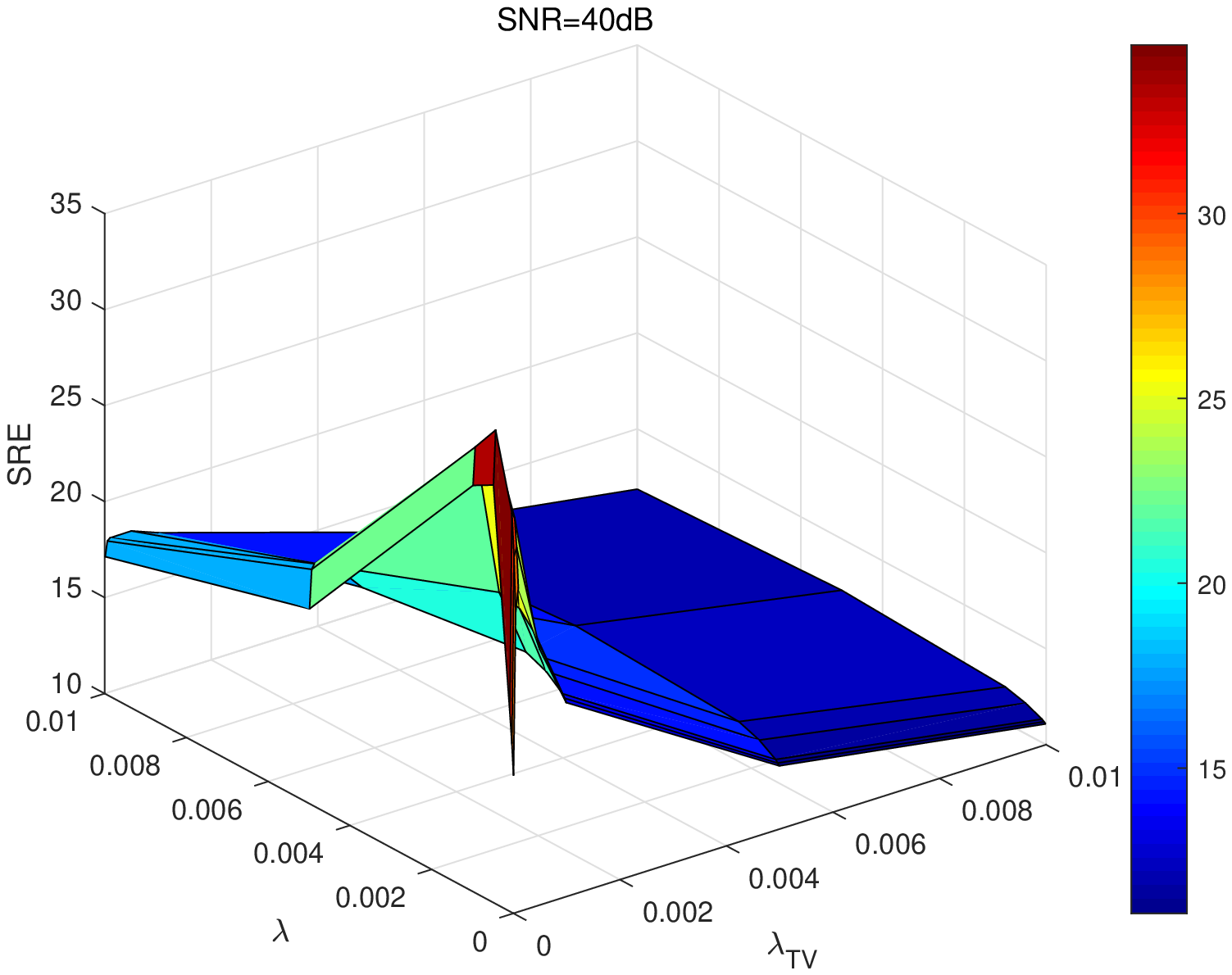}}
        \subfigure[]{\includegraphics[width=1.15in]{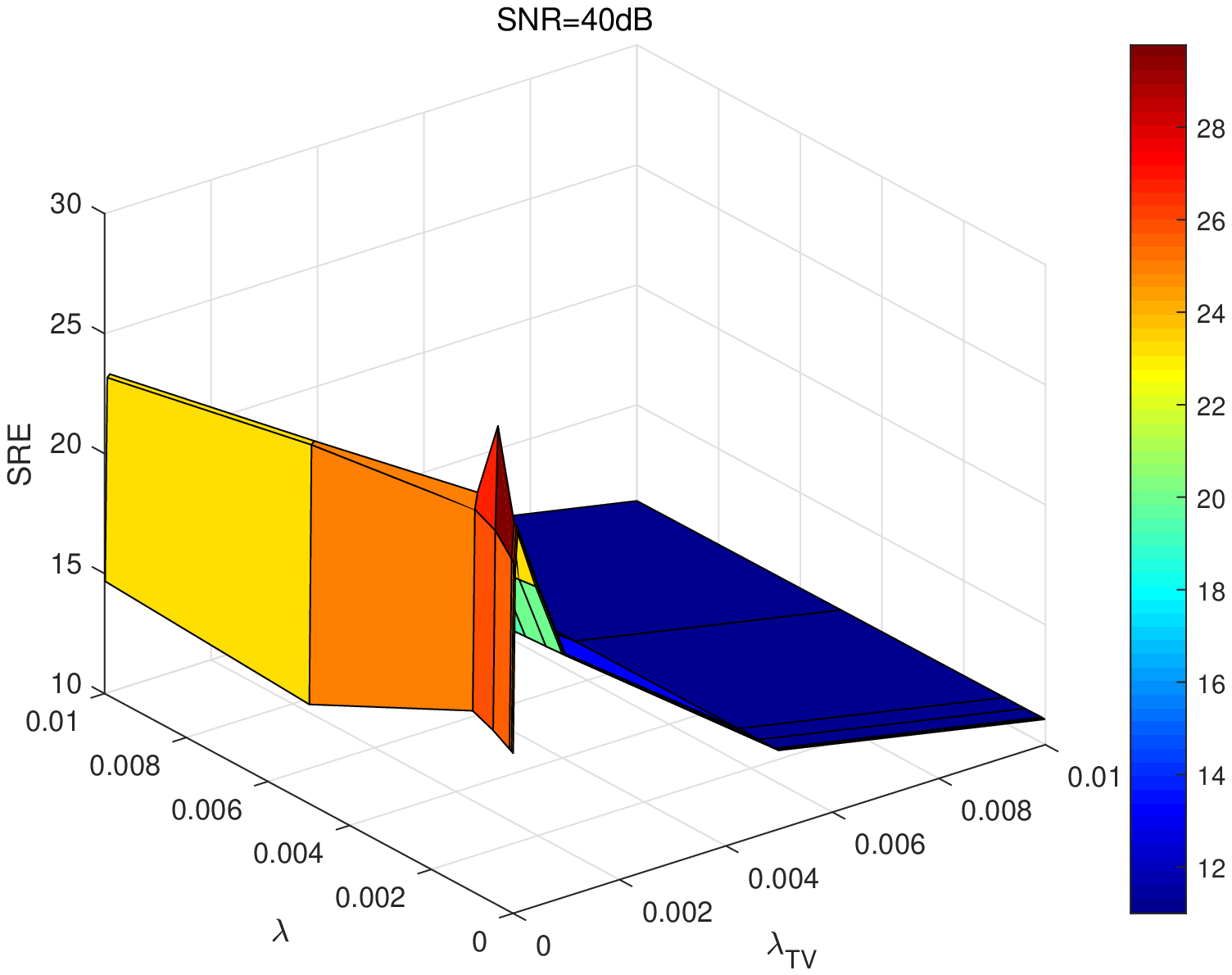}}
	}\\
	\caption{The SREs (dB) obtained by the SUnTV-sGSADMM and CLSUnTV-sGSADMM with respect to the parameters $\lambda$ and $\lambda_{TV}$ for the DC3 with different SNR levels. (a) SUnTV-sGSADMM (white noise, SNR=20 dB); (b) CLSUnTV-sGSADMM (white noise, SNR=20 dB); (c) SUnTV-sGSADMM (correlated noise, SNR=20 dB); (d) CLSUnTV-sGSADMM (correlated noise, SNR=20 dB); (e) SUnTV-sGSADMM for (white noise, SNR=30 dB); (f) CLSUnTV-sGSADMM (white noise, SNR=30 dB); (g) SUnTV-sGSADMM (correlated noise, SNR=30 dB); (h) CLSUnTV-sGSADMM (correlated noise, SNR=30 dB); (i) SUnTV-sGSADMM (white noise, SNR=40 dB); (j) CLSUnTV-sGSADMM for (white noise, SNR=40 dB); (k) SUnTV-sGSADMM (correlated noise, SNR=40 dB); (l) CLSUnTV-sGSADMM (correlated noise, SNR=40 dB).}
	\label{ParameterAanlysis}
\end{figure*}

For the primal ADMM, it introduces several variables. Once a slack variable is introduced, one more equality constraint is added. Then an extra penalty term is added in the augmented Lagrangian function, which nearly halves the step length of the corresponding variable. And the small step length may lead to the fact that we can not obtain an ideal SRE value even if the number of the iterations achieves the maximum. We also provide Fig. 2 to illustrate our claim.
In Fig. \ref{ConvergenceAnalysis}, (a) and (b) show the variations of the SRE (dB) values with respect to the time obtained by the different unmixing algorithms when dealing with the DC1 and DC2 under different SNR types.
We can see that the SRE values of the unmixing based on the dual sGS-ADMM can reach their peak values in relatively less time compared to the primal ADMM.
In Fig. \ref{ConvergenceAnalysis} (a), the peak SRE values based on the dual sGS-ADMM is obviously higher than that based on the primal ADMM,
while in Fig. \ref{ConvergenceAnalysis} (b), the peak SRE values of the unmixing based on the dual sGS-ADMM is close to that based on the primal ADMM.
To visually illustrate the convergence of different algorithms, in Fig. \ref{ConvergenceAnalysis} (c)(d) and Fig. \ref{ConvergenceAnalysis} (e)(f), we plot the variations of the relative errors between the restored abundances and the truth abundances ($\|\mathbf{X}^\text{k}-\mathbf{X}\|_F/\|\mathbf{X}\|_F$) and $\min(\mathbf{R}_P,\mathbf{R}_D)$ with respect to the iterations for different algorithms when dealing with the DC1 and DC2 under different SNR types.
By comparing these convergence curves, we can see that the dual sGS-ADMM is obviously faster than the primal ADMM.

In order to analyze the influences of parameters $\lambda$ and $\lambda_{TV}$, Fig. \ref{ParameterAanlysis} illustrates the SREs (dB) obtained by the SUnTV-sGSADMM and CLSUnTV-sGSADMM with respect to the parameters $\lambda$ and $\lambda_{TV}$ for the DC3 with different SNR levels.
From Fig. \ref{ParameterAanlysis} (a) (b) (e) (f) (i) (j), we can see that a general trend in the white noise cases is that the parameters $\lambda$ and $\lambda_{TV}$ become smaller with the SNR increasing to obtain the highest SRE values. This observation is consistent with the fact that the quadratic term $\frac{1}{2}\|\mathbf{AX}-\mathbf{Y}\|^2$ plays a dominant role in the unmixing model in the high SNR cases.
From Fig. \ref{ParameterAanlysis} (c) (d) (g) (h) (k) (l), we note that in the correlated noise cases, the parameters $\lambda$ and $\lambda_{TV}$ are obviously smaller than those in the white noise cases to obtain the highest SRE values.
This is because the variance of the correlation noise is smaller than that of the white noise in our experiments.
The performances of all algorithms  in our experiments tend to degrade when the parameters $\lambda$ and $\lambda_{TV}$ are very small (less than 0.0001).

By and large, the priority of the dual sGS-ADMM lies in the less computing time and the relatively higher SRE value. The reason is that for the primal ADMM, it introduces several variables which directly lead to smaller iteration steps, while the dual sGS-ADMM takes relatively larger steps.
And very small steps may not achieve an ideal SRE value which the model problem can supply in a reasonable amount of time. In other words, even for the same problem, the algorithm we choose is also important.

\subsection{Numerical results for the real data}

\begin{table*}[!tbp]
\scriptsize
\caption{The numerical experiments on the AVIRIS Cuprite data}
\setlength{\tabcolsep}{1.8mm}{
\centering
\label{table3}
\begin{tabular}{|c|c|c|c|c|c|}
\hline
\multirow{4}{*}{\begin{tabular}[c]{@{}c@{}}AVIRIS\\Cuprite\\Data\end{tabular}} &Parameters  & SUnSAL-TV  & SUnTV-sGSADMM  & CLSUnSAL-TV  & CLSUnTV-sGSADMM  \\ \cline{2-6}
&$\lambda$  &0.001  &0.001  &0.001  &0.001  \\ \cline{2-6}
                                                                     &$\lambda_{TV}$  &0.0001  &0.0001  &0.0001  &0.0001  \\ \cline{2-6}
                                                                     & time(s) &2794.2735 (9.7894)  &239.1765 (2.4095)  &2865.4364 (7.0131)  &243.1312 (0.7710)  \\ \hline
\end{tabular}}
\end{table*}

  \begin{table*}[!tbp]
\scriptsize
\caption{The numerical experiments on the Urban data}
\setlength{\tabcolsep}{2.2mm}{
\centering
\label{table4}
\begin{tabular}{|c|c|c|c|c|c|}
\hline
\multirow{4}{*}{\begin{tabular}[c]{@{}c@{}}Urban\\Data\end{tabular}} &Parameters  & SUnSAL-TV  & SUnTV-sGSADMM  & CLSUnSAL-TV  & CLSUnTV-sGSADMM  \\ \cline{2-6}
&$\lambda$  &0.0001  &0.0001  &0.0001  &0.0001  \\ \cline{2-6}
                                                                     &$\lambda_{TV}$  &0.0001  &0.0001  &0.0001  &0.0001  \\ \cline{2-6}
                                                                     & time(s) &153.6254 (1.2964)  &38.8155 (0.7873)  &155.1793 (1.8359)  &36.4258 (0.5744)  \\ \hline
\end{tabular}}
\end{table*}

 \begin{figure}[!t]
	\centering
     {\includegraphics[width=3.1in]{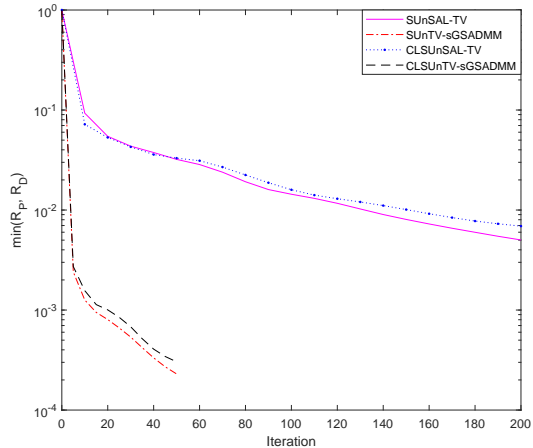}}
	\caption{Variations of $\min(\mathbf{R}_P,\mathbf{R}_D)$ with respect to the iterations  for different algorithms when dealing with the AVIRIS Cuprite data.}
	\label{ConvergenceAnalysisRealData}
\end{figure}

\begin{figure*}[!t]
\centering
\mbox{
{\includegraphics[width=1.5in]{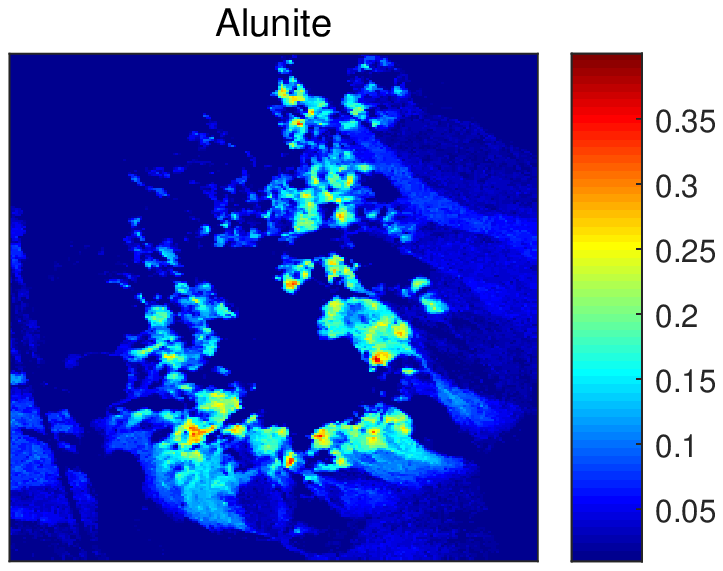}}
{\includegraphics[width=1.5in]{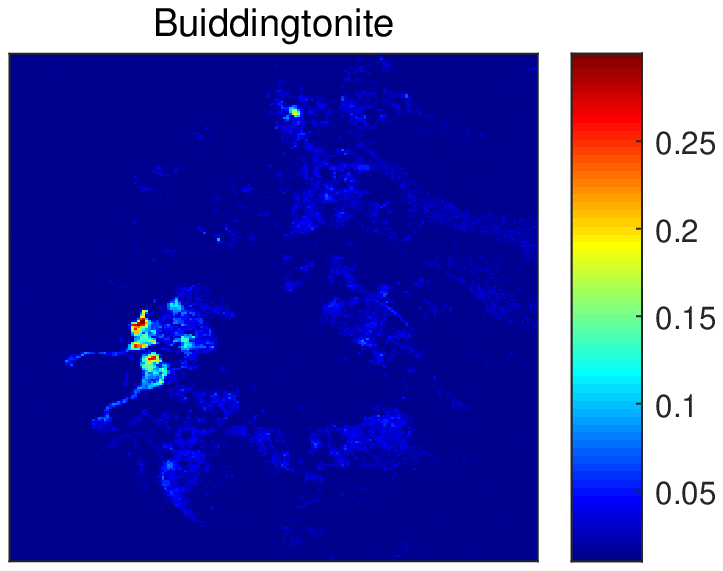}}
{\includegraphics[width=1.5in]{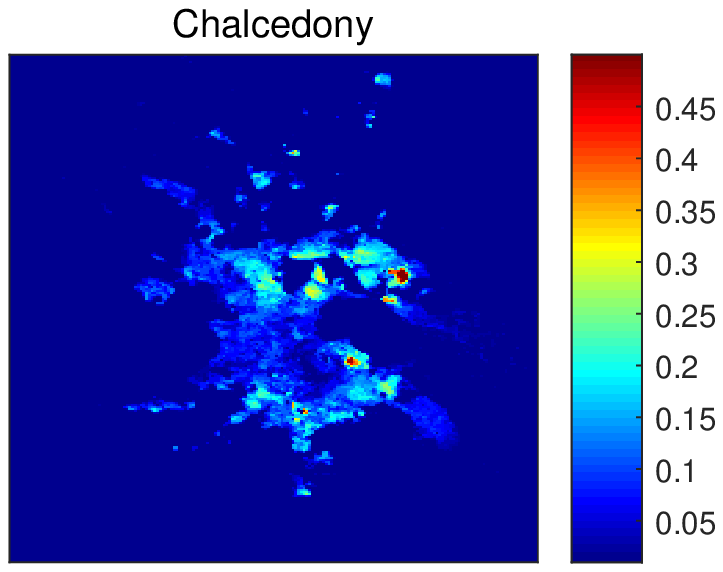}}
}\\(a)

\mbox{
{\includegraphics[width=1.5in]{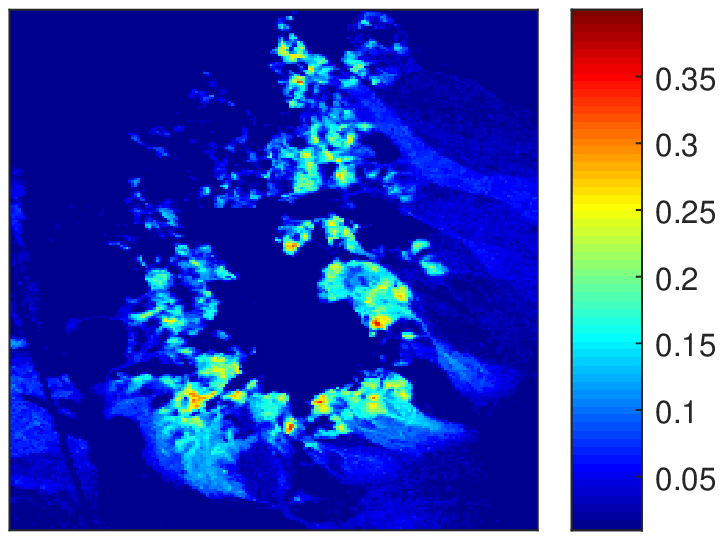}}
{\includegraphics[width=1.5in]{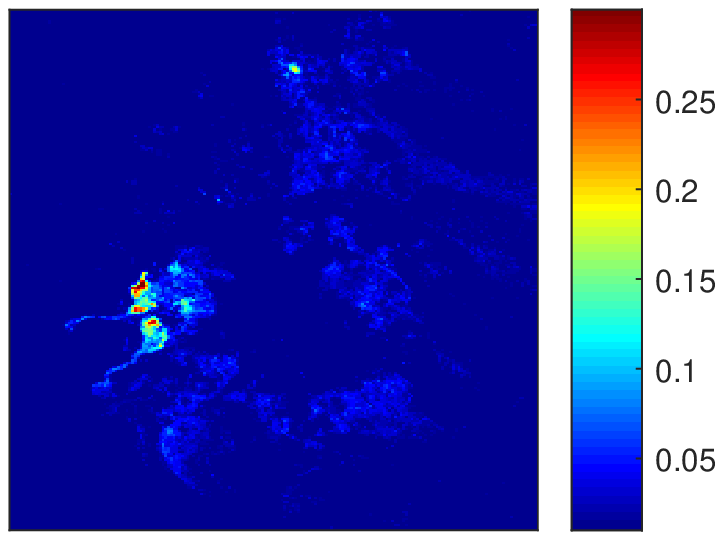}}
{\includegraphics[width=1.5in]{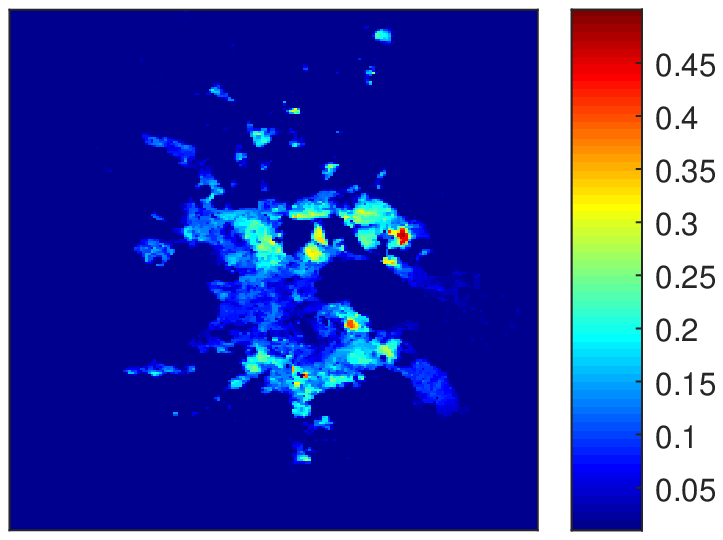}}
}\\(b)

\mbox{
{\includegraphics[width=1.5in]{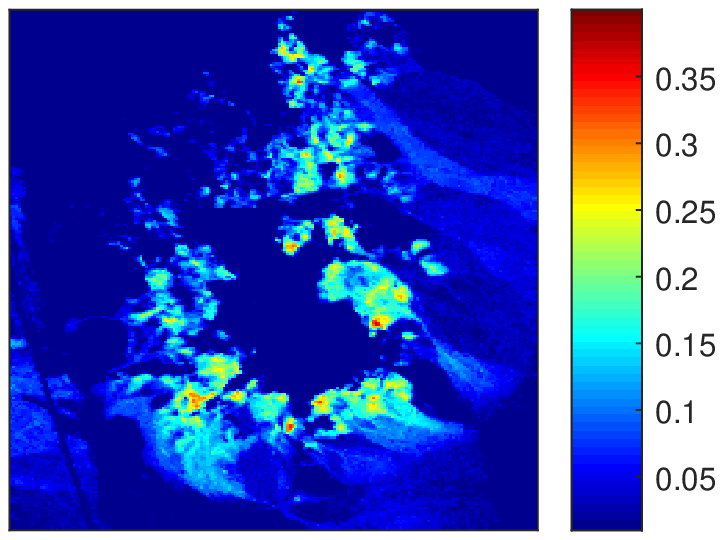}}
{\includegraphics[width=1.5in]{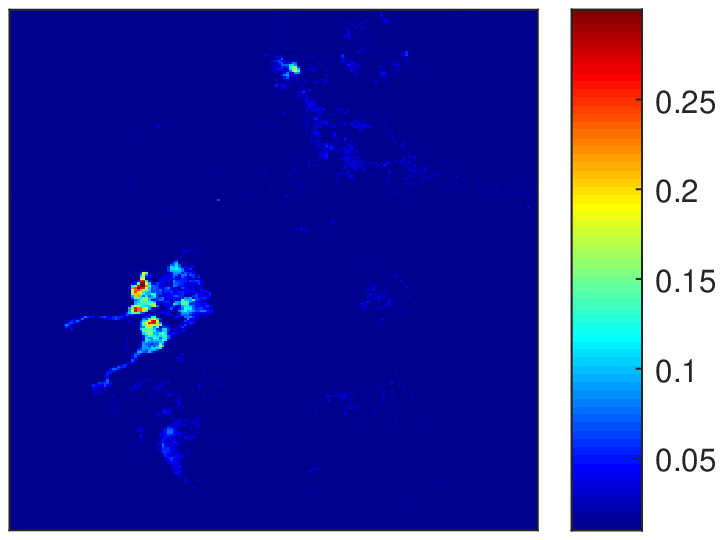}}
{\includegraphics[width=1.5in]{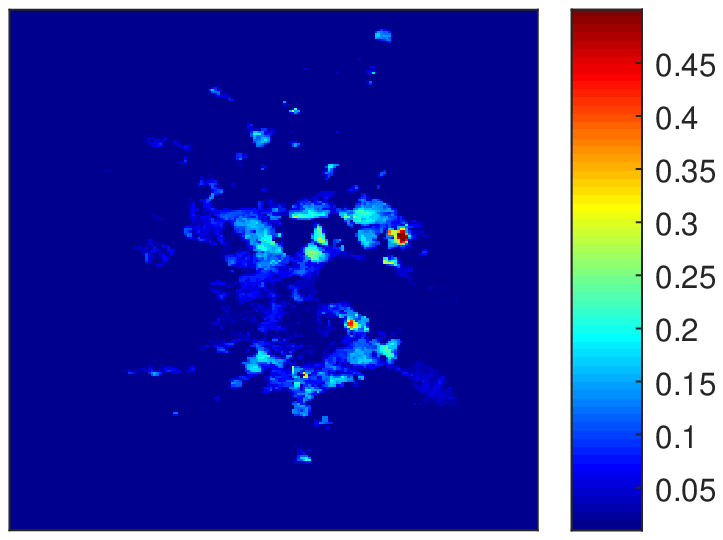}}
}\\(c)

\mbox{
{\includegraphics[width=1.5in]{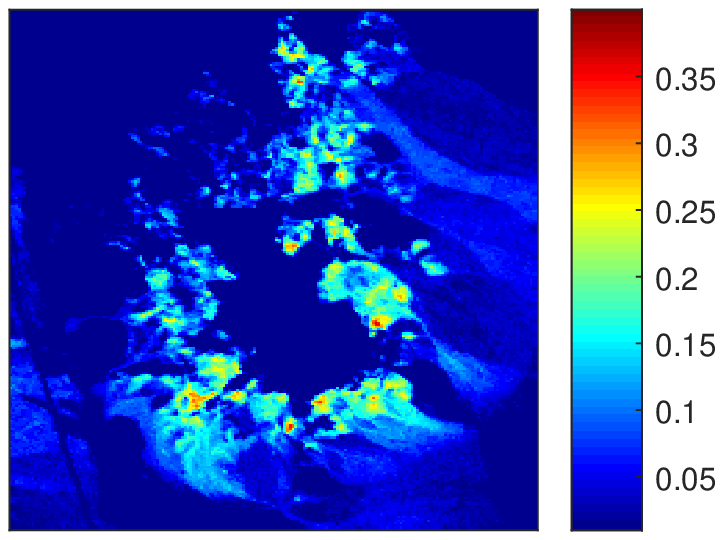}}
{\includegraphics[width=1.5in]{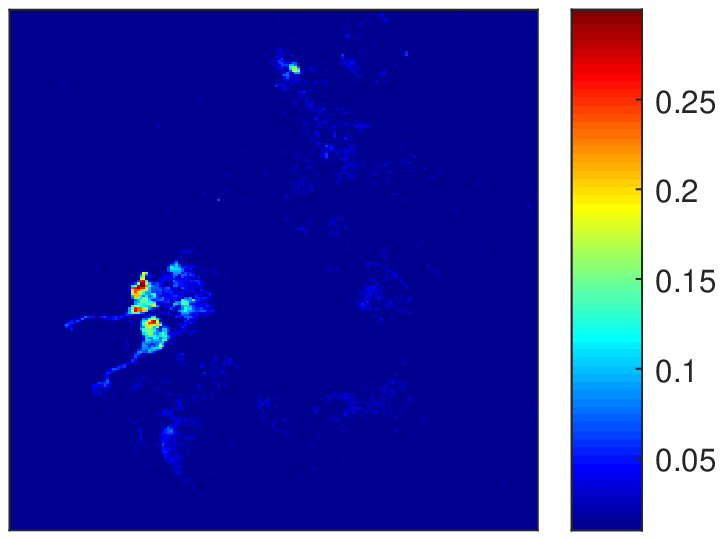}}
{\includegraphics[width=1.5in]{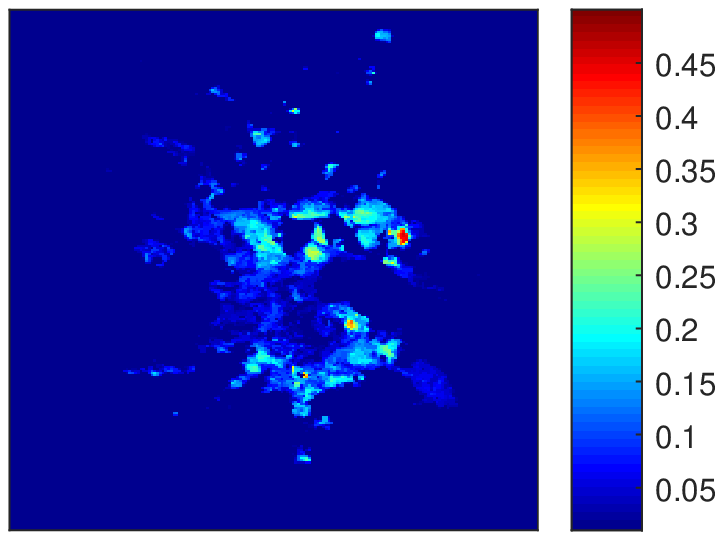}}
}\\(d)

\caption{The estimated abundances obtained by the different unmixing algorithms for the minerals: alunite, buddingtonite, chalcedony. From top to bottom: (a) SUnSAL-TV ($\lambda=0.001,\lambda_{TV}=0.0001$); (b) SUnTV-sGSADMM ($\lambda=0.001,\lambda_{TV}=0.0001$); (c) CLSUnSAL-TV ($\lambda=0.001,\lambda_{TV}=0.0001$); (d) CLSUnTV-sGSADMM ($\lambda=0.001,\lambda_{TV}=0.0001$).}
\label{RealDataImages}
\end{figure*}

\begin{figure*}[!h]
\centering
\mbox{
{\includegraphics[width=2.0in,height=1.5in]{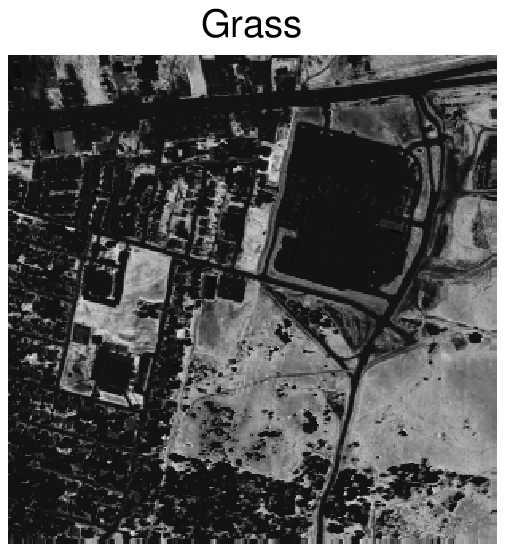}}
{\includegraphics[width=2.0in,height=1.5in]{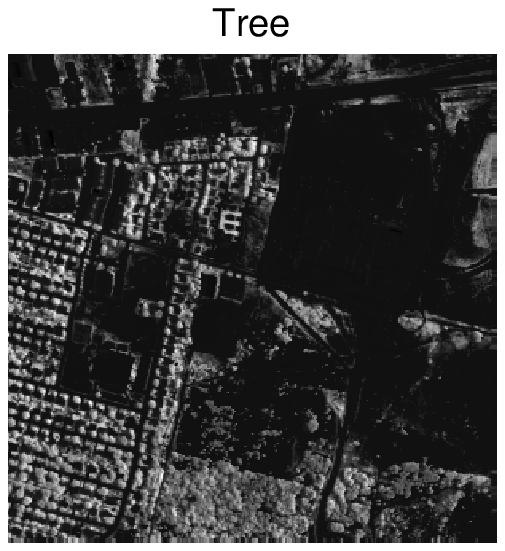}}
{\includegraphics[width=1.4in,height=1.5in]{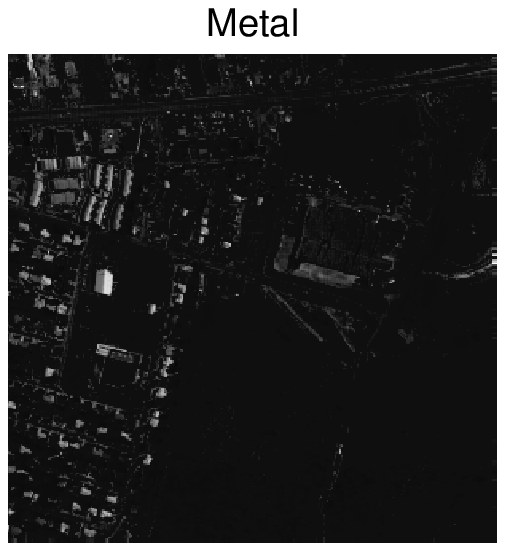}}
}\\(a)

\mbox{
{\includegraphics[width=2.0in,height=1.5in]{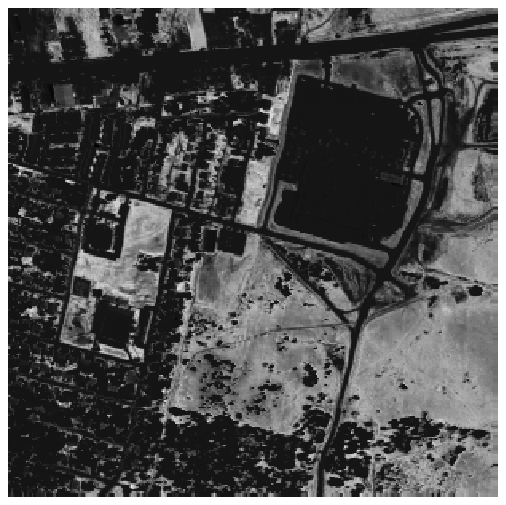}}
{\includegraphics[width=2.0in,height=1.5in]{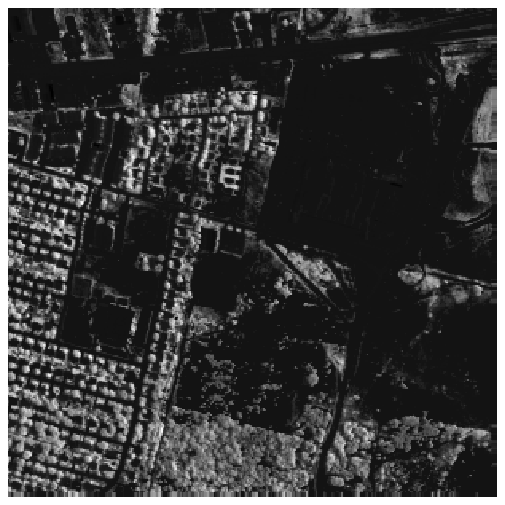}}
{\includegraphics[width=1.4in,height=1.5in]{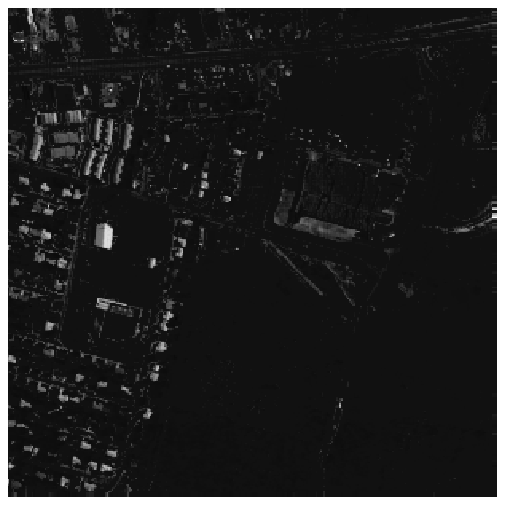}}
}\\(b)

\mbox{
{\includegraphics[width=2.0in,height=1.5in]{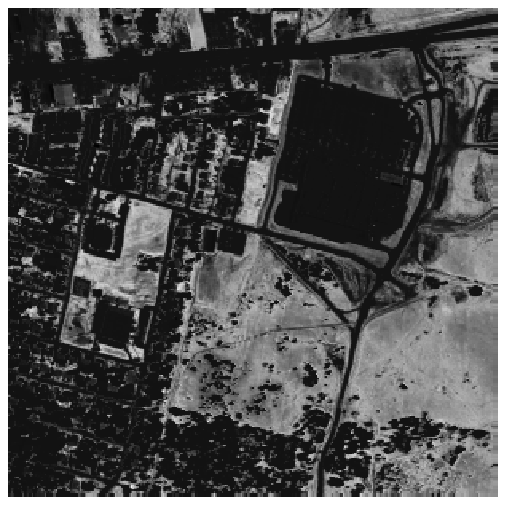}}
{\includegraphics[width=2.0in,height=1.5in]{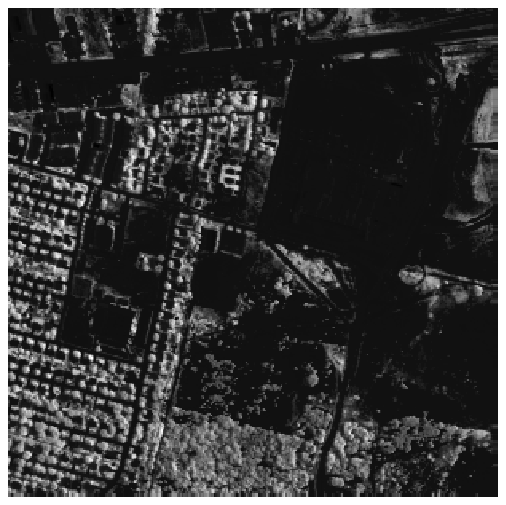}}
{\includegraphics[width=1.4in,height=1.5in]{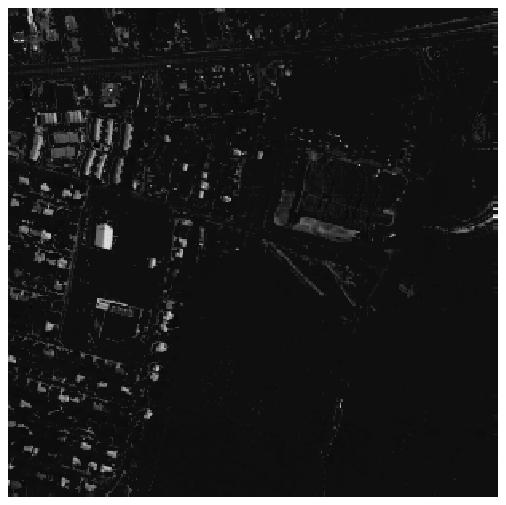}}
}\\(c)

\mbox{
{\includegraphics[width=2.0in,height=1.5in]{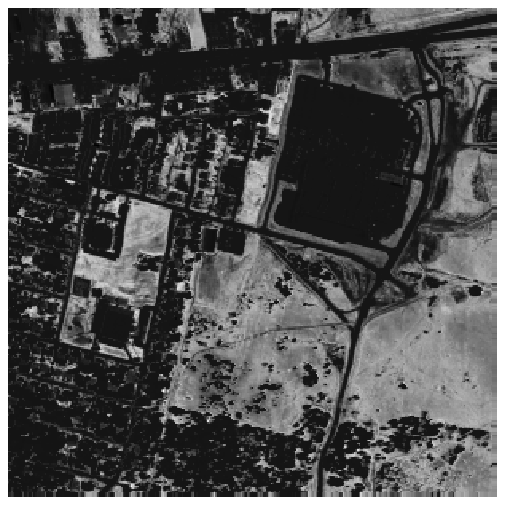}}
{\includegraphics[width=2.0in,height=1.5in]{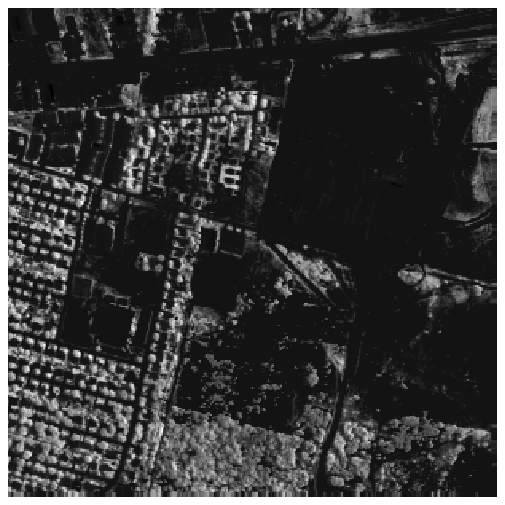}}
{\includegraphics[width=1.4in,height=1.5in]{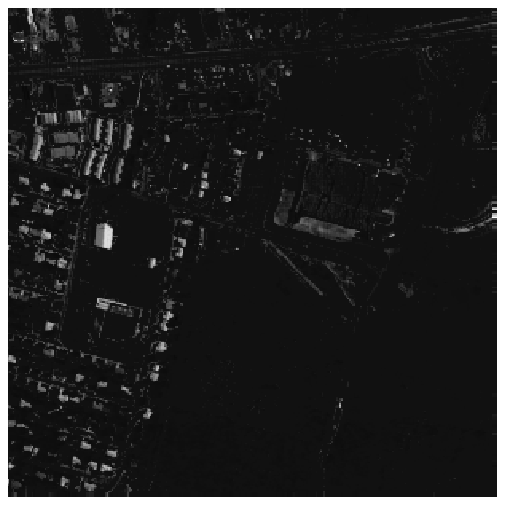}}
}\\(d)

\caption{The estimated abundances obtained by the different unmixing algorithms for the endmembers: grass, tree, metal. From top to bottom: (a) SUnSAL-TV ($\lambda=0.0001,\lambda_{TV}=0.0001$), (b) SUnTV-sGSADMM ($\lambda=0.0001,\lambda_{TV}=0.0001$), (c) CLSUnSAL-TV ($\lambda=0.0001,\lambda_{TV}=0.0001$), (d) CLSUnTV-sGSADMM ($\lambda=0.0001,\lambda_{TV}=0.0001$).}
\label{UrbanDataImages}
\end{figure*}

The first real hyperspectral remote sensing data is from a very famous Airborne Visible Infrared Imaging Spectrometer (AVIRIS) Cuprite data set.\footnote{Available online: http://aviris.jpl.nasa.gov/html/aviris.freedata.html.}
This data set has been widely used to verify the performances of the unmixing algorithms since it is well understood in mineralogy field and also has several exposed minerals of interest \cite{CLSUnSAL}.
The data used in our experiments corresponds to a $250\times191$-pixel subset of the sector labeled as f970619t01p02\_r02\_{s}c03.a.rfl in the online data.
The scene comprises 224 spectral bands between 0.4 and 2.5 $\mu$m, with nominal spectral resolution of 10 nm.
Prior to the analysis, bands 1-2, 105-115, 150-170, 223-224 were cut off because of water absorption and low SNR in those bands, leaving a total of 188 spectral bands.
The spectral library used here is to select 498 spectra from the USGS spectral library and remove the corresponding bands.
The computing time of the different unmixing algorithms for the AVIRIS Cuprite data are given in Table 3. As one can see from Table 3, our proposed algorithm is about 10 times faster than the primal ADMM.
To further illustrate the convergence, in Fig. \ref{ConvergenceAnalysisRealData}, we plot the variations of  $\min(\mathbf{R}_P,\mathbf{R}_D)$ with respect to the iterations for different algorithms when dealing with the AVIRIS Cuprite data.  As one can see that our algorithm converges faster than the primal ADMM.
We can just make a visual comparison on the abundance maps of the minerals since the true abundances of the AVIRIS Cuprite data are unknown.
Fig. \ref{RealDataImages} shows the estimated abundances obtained by the different unmixing algorithms for the minerals: alunite, buddingtonite, chalcedony.
From Fig. \ref{RealDataImages}, it can be observed that the effects by the SUnTV-sGSADMM and the CLSUnTV-sGSADMM are as good as those by the SUnSAL-TV and the CLSUnSAL-TV.

The second real hyperspectral remote sensing data is the Urban data captured by the Hyperspectral Digital Imagery Collection Experiment (HYDICE) sensor over an area located as Copperas Cove near Fort Hood, TX, U.S., in October 1995.
There are 307$\times$307 pixels. The scene comprises 210 spectral bands ranging from 0.4 $\mu$m  to 2.5 $\mu$m, with nominal spectral resolution of 10 nm. After removing the bands 1-4, 76, 87, 101-111, 136-153, 198-210, we remain 162 bands.
The spectral library used here is the endmembers obtained by the method provided in \cite{Jia2007}-\cite{Zhu2009}.
The computing time of the different unmixing algorithms for the Urban data are given in Table 4.
As one can see, our proposed algorithm is about 4 times faster than the primal ADMM.
Fig. 6 shows the estimated abundances obtained by the different unmixing algorithms for the endmembers: grass, tree, metal.
Again, the proposed unmixing algorithms show similar performances when compared to the remaining unmixing algorithms.

\section{Conclusion}
In this paper, we developed an efficient and convergent dual sGS-ADMM for the hyperspectral sparse unmixing with a TV regularization term. As shown in the numerical experiments, this approach can obviously improve the efficiency of the unmixing compared with the state-of-the-art algorithm. More importantly, we can obtain relatively higher SREs for different problems.
Our future work will focus on how to choose the regularization parameters adaptively for our algorithm and how to generalize our algorithm to deal with the nonlinear models.

\section*{Acknowledgements}
We would like to thank Professor Xile Zhao at School of Mathematics, University of Electronic Science and Technology of China for his useful comments and suggestions.
We also thank Professor Heng-Chao Li at School of Information Science and Technology, Southwest Jiaotong University for fruitful discussions.
Besides, We are grateful to the two anonymous referees and the Editor-in-Chief Prof. Claude Brezinski for their constructive and helpful suggestions on improving the quality of the paper.



\begin{thebibliography}{99}
     \bibitem{Bioucas2012} Bioucas-Dias, J.M., Plaza, A., Dobigeon, N., Parente, M., Du, Q., Gader, P., Chanussot, J.: Hyperspectral unmixing overview: Geometrical, statistical, and sparse regression-based approaches, IEEE J. Sel. Topics Appl. Earth Observ. Remote Sens. \textbf{5}(2), 354-379 (2012)

     \bibitem{VCA} Nascimento, J.P., Bioucas-Dias, J.M.: Vertex component analysis: A fast algorithm to unmix hyperspectral data, IEEE Trans. Geosci. Remote Sens. \textbf{43}(4), 898-910 (2005)

     \bibitem{PPI} Boardman, J.W., Kruse, F.A., Green, R.O.: Mapping target signatures via partial unmixing of AVIRIS data, in Proc. JPL Airborne Earth Sci. Workshop. 23-26 (1995)

     \bibitem{SGA} Chang, C.-I., Wu, C.-C., Liu, W., Ouyang, Y.-C.: A new growing method for simplex-based endmember extraction algorithm, IEEE Trans. Geosci. Remote Sens. \textbf{44}(10), 2804-2819 (2006)

     \bibitem{MVES} Chan, T.-H., Chi, C.-Y., Huang, Y.-M., Ma, W.-K.: Convex analysis based minimum-volume enclosing simplex algorithm for hyperspectral unmixing, IEEE Trans. Signal Process. \textbf{57}(11), 4418-4432 (2009)

     \bibitem{ICE} Berman, M., Kiiveri, H., Lagerstrom, R., Ernst, A., Dunne, R., Huntington, J.F.: ICE: A statistical approach to identifying endmembers in hyperspectral images, IEEE Trans. Geosci. Remote Sens. \textbf{42}(10), 2085-2095 (2004)

     \bibitem{MVC-NMF} Miao, L., Qi, H.: Endmember extraction from highly mixed data using minimum volume constrained nonnegative matrix factorization, IEEE Trans. Geosci. Remote Sens. \textbf{45}(3), 765-777 (2007)

     \bibitem{Dobigeon2009} Dobigeon, N., Moussaoui, S., Coulon, M., Tourneret, J.-Y., Hero, A.O.: Joint Bayesian endmember extraction and linear unmixing for hyperspectral imagery, IEEE Trans. Signal Process. \textbf{57}(11), 4355-4368 (2009)

     \bibitem{Bruckstein2009} Bruckstein, A.M., Donoho, D.L., Elad, M.: From sparse solutions of systems of equations to sparse modeling of signals and images, SIAM Rev. \textbf{51}(1), 34-81 (2009)

     \bibitem{MUSIC-CSR} Iordache, M.-D., Bioucas-Dias, J.M., Plaza, A., Somers, B.: MUSIC-CSR: Hyperspectral unmixing via multiple signal classification and collaborative sparse regression, IEEE Trans. Geosci. Remote Sens. \textbf{52}(7), 4364-4382 (2014)

     \bibitem{RIC} Cand\`{e}s, E.J., Romberg, J., Tao, T.: Robust uncertainty principles: Exact signal reconstruction from highly incomplete frequency information, IEEE Trans. Inf. Theory. \textbf{52}(2), 489-509 (2006)

     \bibitem{SUnSAL} Iordache, M.-D., Bioucas-Dias, J.M., Plaza, A.: Sparse unmixing of hyperspectral data, IEEE Trans. Geosci. Remote Sens. \textbf{49}(6), 2014-2039 (2011)

     \bibitem{Rudin1992} Rudin, L.I., Osher, S., Fatemi, E.: Nonlinear total variation based noise removal algorithms, Phys. D Nonlinear Phenom. \textbf{60}(1-4), 259-268 (1992)

     \bibitem{Zhao2013} Zhao, X.-L., Wang, W., Zeng, T.-Y., Huang, T.-Z., Ng, M.K.: Total variation structured total least squares method for image restoration, SIAM J. Sci. Comput. \textbf{35}(6), 1304-1320 (2013)

     \bibitem{Zakharova2017} Zakharova, A.: Total variation reconstruction from quadratic measurements, Numer. Algorith. \textbf{75}(1), 81-92 (2017)

      \bibitem{SUnSAL-TV} Iordache, M.-D., Bioucas-Dias, J.M., Plaza, A.: Total variation spatial regularization for sparse hyperspectral unmixing, IEEE Trans. Geosci. Remote Sens. \textbf{50}(11), 4484-4502 (2012)

     \bibitem{Zhang17} Zhang, S., Li, J., Liu, K., Plaza, A.: Hyperspectral unmixing based on local collaborative sparse regression, IEEE Geosci. Remote Sens. Lett. \textbf{13}(5), 631-635 (2016)

     \bibitem{Zhang16} Zhang, L., Wei, W., Zhang, Y., Yan, H., Li, F., Tian, C.: Locally similar sparsity-based hyperspectral compressive sensing using unmixing, IEEE Trans. Comput. Imag. \textbf{2}(2), 86-100 (2016)

     \bibitem{CLSUnSAL} Iordache, M.-D., Bioucas-Dias, J.M., Plaza, A.: Collaborative sparse regression for hyperspectral unmixing, IEEE Trans. Geosci. Remote Sens. \textbf{52}(1), 341-354 (2014)

     \bibitem{CLSUnSAL-TV} Chen, Y.-J., Ge, W.-D., Sun, L.: A novel linear hyperspectral unmixing method based on collaborative sparsity and total variation, Acta Auto. Sinica. \textbf{44}(1), 116-128 (2018)

     \bibitem{GlowinskiM} Glowinski, R., Marroco, A.: Sur l'approximation, par \'{e}l\'{e}ments finis d'ordre un, et la r\'{e}solution, par p\'{e}nalisation-dualit\'{e} d'une classe de probl\'{e}mes de Dirichlet non lin\'{e}aires, Revue Francaise d'Automatique, Informatique et Recherche Op\'{e}rationelle \textbf{2}(R-2), 41-76 (1975)

     \bibitem{GabayM} Gabay, D., Mercier, B.: A dual algorithm for the solution of nonlinear variational problems via finite element approximation, Comput. Math. Appl. \textbf{2}(1), 17-40 (1976)

     \bibitem{EcksteinY} Eckstein, J., Yao, W.: Understanding the convergence of the alternating direction method of multipliers: Theoretical and computational perspectives, Pac. J. Optim. \textbf{11}(4), 619-644 (2015)

     \bibitem{Glowinski} Glowinski, R.: On alternating direction methods of multipliers: A historical perspective, in W. Fitzgibbon, Y. A. Kuznetsov, P. Neittaanmaki and O. Pironneau (eds.), Modeling, Simulation and Optimization for Science and Technology 59-82 (2014)

     \bibitem{FazelPST} Fazel, M., Pong, T.K., Sun, D., Tseng,  P.: Hankel matrix rank minimization with applications to system identification and realization, SIAM J. Matrix Anal. Appl. \textbf{34}(3), 946-977 (2013)

     \bibitem{linear-rate} Han,  D., Sun, D., Zhang, L.: Linear rate convergence of the alternating direction method of multipliers for convex composite programming. (2015). arXiv:1508.02134

     \bibitem{CHYY} Chen, C., He, B., Ye, Y., Yuan, X.: The direct extension of ADMM for multi-block convex minimization problems is not necessarily convergent, Math. Program. \textbf{155}(1-2), 57-79 (2016)

     \bibitem{LiST2016} Li, X., Sun, D., Toh, K.-C.: A Schur complement based semi-proximal ADMM for convex quadratic conic programming and extensions, Math. Program. \textbf{155}(1-2), 333-373 (2016)

     \bibitem{LiST2019} Li, X., Sun, D., Toh, K.-C.: A block symmetric Gauss-Seidel decomposition theorem for convex composite quadratic programming and its applications, Math. Program. \textbf{175}(1-2), 395-418 (2019)

     \bibitem{Chen2017} Chen, L., Sun, D., Toh, K.-C.: An efficient inexact symmetric Gauss-Seidel based majorized ADMM for high-dimensional convex composite conic programming, Math. Program. \textbf{161}(1-2), 237-270 (2017)

     \bibitem{Chen2018} Chen, L., Sun, D., Toh, K.-C., Zhang, N.: A unified algorithmic framework of symmetric Gauss-Seidel decomposition based proximal ADMMs for convex composite programming, Math. Program. (2018). arXiv:1812.06579

     \bibitem{Beck2009} Beck, A., Teboulle, M.: Fast gradient-based algorithms for constrained total variation image denoising and deblurring problems, IEEE Trans. Image Process. \textbf{18}(11), 2419-2434 (2009)

     \bibitem{Zuo2011} Zuo, W., Lin, Z.: A generalized accelerated proximal gradient approach for total-variation-based image restoration, IEEE Trans. Image Process. \textbf{20}(10), 2748-2759 (2011)

     \bibitem{Rockafellar1970} Rockafellar, R.T.: Convex Analysis, University of Princeton, USA (1970)

     \bibitem{Horn1985} Horn, R.A., Johnson, C.R.: Matrix Analysis, University of Cambridge, UK (1985)

     \bibitem{Sun1986} Sun, J.: On monotropic piecewise quadratic programming, Ph.D. dissertation, University of Washington, USA (1986)

     \bibitem{Robinson1981} Robinson, S.M.: Some contimuity properties of polyhedral multifuntions, Math. Program. Study \textbf{14}, 206-214 (1981)

     \bibitem{Rockafellar2009} Dontchev, A., Rockafellar, R.T.: Implicit Function and Solution Mappings, Springer, USA (2009)

     \bibitem{Rockafellar1998} Rockafellar, R.T., Wets, R.J.-B.: Variational Analysis, Springer, USA (1998)

    \bibitem{LiHengChao2017} Wang, R., Li, H.-C., Liao, W., Huang, X., Philips, W.: Centralized collaborative sparse unmixing for hyperspectral images, IEEE J. Sel. Topics Appl. Earth Observ. Remote Sens. \textbf{10}(5), 1949-1962 (2017)

     \bibitem{Condat2013} Condat, L.: A direct algorithm for 1-D total variation denoising, IEEE Signal Process. Lett. \textbf{20}(11), 1054-1057 (2013)

     \bibitem{Friedman2007} Friedman, J., Hastie, T., H\"{o}fling, H., Tibshirani, R.: Pathwise coordinate optimization, Ann. Appl. Stat. \textbf{1}(2), 302-332 (2007)

     \bibitem{Mattern1999} Kozintsev, B.: Computations with Gaussian random fields, Ph.D. dissertation, University of Maryland, USA (1999)

     \bibitem{Yu2013} Yu, Y.: On decomposing the proximal map, Adv. Neural Inform. Process. Syst \textbf{1}, 91-99 (2013)

     \bibitem{Yin2008} Yin, W.T., Osher, S., Goldfarb, D., Darbon, J.: Bregman iterative algorithms for $l_1$-minimization with applications to compressed sensing, SIAM J. Imag. Sci. \textbf{1}(1), 143-168 (2008)

     \bibitem{Wright2009} Wright, S.J., Nowak, R.D., Figueiredo, M.A.T.: Sparse reconstruction by separable approximation, IEEE Trans. Signal Process. \textbf{57}(7), 2479-2493 (2009)

     \bibitem{Jia2007} Jia, S., Qian, Y.: Spectral and spatial complexity-based hyperspectral unmixing, IEEE Trans. Geosci. Remote Sens. \textbf{45}(12), 3867-3879 (2007)

     \bibitem{Jia2009} Jia, S., Qian, Y.: Constrained nonnegative matrix factorization for hyperspectral unmixing, IEEE Trans. Geosci. Remote Sens. \textbf{47}(1), 161-173 (2009)

     \bibitem{Zhu2009} Zhu, F., Wang, Y., Xiang, S., Fan, B., Pan, C.: Structured sparse method for hyperspectral unmixing, ISPRS J. Photogramm. Remote Sens. \textbf{88}(2), 101-118 (2014)

\end{thebibliography}
\end{document}